\documentclass[a4paper, 10.5pt]{amsart}

\usepackage[left=2.5cm,right=2.5cm,top=3.5cm,bottom=3cm]{geometry}
\numberwithin{equation}{subsection}

\usepackage{tikz-cd}
\usepackage{graphicx}				
\usepackage{amsthm,amssymb,amsmath,amsfonts,mathrsfs,amscd}
\usepackage{relsize}
\usepackage{stmaryrd}
\usepackage{enumitem}
\usepackage{color}
\usepackage{hyperref}
\usepackage{tikz-cd}
\usepackage{url}
\usepackage{soul}

\theoremstyle{plain}
\newtheorem{thm}{Theorem}[subsection]

\newtheorem*{theo}{Theorem}
\newtheorem{theor}{Theorem}

\newtheorem*{con}{Conjecture}
\newtheorem{conj}[thm]{Conjecture}
\newtheorem{lemma}[thm]{Lemma}
\newtheorem{prop}[thm]{Proposition}
\newtheorem{corollary}[thm]{Corollary}
\theoremstyle{definition}
\newtheorem{definition}[thm]{Definition}
\theoremstyle{remark}
\newtheorem{rmk}[thm]{Remark}

\def\Z{\mathbb{Z}}
\def\N{\mathbb{N}}
\def\Q{\mathbb{Q}}
\def\Qp{\mathbb{Q}_p}
\def\Zp{\mathbb{Z}_p}

\def\bQp{\overline{\mathbb{Q}}_p}

\def\bQ{\overline{\mathbb{Q}}}

\def\R{\mathbb{R}}
\def\C{\mathbb{C}}
\def\A{\mathbb{A}}
\def\F{\mathbb{F}}

\def\O{\mathcal{O}}
\def\m{\mathfrak{m}}

\def\B{\mathfrak{B}}
\def\V{\mathbf{V}}
\def\cV{\check{\mathbf{V}}}
\def\H{\mathcal{H}}
\def\D{\mathcal{D}}

\def\p{\mathfrak{p}}
\def\U{\mathbf{U}}

\def\cS{\mathcal{S}}
\def\cT{\mathcal{T}}
\def\cC{\mathcal{C}}
\def\L{\mathcal{L}}
\def\P{\mathbb{P}}
\def\cR{\mathcal{R}}

\def\cU{\mathbb{U}}
\def\cX{\mathfrak{X}}
\def\W{\mathcal{W}}
\def\cR{\mathcal{R}}
\def\I{\mathcal{I}}
\def\iT{\mathit{T}}
\def\Rone{\mathbb{R}^1}
\def\fC{\mathfrak{C}}
\def\fT{\mathfrak{T}}
\def\fD{\mathfrak{D}}

\def\GL{\operatorname{GL}}
\def\SL{\operatorname{SL}}
\def\PGL{\operatorname{PGL}}
\def\Gal{\operatorname{Gal}}
\def\GQp{G_{\Qp}}

\def\Sp{\operatorname{Sp}}
\def\cInd{\operatorname{c-Ind}}
\def\Ind{\operatorname{Ind}}

\def\ind{\operatorname{ind}}
\def\1{\mathbf{1}}
\def\bracketGQp{\llbracket \GQp \rrbracket}

\def\pr{\operatorname{pr}}

\def\Mod{\operatorname{Mod}}
\def\sm{\operatorname{sm}}
\def\Modsm{\operatorname{Mod}^{\operatorname{sm}}_G}
\def\Modlfin{\operatorname{Mod}^{\operatorname{l.fin}}_G}
\def\Modladm{\operatorname{Mod}^{\operatorname{l.adm}}_G}
\def\Modsmz{\operatorname{Mod}^{\operatorname{sm}}_{G, \zeta}}
\def\Modlfinz{\operatorname{Mod}^{\operatorname{l.fin}}_{G, \zeta}}
\def\Modladmz{\operatorname{Mod}^{\operatorname{l.adm}}_{G, \zeta}}
\def\Modpro{\operatorname{Mod}^{\operatorname{pro}}}
\def\Ban{\operatorname{Ban}^{\operatorname{adm}}_{G, \zeta}}
\def\Irr{\operatorname{Irr}_{G, \zeta}}
\def\Modfinz{\operatorname{Mod}^{\operatorname{fin}}_{G, Z}}
\def\ModGQp{\operatorname{Mod}^{\operatorname{fin}}_{\GQp}}

\def\CNL{\operatorname{CNL}}
\def\Set{\operatorname{Sets}}

\def\ad{\operatorname{ad}}
\def\WD{\operatorname{WD}}
\def\HT{\operatorname{HT}}

\def\rec{\operatorname{rec}}

\def\Art{\operatorname{Art}}

\def\ab{\operatorname{ab}}

\def\muG{\mu_{\operatorname{Gal}}}
\def\muA{\mu_{\operatorname{Aut}}}

\def\rhobar{\overline{\rho}}
\def\rhof{\rho^{\Box}}
\def\Df{\mathcal{D}^{\Box}}
\def\Dz{\mathcal{D}^{\psi}}
\def\Dfz{\mathcal{D}^{\Box, \psi}}
\def\DTW{\mathcal{D}^{\text{TW}}}

\def\Rf{R^{\Box}}
\def\Rfz{R^{\Box, \psi}}
\def\Rz{R^{\psi}}

\def\RSt{R^{St}}
\def\DSt{\D^{St}}
\def\Rord{R^{\Delta}}
\def\Dord{\D^{\Delta}}
\def\Rss{R^{\lambda_v, \tau_v}}
\def\Rcr{R^{\lambda_v, \tau_v, cr}}
\def\Dss{\D^{\lambda_v, \tau_v, ss}}
\def\Dcr{\D^{\lambda_v, \tau_v, cr}}
\def\Rssord{R^{\Delta, \lambda_v, \tau_v}}

\def\Rodd{R^{odd}}
\def\Dodd{\D^{odd}}
\def\Rur{R^{ur}}
\def\Dur{\D^{ur}}
\def\ver{\operatorname{ver}}
\def\ps{\operatorname{ps}}
\def\tRfz{\tilde{R}^{\Box, \psi}}
\def\univ{\operatorname{univ}}

\def\Af{\A^{\infty}_F}
\def\Aft{(\A^{\infty}_F)^{\times}}
\def\T{\mathbb{T}}

\def\Norm{\mathbf{N}}
\def\Iw{\operatorname{Iw}}
\def\Mord{M^{\ord}_{\psi}}
\def\Sord{S^{\ord}_{\psi}}

\def\ord{\operatorname{ord}}

\def\inv{\operatorname{inv}}

\def\Rloc{R^{\loc}}
\def\a{\mathfrak{a}}
\def\Rinf{R_{\infty}}
\def\Oinf{\O_{\infty}}
\def\Minf{M_{\infty}}
\def\Deltainf{\Delta_{\infty}}
\def\Rinvinf{R^{\inv}_{\infty}}
\def\Piinf{\Pi_{\infty}}
\def\Gm{\mathbb{G}_m}
\def\Gmc{\hat{\mathbb{G}}_m}
\def\filter{\mathfrak{F}}

\def\Roinf{R^{\Delta}_{\infty}}
\def\Roinfp{R^{\Delta, \prime}_{\infty}}
\def\Moinf{M^{\Delta}_{\infty}}
\def\Roinvinf{R^{\Delta, \inv}_{\infty}}

\def\Sinf{S_{\infty}}
\def\Moinf{M^{\ord}_{\infty}}
\def\Ann{\operatorname{Ann}}
\def\red{\operatorname{red}}

\def\LL{\operatorname{LL}}
\def\bsigma{\overline{\sigma}}
\def\rbar{\overline{r}}

\def\pialg{\pi_{\operatorname{alg}}}
\def\pism{\pi_{\operatorname{sm}}}
\def\lalg{\operatorname{l.alg}}

\def\Supp{\operatorname{Supp}}

\def\Rzp{R^{\zeta}_{\rbar}}

\def\tMinf{\tilde{M}_{\infty}}
\def\tRinf{\tilde{R}_{\infty}}
\def\tRinvinf{\tilde{R}^{\inv}_{\infty}}
\def\tRinvinfp{\tilde{R}^{\inv, \prime}_{\infty}}
\def\rinf{r_{\infty}}

\def\loc{\operatorname{loc}}
\def\tPiinf{\tilde{\Pi}_{\infty}}

\def\tAinf{\tilde{A}_{\infty}}
\def\tAinvinf{\tilde{A}^{\inv}_{\infty}}
\def\tAinvinfp{\tilde{A}^{\inv, \prime}_{\infty}}
\def\tNinf{\tilde{N}_{\infty}}

\def\rig{\operatorname{rig}}

\def\hT{\hat{T}}
\def\unr{\operatorname{unr}}

\def\tri{\operatorname{tri}}

\def\Xtri{X^{\tri}}
\def\Xinf{\cX_{\infty}}
\def\an{\operatorname{an}}
\def\Xtrinf{X^{\tri}_{\infty}}
\def\sign{\operatorname{sign}}
\def\Rfsign{R^{\Box, \sign}}

\def\Sym{\operatorname{Sym}}
\def\det{\operatorname{det}}
\def\tr{\operatorname{tr}}
\def\Hom{\operatorname{Hom}}
\def\Homc{\operatorname{Hom}^{\operatorname{cont}}}
\def\Ext{\operatorname{Ext}}
\def\End{\operatorname{End}}
\def\Ker{\operatorname{ker}}

\def\ctimes{\hat{\otimes}}
\def\Spec{\operatorname{Spec}}
\def\mSpec{\operatorname{m-Spec}}
\def\Spf{\operatorname{Spf}}

\def\Frob{\operatorname{Frob}}
\def\sm{\operatorname{sm}}

\title{On the modularity of 2-adic potentially semi-stable deformation rings}
\author{Shen-Ning Tung}

\address{Fakult\"at f\"ur Mathemati,\\
  Universit\"at Duisburg-Essen\\
  45127 Essen, Germany}
\email{shen-ning.tung@stud.uni-due.de}
				
\begin{document}

\begin{abstract}
Using $p$-adic local Langlands correspondence for $\GL_2(\mathbb{Q}_2)$ and an ordinary $R = \T$ theorem, we prove that the support of patched modules for quaternionic forms meet every irreducible component of the potentially semi-stable deformation ring. This gives a new proof of the Breuil-M\'{e}zard conjecture for 2-dimensional representations of the absolute Galois group of $\mathbb{Q}_2$, which is new in the case $\rbar$ a twist of an extension of the trivial character by itself. As a consequence, a local restriction in the proof of Fontaine-Mazur conjecture in \cite{MR3544298} is removed.
\end{abstract}

\maketitle

\section*{Introduction}
Let $p$ be a prime number and $\O$ be the ring of integers of a sufficiently large finite extension over $\Qp$. Let $f$ be a normalized cuspidal eigenform of weight $k \geq 2$ and level $N \geq 1$, normalized so that $f$ has Fourier expansion $f = \sum_1^{\infty} a_n q^n$, with $a_1 = 1$. It is proved that there exists a Galois representation
\[
\rho_f : \Gal(\bQ / \Q) \rightarrow \GL_2(\O)
\]
by Eichler and Shimura for $k=2$, and Deligne for $k \geq 2$, characterized by the following property: $\rho_f$ is unramified at primes $l \nmid pN$ with $\tr(\rho_f(\Frob_l)) = a_l$. Due to the work of many people, the representation is known to be irreducible, odd (i.e. $\det \rho_f(c) = -1$ with $c$ the complex conjugation), and de Rham (in the sense of Fontaine) at $p$ with Hodge-Tate weights $(0, k-1)$.

In \cite{MR1363495} Fontaine and Mazur made a conjecture which asserts the converse: 

\begin{con}[Fontaine-Mazur] \label{cong:FM}
Let
\[
\rho: \Gal(\bQ / \Q) \rightarrow \GL_2(\O)
\]
be a continuous, irreducible representation such that
\begin{itemize}
\item $\rho$ is odd;
\item $\rho$ is unramified outside all but finitely many places;
\item the restriction of $\rho$ at the decomposition group at $p$ is de Rham with distinct Hodge-Tate weights.
\end{itemize}
Then (up to a twist) $\rho \cong \rho_f$ for some cuspidal eigenform $f$.
\end{con}

We will say that $\rho$ is modular if it is isomorphic to a twist of $\rho_f$ by a character. Similarly, we will say that $\rhobar: \Gal(\bQ/ \Q) \rightarrow \GL_2(k)$ is modular if $\rhobar \cong \rhobar_f$ up to a twist, where $k$ is the residue field of $\O$ and $\rhobar$ is obtained by reducing the matrix entries of $\rho_f$ modulo the maximal ideal of $\O$. This conjecture has been proved in several cases under different assumptions, e.g. \cite{MR2251474, Emerton2006localglobal}. We will only focus on those related to the groundbreaking work of Kisin in \cite{MR2505297}.

\begin{theo}[Kisin, Pa\v{s}k\=unas, Hu-Tan, Tung]
Let $\rho$ be as in the conjecture. Let $\rhobar: \Gal(\bQ / \Q) \rightarrow \GL_2(k)$ be the reduction of $\rho$ modulo the maximal ideal of $\O$. Assume furthermore that
\begin{itemize}
\item $\rho |_{\Gal(\bQp/ \Qp)}$ has distinct Hodge-Tate weights.
\item $\rhobar$ is modular.
\item $\rhobar$ has non-solvable image if $p=2$; $\rhobar |_{\Gal(\Q(\zeta_p) / \Q)}$ is absolutely irreducible if $p > 2$.
\item if $p = 2$, then $\rhobar |_{\Gal(\bQp/ \Qp)} \not \sim (\begin{smallmatrix} \chi & * \\ 0 & \chi \end{smallmatrix})$ for any character $\chi: \Gal(\bQp/ \Qp) \rightarrow k^{\times}$.
\end{itemize}
Then $\rho$ is modular.
\end{theo}

Such a result is known as a modularity lifting theorem, which says that if $\rhobar$ is modular, then any lift $\rho$ of $\rhobar$ satisfying necessary local conditions is also modular. We note that since we work over $\Q$, the condition on the modularity of $\rhobar$ follows from a deep theorem of Khare-Wintenberger \cite{MR2551764} and Kisin \cite{MR2551765}. 
Establishing a modularity lifting theorem comes down to proving that a certain surjection ${\tRinf} \twoheadrightarrow \T_{\infty}$ of a patched global deformation ring $\tRinf$ onto a patched Hecke algebra $\T_{\infty}$ is an isomorphism after inverting $p$, both of which act on a patched module $\tMinf$ coming from applying the Taylor-Wiles-Kisin method, which uses the third assumption essentially, to algebraic modular forms on a definite quaternion algebra.

A key ingredient in Kisin's approach to the Fontaine-Mazur conjecture is a purely local statement, known as the Breuil-M\'ezard conjecture \cite{MR1944572}, which predicts that $\muG$, the Hilbert-Samuel multiplicity of certain quotients of the framed deformation ring of $\rhobar |_{\Gal(\bQp/ \Qp)}$ parametrizing deformations subjected to $p$-adic Hodge theoretical conditions modulo the maximal ideal of $\O$, is equal to $\muA$, an invariant which can be computed from the representation theory of $\GL_2(\Zp)$ over $k$.  A refined version of this conjecture replacing multiplicities with cycles was formulated by Emerton and Gee in \cite{MR3134019}.

In his work, Kisin establishes a connection between $\tRinf[1/p] \cong \T_{\infty}[1/p]$ and the Breuil-M\'ezard conjecture (when $p > 2$). He shows that $\tRinf \twoheadrightarrow \T_{\infty}$ implies $\muG \geq \muA$, with equality if and only if $\tRinf[1/p] \cong \T_{\infty}[1/p]$. It follows that in each case where one can prove the reverse inequality, one would simultaneously obtain both the Breuil-M\'ezard conjecture and a modularity lifting theorem. A similar argument when $p=2$ was carried out in \cite{MR3544298} using the results of Khare-Wintenberger \cite{MR2551764}.

The key ingredient to prove the reverse inequality $\muG \leq \muA$ is the $p$-adic local Langlands correspondence for $\GL_2(\Qp)$ due to Breuil, Berger, Colmez, Emerton, Kisin and Pa\v{s}k\=unas. The correspondence is given by Colmez's Montreal functor in \cite{MR2642409}, which is an exact, covariant functor $\cV$ sending certain $\GL_2(\Qp)$-representations on $\O$-modules to finite $\O$-modules with a continuous action of $\Gal(\bQp / \Qp)$. Moreover, via reduction modulo $p$ it is compatible with Breuil's (semi-simple) mod $p$ Langlands correspondence in \cite{MR2018825}.

By using the $p$-adic local Langlands correspondence, \cite{MR2505297} deduces the inequality $\muA \geq \muG$ (and thus the Breuil-M\'ezard conjecture) in the cases that $p$ is odd and $\rbar$ ($:= \rhobar  |_{\Gal(\bQp/ \Qp)}$) is not (a twist of) an extension of $\1$ by $\omega$, where $\omega$ is the mod $p$ cyclotomic character. Later on, a purely local proof of the Breuil-M\'ezard conjecture for all continuous representations $\rbar$, which has only scalar endomorphism and is not (a twist of) an extension of $\1$ by $\omega$  if $p = 2, 3$,  is given in \cite{MR3306557, MR3544298} using the results in \cite{MR3150248}. The cases that $\rbar$ is a direct sum of two distinct characters whose ratios are not $\omega$ when $p = 2, 3$ are proved in \cite{MR3429471, MR3671561} by a similar local method. The combined work of Kisin, Hu-Tan and Pa{\v{s}}k\=unas handle the Breuil-M\'ezard conjecture in all cases except when $p=2$ or $3$ and $\rbar \sim (\begin{smallmatrix} \omega\chi & * \\ 0 & \chi \end{smallmatrix})$. 

In \cite{2018arXiv180307451T}, the author gives another proof of this theorem when $p>2$. Instead of proving $\muA \geq \muG$ (or the Breuil-M\'ezard conjecture), we prove $\tRinf[1/p] \cong \T_{\infty}[1/p]$ for automorphic forms on definite unitary groups directly. As a result, the Breuil-M\'ezard conjecture for 2-dimensional Galois representations of $\Gal(\bQp/\Qp)$ follows by a similar equivalence in this setting due to \cite{MR3134019}, which is new in the cases that $p=3$ and $\rbar$ is a twist of the $\1$ by $\omega$. As a result, the theorem is proved.

In this paper, we follow the strategy in \cite{2018arXiv180307451T} to remove the restriction on $\rhobar |_{\Gal(\bQp/ \Qp)}$ when $p=2$. Here is our result:

\begin{theor} \label{thm:FM}
Assume $p=2$. Let $\rho$ be as in the conjecture. Let $\rhobar: \Gal(\bQ / \Q) \rightarrow \GL_2(k)$ be the reduction of $\rho$ modulo the maximal ideal of $\O$. Assume furthermore that
\begin{itemize}
\item $\rho |_{\Gal(\bQp/ \Qp)}$ has distinct Hodge-Tate weights.
\item $\rhobar$ is modular.
\item $\rhobar$ has non-solvable image.
\end{itemize}
Then $\rho$ is modular.
\end{theor} 

Indeed we prove the theorem in a more general context, i.e. $F$ is a totally real field in which $p$ splits completely and $\rho: \Gal(\overline{F} / F) \rightarrow \GL_2(\O)$ (see Theorem \ref{thm:FMC} for the precise statement). We explain our method in more detail below.

Let $p=2$, $G_{\Qp} = \Gal(\bQp/ \Qp)$ be the absolutely Galois group of the field of $p$-adic numbers $\Qp$ and $\rbar : G_{\Qp} \rightarrow \GL_2(k)$ be a continuous representation. We denote the fixed determinant universal framed deformation ring of $\rbar$ by $\Rf_p$. It can be shown that $\rbar$ is isomorphic to the restriction to a decomposition group at $p$ of a mod $p$ Galois representation $\rhobar$ associated to an algebraic modular form on some definite quaternion algebra. By applying the Taylor-Wiles-Kisin patching method in \cite{MR3529394} to algebraic modular forms on a definite quaternion algebra, we construct an $\Rinf$-module $\Minf$ equipped with a commuting action of $\GL_2(\Qp)$, where $\Rinf$ is a complete local noetherian $\Rf_p$-algebra with residue field $k$. For simplicity, one may think of $\Rinf$ as $\Rf_p \llbracket x_1, \cdots, x_m \rrbracket$. In particular, there is no local deformation condition at the place $p$.

If $y \in \mSpec \Rinf[1/p]$, then
\[
\Pi_y := \Homc_{\O}(\Minf \otimes_{\Rinf, y} E_y, E)
\]
is an admissible unitary $E$-Banach space representation of $G$, where $\mSpec(\Rinf[1/p])$ is the set of maximal ideals of $\Rinf[1/p]$ and $E_y$ is the residue field at $y$. Since $\Pi_y$ lies in the range of $p$-adic local Langlands, we may apply the Colmez's functor $\cV$ to $\Pi_y$ and obtain a $\Rinf$-module $\cV(\Pi_y)$ equipped with an action of $G_{\Qp}$. On the other hand, the composition $x: \Rf_p \rightarrow \Rinf \xrightarrow{y} E_y$ defines a continuous Galois representation $r_x: G_{\Qp} \rightarrow \GL_2(E_y)$. It is expected that the Banach space representation $\Pi_y$ depends only on $x$ (see \cite{MR3732208}) and that it should be related to $r_x$ by the $p$-adic local Langlands correspondence (see Theorem \ref{thm:ColMinf} below).

Our patched module $\Minf$ is related to Kisin's $\tMinf$ as follows. The patching in Kisin's paper is always with fixed Hodge-Tate weights and a fixed inertial type. This information can be encoded in an irreducible locally algebraic representation $\sigma$ of $\GL_2(\Zp)$ over $E$. Let $\Rf_p(\sigma)$ be quotient of $\Rf_p$ parameterizing the lifts of $\rhobar$ of type $\sigma$. We define $\Rinf(\sigma) = \Rinf \otimes_{\Rf_p} \Rf_p(\sigma)$ (which is Kisin's patched global deformation ring $\tRinf$ introduced before) and $\Minf(\sigma^{\circ}) = \Minf \ctimes_{\O \llbracket \GL_2(\Zp) \rrbracket} \sigma^{\circ}$ with $\sigma^{\circ}$ a $\GL_2(\Zp)$-stable $\O$-lattice of $\sigma$. Then $\Minf(\sigma^{\circ})$ is a finitely generated $\Rinf$-module with the action of $\Rinf$ factoring through $\Rinf(\sigma)$. Moreover, an argument using the Auslander-Buchsbaum formula shows that the support of $\Minf(\sigma^{\circ})$ is equal to a union of irreducible components of $\Rinf(\sigma)$. It can be shown that Kisin's patched module $\tMinf$ is isomorphic to $\Minf(\sigma^{\circ})$. The main theorem in this paper is the following:

\begin{theor} \label{theo:main}
Every irreducible component of $\tRinf$ is contained in the support of $\tMinf$.
\end{theor}

By the local-global compatibility for the patched module $\Minf$, this amounts to showing that if $r_x$ is de Rham with distinct Hodge-Tate weights, then (a subspace of) locally algebraic vectors in $\Pi_y$ can be related to $\WD(r_x)$ via the classical local Langlands correspondence, where $\WD(r_x)$ is the Weil-Deligne representation associated to $r_x$ defined by Fontaine. 

One of the ingredients to show this is a result in \cite{2018arXiv180906598E}, which implies that the action of $\Rinf$ on $\Minf$ is faithful. Note that this does not imply that $\Pi_y \neq 0$ since $\Minf$ is not finitely generated over $\Rinf$. In \cite{2018arXiv180307451T}, this issue has been overcome by applying Colmez's functor $\cV$ to $\Minf$ and showing that $\cV(\Minf)$ is a finitely generated $\Rinf$-module. Let us note that a similar finiteness result has been proved in \cite{2019arXiv190107166P} using results of \cite{MR3150248}. Our proof is different since results of \cite{MR3150248, MR3544298} are not available when $p=2$ and $\rbar$ has scalar semisimplification.

Since $\cV(\Minf)$ is a finitely generated $\Rinf$-module, the specialization of $\cV(\Minf)$ at any $y \in \mSpec \Rinf[1/p]$ is non-zero by Nakayama's lemma, which in turn implies that $\Pi_y$ is nonzero. Combining these, results from $p$-adic local Langlands, and a result in \cite{MR1094193} which says that a $2$-dimensional absolutely irreducible Galois representation is isomorphic to its associated Cayley-Hamilton algebra, we prove the following:

\begin{theor} \label{thm:ColMinf}
If $r_x$ is absolutely irreducible, then $\cV(\Pi_y) \cong r_x^{\oplus n_y}$ for some positive integer $n_y$. Moreover, $n_y = 1$ in a dense subset of $\mSpec \Rinf[1/p]$.
\end{theor}

This shows that Kisin's patched module $\tMinf$ is supported at every generic point whose associated local Galois representation at place $p$ is absolutely irreducible. So we only have to handle the reducible (thus ordinary) locus, which can be shown to be modular by using an ordinary modularity lifting theorem, which is an analog of \cite{MR2941425, MR3252020, MR3904451, Sasaki2018II} in our setting. This finishes the proof of Theorem \ref{theo:main} and gives a new proof of the Breuil-M\'ezard conjecture by the formalism in \cite{MR2505297, MR3292675, MR3134019, MR3306557}, which is new in the cases that $p=2$ and $\rbar$ is a twist of $\1$ by itself (note that $\omega \cong \1$ when $p=2$). As a consequence, we prove new cases of Fontaine-Mazur conjecture. We remark that by using the patching in \cite{MR2505297}, our method applies to the case $p>2$ without any change. We focus only on the case $p=2$ since this is the only remaining case with the restriction on $\rhobar |_{\Gal(\bQp/ \Qp)}$.

Note that our method for Theorem \ref{thm:ColMinf} doesn't apply to the case that $r_x$ is reducible since the characteristic polynomial only determines a Galois representation up to semi-simplification. Nevertheless, the same conclusion can be deduced from existing local-global compatibility results when $r_x$ is crystabelline \cite{MR3347316} or when $r_x$ is semi-stable \cite{MR3510330}.

The paper is organized as follows. We first recall some background knowledge and properties in Sects. \ref{section:prerepthory}, \ref{section:autoform} and \ref{section:Galoisdeform} on representation theory, automorphic forms and Galois deformation theory respectively. In Sect. \ref{section:patchingargument}, we introduce completed cohomology and construct the patched module. We relate our patched module to the Breuil-M\'ezard conjecture in Sect. \ref{section:patchBM} and to the $p$-adic Langlands correspondence in Sect. \ref{section:patchpLLC} using a faithfulness result in \cite{2018arXiv180906598E}. In Sect. \ref{section:ordpatch}, we construct some partially ordinary Galois representations by an ordinary $R = \T$ theorem. In Sect. \ref{section:main}, we put all these results together and prove our main theorem, and use it to give a new proof of the Breuil-M\'ezard conjecture and the Fontaine-Mazur conjecture.

\subsection*{Acknowledgement}
I would like to thank my advisor Vytautas Pa{\v{s}}k\=unas, for suggesting me to work on this project and sharing with his profound insight and ideas. I also thanks Shu Sasaki for many helpful discussions on modularity lifting theorems and for pointing out many inaccuracies in an earlier draft, and Jack Thorne for his hospitality during my visit to Cambridge in May 2018 and for answering my questions regarding $2$-adic modularity lifting theorems. I would also like to thank Patrick Allen and the anonymous referee for many useful suggestions, comments, and corrections. This research was funded in part by the DFG, SFB/TR 45 "Periods, moduli spaces and arithmetic of algebraic varieties”.

\section*{Notations} \label{Notation}
If $F$ is a field with a fixed algebraic closure $\overline{F}$, then we write $G_F = \Gal(\overline{F}/F)$ for its absolutely Galois group. We write $\varepsilon: G_F \rightarrow \Zp^{\times}$ for the $p$-adic cyclotomic character, and $\omega$ for the mod $p$ cyclotomic character. If $F$ is a finite extension of $\Qp$, we write $I_F$ for the inertia subgroup of $G_F$, $\varpi_F$ for a uniformizer of the ring of integers $\O_F$ of $F$ and $k_F = \O_F / \varpi_F$ its residual field.

If $F$ is a number field and $v$ is a place of $F$, we let $F_v$ be the completion of $F$ at $v$ and $\A_F$ its ring of adeles. If $S$ is a finite set of places of $F$, we let $\A_F^S$ denote the resticted tensor product $\prod'_{v \notin S} F_v$. In particular, $\Af$ denotes the ring of finite adeles. For each finite place $v$ of $F$, we will denote by $q_v$ the order of residue field at $v$, and by $\varpi_v \in F_v$ a uniformizer and $\Frob_v$ an arithmetic Frobenius element of $G_{F_v}$. 

We let 
\[
\Art_F = \prod_v \Art_{F_v}: \A_F^{\times} / \overline{F^{\times} (F_{\infty}^{\times})^{\circ} } \xrightarrow{\sim} G_F^{\ab}
\]
be the global Artin map, where the local Artin map $\Art_{F_v}: F^{\times}_v \rightarrow W_{F_v}^{ab}$ is the isomorphism provided by local class field theory, which sends our fixed uniformizer to a geometric Frobenius element.

We fix a finite extension $E / \Qp$ sufficiently large in the sense that all embeddings $F \rightarrow \bQp$ have image lying in $E$. We denote $\O$ the ring of integers of $E$ and $k$ its residue field.

We will consider a locally algebraic character $\psi: \A_F^{\times} / \overline{F^{\times} (F_{\infty}^{\times})^{\circ} } \rightarrow \O^{\times}$ in the sense that there exists an open compact subgroup $U$ of $\Aft$ such that $\psi(u) = \prod_{v | p} \Norm_v(u_v)^{t_v}$ for $u \in U$, where $u_v$ is the projection of $u$ to the place $v$, $\Norm_v$ the local norm, and $t_v$ an integer. When $\overline{F^{\times} (F_{\infty}^{\times})}$ lies in the kernel of $\psi$, we consider $\psi$ as a character $\psi: \Aft / F^{\times} \rightarrow \O^{\times}$, whose corresponding Galois character is totally even.

Let $W$ be a de Rham representation of $G_{\Qp}$ over $E$. We will write $\HT(W)$ for the set of Hodge-Tate weights of $W$ normalized by $\HT(\varepsilon) = \{ -1 \}$. We say that $W$ is regular if $\HT(W)$ are pairwise distinct. Let $\Z^2_+$ denote the set of tuples $(\lambda_1, \lambda_2)$ of integers with $\lambda_1 \geq \lambda_2$. If $W$ be a 2-dimensional de Rham representation which is regular, then there is a $\lambda = (\lambda_1, \lambda_2) \in \Z^2_+$ such that $\HT(W) = \{\lambda_{2}, \lambda_{1} + 1\}$, and we say that $W$ is regular of weight $\lambda$. 

For any $\lambda \in \Z^2_+$, we write $\Xi_{\lambda} = \Sym^{\lambda_1 - \lambda_2} \otimes \det^{\lambda_2}$ for the algebraic $\Zp$-representation of $\GL_2$ with highest weight $\lambda$ and $M_{\lambda}$ for the $\O$-representation of $\GL_2(\O_{\Qp})$ obtained by evaluating $\Xi_{\lambda}$ on $\Zp$.

An inertial type is a representation $\tau: I_{\Qp} \rightarrow \GL_2(\bQp)$ with open kernel which extends to the Weil group $W_{\Qp}$. We say a de Rham representation $\rho: G_{\Qp} \rightarrow \GL_2(E)$ has inertial type $\tau$ if the restriction to $I_{\Qp}$ of the Weil-Deligne representation $\WD(\rho)$ associated to $\rho$ (see \cite{MR1293977} for the precise definition) is equivalent to $\tau$. Given an inertia type $\tau$, by a result of Henniart in the appendix of \cite{MR1944572}, there is a (unique if $p>2$) finite dimensional smooth irreducible $\bQp$-representation $\sigma(\tau)$ \big(resp. $\sigma^{cr}(\tau)$\big) of $\GL_2(\Z_p)$, such that for any infinite dimensional smooth absolutely irreducible representation $\pi$ of $G$ and the associated Weil-Deligne representation $\LL(\pi)$ attached to $\pi$ via the classical local Langlands correspondence, we have $\Hom_K(\sigma(\tau), \pi) \neq 0$ (resp. $\Hom_K(\sigma^{cr}(\tau), \pi) \neq 0$) if and only if $\LL(\pi) \vert_{I_{\Qp}} \cong \tau$ (resp. $\LL(\pi) \vert_{I_{\Qp}} \cong \tau$ and the monodromy operator $N$ is trivial). Enlarging $E$ if needed, we may assume $\sigma(\tau)$ is defined over $E$.

If $L$ be a finite extension of $\Qp$, we let $\rec$ for the local Langlands correspondence for $\GL_2(L)$, as defined in \cite{MR2234120, MR1876802}. By definition, it is a bijection between the set of isomorphism classes of irreducible admissible representation of $\GL_2(L)$ over $\C$, and the set of Frobenius semi-simple Weil-Deligne representation of $W_L$ over $\C$. Fix once and for all an isomorphism $\iota: \bQp \xrightarrow{\sim} \C$. We define the local Langlands correspondence $\rec_p$ over $\bQp$ by $\iota \circ \rec_p = \rec \circ \iota$, which depends only on $\iota^{-1}(\sqrt{p})$. If we set $r_p(\pi):= \rec_p(\pi \otimes |\det|^{-1/2})$, then $r_p$ is independent of the choice of $\iota$. Furthermore, if $V$ is a Frobenius semi-simple Weil-Deligne representation Weil-Deligne representation of $W_L$ over $E$, then $r_p^{-1}(V)$ is also defined over $E$.

If $r: G_{\Qp} \rightarrow \GL_2(E)$ is de Rham of regular weight $\lambda$, then we write $\pialg(r) = M_\lambda \otimes_{\O} E$, $\pism(r) = r_p^{-1}(\WD(r_x)^{F-ss})$ and $\pi_{\lalg}(r) = \pialg(r) \otimes \pism(r)$, all of which are $E$-representations of $\GL_2(\Qp)$.

Recall that a linearly topological $\O$-module is a topological $\O$-module which has a fundamental system of open neighborhoods of the identity which are $\O$-submodules. If $A$ is a linear topological $\O$-module, we write $A^{\vee}$ for its Pontryagin dual $\Homc_{\O}(A, E/\O)$, where $E / \O$ has the discrete topology, and we give $A^{\vee}$ the compact open topology. We write $A^d$ for the Schikhof dual $\Homc_{\O}(A, \O)$, which induces an anti-equivalence of categories between the category of compact, $\O$-torsion free linear-topological $\O$-modules $A$ and the category of $\varpi$-adically complete separated $\O$-torsion free $\O$-modules. A quasi-inverse is given by $B \mapsto B^d := \Hom_{\O}(B, \O)$, where the target is given the weak topology of pointwise convergence. Note that if $A$ is an $\O$-torsion free profinite linearly topological $\O$-module, then $A^d$ is the unit ball in the $E$-Banach space $\Hom_{\O}(A, E)$.

For $R$ a Noetherian local ring with maximal ideal $\m$ and $M$ a finite $R$-module, let $e(M, R)$ denote the Hilbert-Samuel multiplicity of $M$ with respect to $\m$. We abbreviate $e(R,R)$ for $e(R)$. For $R$ a Noetherian ring and $M$ a finite $R$-module of dimension at most $d$., let $\ell_{R_{\p}}(M_{\p})$ denote the length of the $R_{\p}$-module $M_{\p}$, and let $Z_d(M) = \sum_{\p} \ell_{R_{\p}}(M_{\p}) \p$ for all $\p \in \Spec R$ such that $\dim R / \p = d$.  If $M$ and $N$ are finitely generated $R$- and $S$-module of dimension at most $d$ and $e$ respectively, then the completed tensor product $M \ctimes_k N$ is of dimension $d+e$, and $Z_d(M) \times_k Z_e(N)$ is equal to $Z_{d+e}(M \ctimes_k N)$. We refer the reader to \cite[\S 2]{MR3134019} for details.

Let $(A, \m)$ be a complete local $\O$-algebra with maximal ideal $\m$ and residue field $k = A / \m$, we will denote $\CNL_{A}$ the category of complete local $A$-algebra with residue field $k$.
\section{Preliminaries in representation theory} \label{section:prerepthory}
\subsection{Generalities}
Let $G$ be a $p$-adic analytic group, $K$ be a compact open subgroup of $G$, and $Z$ be the center of $G$.

Let $(A, \m_A) \in \CNL_{\O}$. We denote by $\Mod_G(A)$ the category of $A[G]$-modules and by $\Modsm(A)$ the full subcategory with objects $V$ such that $V = \cup_{H, n} V^H[\m^n]$, where the union is taken over all open subgroups of $G$ and integers $n \geq 1$ and $V[\m^n]$ denotes elements of $V$ killed by all elements of $\m^n$. Let $\Modlfin(A)$ be the full subcategory of $\Modsm(A)$ with objects smooth $G$-representation which are locally of finite length, this means for every $v \in V$, the smallest $A[G]$-submodule of $V$ containing $v$ is of finite length.

An object $V$ of $\Modsm(A)$ is called admissible if $V^H[\m^i]$ if a finitely generated $A$-module for every open subgroup $H$ of $G$ and every $i \geq 1$; $V$ is called locally admissible if for every $v \in V$ the smallest $A[G]$-submodule of $V$ containing $v$ is admissible. Let $\Modladm(A)$ be the full subcategory of $\Modsm(A)$ consisting of locally admissible representations. 

For a continuous character $\zeta : Z \rightarrow A^\times$, adding the subscript $\zeta$ in any of the above categories indicates the corresponding full subcategory of $G$-representations with central character $\zeta$. These categories are abelian and are closed under direct sums, direct limits and subquotients. Note that if $G = \GL_2(\Qp)$ or $G$ is a torus then $\Modlfinz(A) = \Modladmz(A)$ \cite[Theorem 2.3.8]{MR2667882}.

Let $H$ be a compact open subgroup of $G$ and $A \llbracket H \rrbracket$ the completed group algebra of $H$. Let $\Modpro_G(A)$ be the category of profinite linearly topological $A \llbracket H \rrbracket$-modules with an action of $A[G]$ such that the two actions are the same when restricted to $A[H]$ with morphisms $G$-equivariant continuous homomorphisms of topological $A \llbracket H \rrbracket$-modules. The definition does not depend on $H$ since any two compact open subgroups of $G$ are commensurable. By \cite[Lemma 2.2.7]{MR2667882}, this category is anti-equivalent to $\Modsm(A)$ under the Pontryagin dual $V \mapsto V^\vee := \Hom_{\O}(V, E / \O)$ with the former being equipped with the discrete topology and the latter with the compact-open topology. We denote $\fC(A)$ the full subcategory of $\Modpro_{G}(A)$ anti-equivalent to $\Modlfinz(A)$.

An $E$-Banach space representation $\Pi$ of $G$ is an $E$-Banach space $\Pi$ together with a $G$-action by continuous linear automorphisms such that the inducing map $G \times \Pi \rightarrow \Pi$ is continuous. A Banach space representation $\Pi$ is called unitary if there is a $G$-invariant norm defining the topology on $\Pi$, which is equivalent to the existence of an open bounded $G$-invariant $\O$-lattice $\Theta$ in $\Pi$. An unitary $E$-Banach space representation is admissible if $\Theta \otimes_{\O} k$ is an admissible smooth representation of $G$, which is independent of the choice of $\Theta$. We denote $\Ban(E)$ the category of admissible unitary $E$-Banach space representations on which $Z$ acts by $\zeta$.

\subsection{Representations of $\GL_2(\Qp)$} \label{RepGL2}
In this subsection, we assume $p=2$, $G = \GL_2(\Qp)$, $K = \GL_2(\Zp)$, and thus $Z \simeq \Qp^\times$. Let $B$ be the subgroup of upper triangular matrices in $G$. If $\chi_1$ and $\chi_2$ are characters of $\Qp^{\times}$, then we write $\chi_1 \otimes \chi_2$ for the character of $B$ which maps $(\begin{smallmatrix}a & b \\ 0 & d \end{smallmatrix})$ to $\chi_1(a) \chi_2(d)$.

By a Serre weight we mean an absolutely irreducible representation of $K$ on an $k$-vector space. It is of the form $\bsigma_a := \Sym^{a_1 - a_2} k^2 \otimes \det^{a_2}$ for a unique $a = (a_1, a_2) \in \Z^2$ with $a_1 - a_2 \in \{0 , \dots, p-1 \}$ and $a_2 \in \{ 0, \dots p-2\}$. We call such pairs $a$ Serre weights also.

Let $\sigma$ be a Serre weight. There exists an isomorphism of algebras
\[
\End_G(\cInd^G_K \sigma) \cong k[T, S^{\pm 1}]
\]
for certain Hecke operators $T, S \in \End_G(\cInd^G_K \sigma)$. It follows from \cite[Theorem 33]{MR1290194} and \cite[Theorem 1.6]{MR2018825} that the absolutely irreducible smooth $k$-representations of $G$ with a central character fall into four disjoint classes:
\begin{itemize}
\item characters $\eta \circ \det$;
\item special series $\Sp \otimes \eta \circ \det$;
\item principal series $\Ind_B^G(\chi_1 \otimes \chi_2)$, with $\chi_1 \neq \chi_2$;
\item supersingular $\cInd_K^G(\sigma) / (T, S-\lambda)$, with $\lambda \in k^{\times}$,
\end{itemize}
where the Steinberg representation $\Sp$ is defined by the exact sequence
\[
0 \rightarrow \1 \rightarrow \Ind_B^G \1 \rightarrow \Sp \rightarrow 0.
\]

\subsubsection{Blocks}
Let $\Irr$ be the set of equivalent classes of smooth irreducible $k$-representations of $G$ with central character $\zeta$. We write $\pi \leftrightarrow \pi'$ if $\pi \cong \pi'$ or $\Ext^1_{G, \zeta}(\pi, \pi') \neq 0$ or $\Ext^1_{G, \zeta}(\pi', \pi) \neq 0$, where $\Ext^1_{G, \zeta}(\pi, \pi')$ is the Yoneda extension group of $\pi'$ by $\pi$ in $\Modlfinz(k)$. We write $\pi \sim \pi'$ if there exists $\pi_1, \cdots, \pi_n \in \Irr$ such that $\pi \cong \pi_1$, $\pi' \cong \pi_n$ and $\pi_i \leftrightarrow \pi_{i+1}$ for $1 \leq i \leq n-1$. The relation $\sim$ is an equivalence relation on $\Irr$. A block is an equivalence class of $\sim$. The classification of blocks can be found in \cite[Corollary 1.2]{MR3444235}. Moreover, by \cite[Proposition 5.34]{MR3150248}, the category $\Modlfinz(\O)$ decomposes into a direct sum of subcategories
\begin{align}
\Modlfinz(\O) \cong \prod_{\B} \Modlfinz(\O)[\B]
\end{align}
where the product is taken over all the blocks $\B$ and the objects of $\Modlfinz(\O)[\B]$ are representations with all the irreducible subquotients in $\B$. Dually we obtain
\begin{align} \label{equation:Cblockdecomp}
\fC(\O) \cong \prod_{\B} \fC(\O)[\B],
\end{align}
where $\fC(\O)[\B]$ is the full subcategory of $\fC(\O)$ defined by $\Modlfinz(\O)[\B]$ under the anti-equivalence.

\begin{lemma} \label{lemma:SL2thick}
Let $0 \rightarrow \pi_1 \rightarrow \pi_2 \rightarrow \pi_3 \rightarrow 0$ be an extension in $\Modsm(\O)$ then $\SL_2(\Qp)$ acts trivially on $\pi_1$ and $\pi_3$ if an only if it acts trivially on $\pi_2$.
\end{lemma}

\begin{proof}
If $\SL_2(\Qp)$ acts trivially on $\pi_1$ and $\pi_3$, then $\pi_1 \subset \pi_2^{\SL_2(\Qp)}$ and thus $\pi_2 / \pi_2^{\SL_2(\Qp)}$ is a quotient of $\pi_3$. It follows that $\SL_2(\Qp)$ acts trivially on $\pi_2 / \pi_2^{\SL_2(\Qp)}$. On the other hand, it is proved in \cite[Lemma \uppercase\expandafter{\romannumeral3}.40]{MR3267142} that $\pi_2 / \pi_2^{\SL_2(\Qp)}$ has no $\SL_2(\Qp)$-invariant. Hence $\pi_2 / \pi_2^{\SL_2(\Qp)} = 0$. The other implication is trivial. 
\end{proof}

Let $\fT(\O)$ be the full subcategory of $\fC(\O)$ whose objects have trivial $\SL_2(\Qp)$-action. It follows from Lemma \ref{lemma:SL2thick} that $\fT(\O)$ is a thick subcategory of $\fC(\O)$ and hence we may consider the quotient category $\fD(\O) := \fC(\O) / \fT(\O)$. Note that the objects of $\fD(\O)$ is same as the objects of $\fC(\O)$ and the morphisms are given by
\[
\Hom_{\fD}(M, N) := \varinjlim \Hom_{\fC}(M', N/N'),
\]
where the limit is taken over all subobjects $M'$ of $M$ and $N'$ of $N$ such that $\SL_2(\Qp)$ acts trivially on $M/M'$ and $N'$. Let $\iT: \fC(\O) \rightarrow \fD(\O)$ be the functor $\iT M = M$ for every object of $\fC(\O)$ and $\iT f$ the image of $f: M \rightarrow N$ in $\varinjlim \Hom_{\fC}(M', N/N')$ under the natural map. Moreover, $\fD(\O)$ is an abelian category and $\iT$ is an exact functor. We denote $\fD(k)$ the full subcategory of $\fD(\O)$ consisting of objects killed by $\varpi$.

Let $\overline{\zeta}$ be the reduction modulo $\varpi$ of $\zeta$. Note that $(\overline{\zeta} \circ \det)^{\vee}$ is the only absolutely irreducible object in $\fC(\O)$ with trivial $\SL_2(\Qp)$-action. The following proposition is an easy variant of \cite[Lemma 10.26, Lemma 10.27, Lemma 10.28, Lemma 10.29]{MR3150248}. We leave the proof to the reader.

\begin{prop} \label{prop:quotcat} {\ }
\begin{enumerate}
\item Let $M$ and $N$ be objects of $\fC(\O)$. We have
\[
\Hom_{\fD(\O)}(\iT M, \iT N) \cong \Hom_{\fC(\O)}(I_{\SL_2(\Qp)}(M), N/N^{\SL_2(\Qp)}),
\]
where $I_{\SL_2(\Qp)}(M) = \big(M^{\vee} / (M^{\vee})^{\SL_2(\Qp)} \big)^\vee$.
\item If $P$ is a projective object of $\fC(\O)$ with $\Hom_{\fC(\O)}(P, (\overline{\zeta} \circ \det)^{\vee}) = 0$ then $\iT P$ is a projective object of $\fD(\O)$ and
\[
\Hom_{\fC(\O)}(P, N) \cong \Hom_{\fD(\O)}(\iT P, \iT N)
\]
for all $N$. Moreover, the category $\fD(\O)$ has enough projectives.
\item If $\Hom_{\fC(\O)}(N, (\overline{\zeta} \circ \det)^{\vee}) = 0$ then for every essential epimorphism $q: M \twoheadrightarrow N$, $\iT q: \iT M \twoheadrightarrow \iT N$ is an essential epimorphism in $\fD(\O)$.
\end{enumerate}
\end{prop}

Since $\fT(\O)$ is contained in $\fC(\O)[\B]$ with $\B = \{\overline{\zeta} \circ \det, \Sp \otimes \overline{\zeta} \circ \det \}$, we may build the quotient category $\fD(\O)[\B] / \fT(\O)$. We write $\fD(\O)[\B]$ for $\fC(\O)[\B]$ for other blocks and thus (\ref{equation:Cblockdecomp}) induces a decomposition of categories
\[
\fD(\O) \cong \prod_{\B} \fD(\O)[\B].
\]

\subsubsection{Colmez's Montreal functor}
Let $\Modfinz(\O)$ be the full subcategory of $\Modsm(\O)$ consisting of representations of finite length with a central character. Let $\ModGQp(\O)$ be the category of continuous $\GQp$-representations on $\O$-modules of finite length with the discrete topology. In \cite{MR2642409}, Colmez has defined an exact and covariant functor $\V: \Modfinz(\O) \rightarrow \ModGQp(\O)$. If $\psi: \Qp^{\times} \rightarrow \O^{\times}$ is a continuous character, then we may also consider it as a continuous character $\psi: \GQp \rightarrow \O^{\times}$ via class field theory and for all $\pi \in \Modsmz(\O)$ of finite length we have $\V(\pi \otimes \psi \circ \det) \cong \V(\pi) \otimes \psi$. 

Moreover, it follows from the construction in the loc. cit. that $\V(\1) = 0$, $\V(\Sp) = \omega$, $\V(\Ind^G_B \chi_1 \otimes \chi_2) \cong \chi_2$, and $\V(\cInd \Sym^r k^2 / (T, S-1)) \cong \ind{\omega_2^{r+1}}$, where $\omega_2: I_{\Qp} \rightarrow k^{\times}$ is Serre's fundamental character of level $2$, and $\ind{\omega_2^{r+1}}$ is the unique irreducible representation of $\GQp$ of determinant $\omega^r$ and such that $\ind{\omega_2^{r+1}} |_{I_{\Qp}} \cong \omega_2^{r+1} \oplus \omega_2^{2(r+1)}$ with $0 \leq r \leq 1$. Note that this determines the image of supersingular representations under $\V$ completely since every supersingular representation is isomorphic to $\cInd \Sym^r k^2 / (T, S-1)$ for some $0 \leq r \leq 1$ after twisting by a character.

Let $\Modpro_{\GQp}(\O)$ be the category of continuous $\GQp$-representations on compact $\O$-modules. Following \cite[\S 3]{MR3306557}, we define an exact covariant functor $\cV: \fC(\O) \rightarrow \Modpro_{\GQp}(\O)$ as follows: Let $M$ be in $\fC(\O)$, if it is of finite length, we define $\cV(M) := \V(M^\vee)^\vee(\varepsilon \psi)$ where $\vee$ denotes the Pontryagin dual. For general $M \in \fC(\O)$, write $M \cong \varprojlim M_i$, with $M_i$ of finite length in $\fC(\O)$ and define $\cV(M) := \varprojlim \cV(M_i)$. With this normalization, we have 
\begin{itemize}
\item $\cV(\pi^{\vee}) = 0$ if $\pi \cong \eta \circ \det$;
\item $\cV(\pi^{\vee}) \cong \chi_1$ if $\pi \cong \Ind^G_B \chi_1 \otimes \chi_2$;
\item $\cV(\pi^{\vee}) \cong \eta$ if $\pi \cong \Sp \otimes \eta \circ \det$;
\item $\cV(\pi^{\vee}) \cong \V(\pi)$ if $\pi$ is supersingular.
\end{itemize}
The functor $\cV: \fC(\O) \rightarrow \Modpro_{\GQp}(\O)$ kills characters and hence every object in $\fT(\O)$. Hence $\cV$ factors through $\iT: \fC(\O) \rightarrow \fD(\O)$. We denote $\cV: \fD(\O) \rightarrow \Modpro_{\GQp}(\O)$ by the same letter.

Let $\Pi \in \Ban(E)$, we define $\cV(\Pi) = \cV(\Theta^d) \otimes_\O E$ with $\Theta$ any open bounded $G$-invariant $\O$-lattice in $\Pi$, so that $\cV$ is exact and contravariant on $\Ban(E)$. Note that $\cV(\Pi)$ does not depend on the choice of $\Theta$.

\subsubsection{Extension Computations when $p=2$ and $\B =\{\1, \Sp\}$}
In this subsection, we do some similar computations as in \cite[\S 10]{MR3150248} when $p=2$, $\B =\{\1, \Sp \}$ and $\zeta = \1$. We write $\operatorname{Mod}^{\operatorname{l.fin}}_{G/Z}(k)$ for $\operatorname{Mod}^{\operatorname{l.fin}}_{G, \1}(\O)$ and $e(\pi', \pi) := \dim_k \Ext^1_{G/Z}(\pi', \pi)$ with $\pi', \pi \in \operatorname{Mod}^{\operatorname{l.fin}}_{G/Z}(k)$.

\begin{lemma} \label{lemma:exttristeinberg}
We have $e(\Sp, \1) = 1$. In particular, the unique non-split extension of $\Sp$ by $\1$ is $\Ind^G_B \1$.
\end{lemma}

\begin{proof}
Applying $\Hom_{G/Z}(-,\1)$ to the short exact sequence
\begin{align} \label{equation:sesInd}
    0 \rightarrow \1 \rightarrow \Ind^G_B \1 \rightarrow \Sp \rightarrow 0,
\end{align}
we obtain the following long exact sequence
\[
0 \rightarrow \Hom_G(\1, \1) \rightarrow \Ext^1_{G/Z}(\Sp, \1) \rightarrow \Ext^1_{G/Z}(\Ind^G_B \1, \1) \xrightarrow{f} \Ext^1_{G/Z}(\1, \1).
\]
Since $e(\Ind^G_B \1, \1) = 1$ by \cite[Theorem 4.3.13 (2)]{MR2667883}, we have $e(\Sp, \1)$ is 2 if $f$ is the zero map and 1 otherwise.

On the other hand, we have the exact sequence
\[
0 \rightarrow \Ext^1_{\H}(\I(\Ind^G_B \1), \I(\1)) \rightarrow \Ext^1_{G/Z}(\Ind^G_B \1, \1) \rightarrow \Hom_{\H}(\I(\Ind^G_B \1), \Rone \I(\1))
\]
coming from low degree terms associated to the $E_2$-spectral sequence given by the pro-$p$ Iwahori invariant functor $\I$ \cite[Proposition 9.1]{MR2667891}, where $\H$ is the (fixed determinant) pro-$p$ Iwahori Hecke algebra (same as the Iwahori Hecke algebra since Iwahori subgroups are pro-$p$ when $p=2$) and $\I$ is the pro-$p$ Iwahori invariant functor. We claim that $\Ext^1_{\H}(\I(\Ind^G_B \1), \I(\1))$ is nonzero. 

Suppose the claim holds. Note that there is a short exact sequence
\begin{align} \label{equation:IsesInd}
    0 \rightarrow  \I(\1) \rightarrow  \I(\Ind^G_B \1) \rightarrow \I(\Sp) \rightarrow 0
\end{align}
coming from applying $\I$ to (\ref{equation:sesInd}) by \cite[Corollary 6.4]{MR2931521}. Applying $\Hom_{\H}(-, \I(\1))$ to (\ref{equation:IsesInd}), we obtain the following exact sequence
\[
0 \rightarrow \Hom_{\H}(\I(\1), \I(\1)) \rightarrow \Ext^1_{\H}(\I(\Sp), \I(\1)) \rightarrow \Ext^1_{\H}(\I(\Ind^G_B \1), \I(\1)) \rightarrow \Ext^1_{\H}(\I(\1), \I(\1)).
\]
Since $\Ext^1_{\H}(\I(\Sp), \I(\1))$ is 1-dimensional \cite[Lemma 11.3]{MR2667891}, we see that the last map is an injection. It follows that we have the following commutative diagram
\[
  \begin{tikzcd}
  &\Ext^1_{\H}(\I(\Ind^G_B \1), \I(\1)) \arrow[r, hookrightarrow] \arrow[d, hookrightarrow] 
  &\Ext^1_{\H}(\I(\1), \I(\1)) \arrow[d, hookrightarrow] \\   
  &\Ext^1_{G/Z}(\Ind^G_B \1, \1) \arrow[r, "f"]
  &\Ext^1_{G/Z}(\1, \1),
  \end{tikzcd}
\]
where the horizontal maps are induced by functoriality and the vertical maps come from the low degree terms associated to the $E_2$-spectral sequence given by $\I$. This proves the lemma since any nonzero element in $\Ext^1_{\H}(\I(\Ind^G_B \1), \I(\1))$ would give rise to an element of $\Ext^1_{G/Z}(\Ind^G_B \1, \1)$ whose image under $f$ is nonzero.

To prove the claim, we construct a non-trivial extension of $\I(\Ind^G_B \1)$ by $\I(\1)$ explicitly. Note that $\H$ is the $k$-algebra with two generators $T, S$ satisfying two relations $T^2 = 1$ and $(S+1)S = 0$. Moreover, $\I(\1)$ is the simple (right) $\H$-module given by $vT = v; \ vS = 0$, $\I(\Sp)$ is the simple $\H$-module given by $vT = v; \ vS = v$, and $\I(\Ind^G_B \1)$ is the $\H$-module given by $ v_1 T = v_1; \ v_2 T = v_2; \ v_1 S = 0; \ v_2 S = v_1 + v_2$ (c.f. \cite[\S 1.1]{MR2027193}). Since the unique non-split extension of $\I(\1)$ by itself is given by $v_1 T = v_1; \ v_2 T = v_1 + v_2; \ v_1 S = 0; \ v_2 S = 0$ (note that $2=0$ in $k$), it follows that
\begin{align*}
v_1 T = v_1 \quad v_2 T = v_1 + v_2 \quad v_3 T = v_3; \\
v_1 S = 0 \quad v_2 S = 0 \quad v_3 S = v_2 + v_3
\end{align*}
gives a desired non-trivial element in $\Ext^1_{\H}(\I(\Ind^G_B \1), \I(\1))$.
\end{proof}

By \cite[Proposition 4.3.21, Proposition 4.3.22]{MR2667883}, \cite[Proposition \uppercase\expandafter{\romannumeral7}.4.18]{MR2642409} and the above lemma, we have the following table for $e(\pi', \pi)$:\\
\begin{tabular}{c|cc} 
$\pi' \backslash \pi$ & $\1$ & $\Sp$ \\ \hline 
$\1$ & 3 & 3 \\ 
$\Sp$ & 1 & 3 \\  
\end{tabular}

\begin{lemma} \label{lemma:extIndSp}
The natural map $\Ext^1_{G/Z}(\Sp, \Sp) \rightarrow \Ext^1_{G/Z}(\Ind^G_B \1, \Sp)$ is a bijection.
\end{lemma}

\begin{proof}
Consider the exact sequence
\[
0 \rightarrow \Ext^1_{G/Z}(\Sp, \Sp) \rightarrow \Ext^1_{G/Z}(\Ind^G_B \1, \Sp) \rightarrow \Ext^1_{G/Z}(\1, \Sp).
\]
coming from applying $\Hom_{G}(-, \Sp)$ to the short exact sequence $0 \rightarrow \1 \rightarrow \Ind^G_B \1 \rightarrow \Sp \rightarrow 0$. Since $e(\Ind^G_B \1, \Sp) = 3$ by \cite[Theorem 4.3.12 (2)]{MR2667883}, we see that the first map is a bijection and the second map is identically zero.
\end{proof}

Since $e(\1, \Sp) = 3$ there exists a unique smooth $k$-representation $\kappa$ with socle $\Sp$ and have an exact sequence:
\begin{equation} \label{equation:tau}
0 \rightarrow \Sp \rightarrow \kappa \rightarrow \1^{\oplus 3} \rightarrow 0.
\end{equation}

\begin{lemma} \label{lemma:extcompu}
$e(\1, \kappa) = 0$ and $e(\Sp, \kappa) = 3$.
\end{lemma}

\begin{proof}
Applying $\Hom_{G/Z}(\1, -)$ to (\ref{equation:tau}), we obtain the exact sequence
\[
0 \rightarrow \Hom_{G/Z}(\1, \1^{\oplus 3}) \rightarrow \Ext^1_{G/Z}(\1, \Sp) \rightarrow \Ext^1_{G/Z}(\1, \kappa) \xrightarrow{f} \Ext^1_{G/Z}(\1, \1^{\oplus 3}).
\]
Thus to prove the first assertion, it suffices to show that $f$ is identically zero. Suppose not, then there exists a non-split extension of $\1$ by $\kappa$ whose image under $f$ is nonzero, and thus has nonzero image under at least one of the maps 
\[
f_i: \Ext^1_{G/Z}(\1, \kappa) \xrightarrow{f} \Ext^1_{G/Z}(\1, \1^{\oplus 3}) \cong \bigoplus_{i=1}^3 \Ext^1_{G/Z}(\1, \1) \xrightarrow{\pr_i} \Ext^1_{G/Z}(\1, \1)
\]
defined by projecting to $i$-th component. Note that via pullback along $f_i$, such an extension would give rise to a non-split extension of $\1$ by $E_\tau$ (as a subrepresentation), where $E_\tau$ is a non-split extension of $\1$ by $\Sp$ given by some $\tau \in \Hom(\Qp^{\times}, k) \cong \Ext^1_{G/Z}(\1, \Sp)$ defined in \cite[\S \uppercase\expandafter{\romannumeral7}.1]{MR2642409}. This implies that the natural map $\Ext^1_{G/Z}(\1, E_\tau) \rightarrow \Ext^1_{G/Z}(\1, \1)$ is nonzero, which contradicts \cite[Proposition \uppercase\expandafter{\romannumeral7}.5.4]{MR2642409}.

By applying $\Hom_{G/Z}(\Sp, -)$ to (\ref{equation:tau}), we obtain the exact sequence
\[
0 \rightarrow \Ext^1_{G/Z}(\Sp, \Sp) \rightarrow \Ext^1_{G/Z}(\Sp, \kappa) \xrightarrow{g} \Ext^1_{G/Z}(\Sp, \1^{\oplus 3}).
\]
Thus to prove the second assertion, it suffices to show that $g$ is identically zero. Suppose not, then there exists a non-split extension $\kappa'$ of $\Sp$ by $\kappa$ whose image under $f$ is nonzero, and thus has nonzero image under at least one of the maps
\[
g_i: \Ext^1_{G/Z}(\Sp, \kappa) \xrightarrow{g} \Ext^1_{G/Z}(\Sp, \1^{\oplus 3}) \cong \bigoplus_{i=1}^3 \Ext^1_{G/Z}(\Sp, \1) \xrightarrow{\pr_i} \Ext^1_{G/Z}(\Sp, \1)
\]
defined by projecting to $i$-th component. Note that via pullback along $g_i$, such an element would give rise to a non-split extension $\kappa_i$ of $\Ind^G_B \1$ by $\Sp$ (as a subrepresentation of $\kappa'$) by Lemma \ref{lemma:exttristeinberg}. Note that Lemma \ref{lemma:extIndSp} implies that $\Hom_G(\1, \kappa_i) \neq 0$. Hence $\Hom_G(\1, \kappa') \neq 0$, which gives a contradiction since $\Hom_G(\1, \kappa) = \Hom_G(\1, \Sp) = 0$.
\end{proof}

Denote $T_\1 := \iT((\Ind^G_B \1)^{\vee})$, which lies in $\fD(k)$. Note that since $\iT(\1) \cong 0$ in $\fD(k)$ and $\iT$ is exact, we have
\begin{align*}
T_\1 \cong \iT \Sp^{\vee} \cong \iT \tau^{\vee}, \quad \cV(T_\1) \cong \cV(\Sp^\vee) \cong \cV(\tau^{\vee}) \cong \1.
\end{align*}

\begin{lemma} \label{lemma:quotextcompu}
$\Ext^1_{\fD(k)}(T_\1, T_\1)$ is 3-dimensional.
\end{lemma}

\begin{proof}
Replacing \cite[Lemma 10.12]{MR3150248} with Lemma \ref{lemma:extcompu}, the proof of \cite[Lemma 10.34]{MR3150248} works verbatim in our setting. We include the proof for the sake of completeness. Let $J_{\Sp}$ be the injective envelope of $\Sp$ in $\operatorname{Mod}^{\operatorname{l.fin}}_{G/Z}(k)$. It follows from Lemma \ref{lemma:extcompu} that we have an exact sequence:
\begin{equation} \label{equation:injenv}
0 \rightarrow \tau \rightarrow J_{\Sp} \rightarrow J_{\Sp}^{\oplus 3}.
\end{equation}
Moreover, if we let $\theta$ be the cokernel of the second arrow then the monomorphism $\theta \hookrightarrow J_{\Sp}^{\oplus 3}$ induced by the first arrow is essential. We know from Proposition \ref{prop:quotcat} (2) that $\iT J_{\Sp}^\vee$ is the projective envelope of $\Sp^\vee$ in $\fD(k)$. By dualizing (\ref{equation:injenv}), applying $\iT$ and then $\Hom_{\fD(k)}(-, \iT \Sp^\vee)$ we obtain
\[
\Ext^1_{\fD(k)}(T_\1, \iT \Sp^\vee) \cong \Hom_{\fD(k)}(\iT \theta^\vee, \iT \Sp^\vee) \cong \Hom_{\fD(k)}(\iT (J_{\Sp}^{\oplus 3})^{\vee}, \iT \Sp^\vee).
\]
The last isomorphism follows from the fact that $\iT \Sp^\vee$ is irreducible, and $\iT J_{\Sp}^\vee \twoheadrightarrow \iT \theta^\vee$ is essential (Proposition \ref{prop:quotcat} (3)). Hence $\Ext^1_{\fD(k)}(T_\1, T_\1)$ is 3-dimensional.
\end{proof}

\begin{lemma} \label{lemma:injExtV}
The functor $\cV$ induces an injection
\[
\cV: \Ext^1_{\D(\O)}(T_\1, T_\1) \hookrightarrow \Ext^1_{\GQp}(\cV(T_\1), \cV(T_\1)).
\]
\end{lemma}

\begin{proof}
Note that \cite[Proposition \uppercase\expandafter{\romannumeral7}.4.12]{MR2642409} holds when $p=2$. Thus the proof of \cite[Lemma 10.35]{MR3150248} works verbatim in our setting with Lemma 10.34 of loc. cit. replaced by Lemma \ref{lemma:quotextcompu} above.
\end{proof}

\subsection{A finiteness lemma}
\begin{lemma} \label{lemma:HomVinj}
Let $M, N \in \fD(\O)$ be of finite length. Then $\cV$ induces:
\begin{align*}
\Hom_{\fD(\O)}(M, N) &\cong \Hom_{\GQp}\big(\cV(M), \cV(N)\big), \\
\Ext^1_{\fD(\O)}(M, N) &\hookrightarrow \Ext^1_{\GQp}\big(\cV(M), \cV(N)\big).
\end{align*}
\end{lemma}

\begin{proof}
This is proved in \cite[Lemma A1]{MR2667891} for supersingular blocks and \cite[\S 8]{MR3150248} for principal series blocks. So the only remaining case is when $\B =\{\1, \Sp\} \otimes \delta \circ \det$, where $\delta: \Qp^{\times} \rightarrow k^{\times}$ is a smooth character. The argument in P\v{a}sk\={u}nas' proof is by induction on $\ell(M) + \ell(N)$, where $\ell$ denotes the number of irreducible subquotients, and thus reduces the assertion to the case that both $M$ and $N$ are irreducible. Note that in the exceptional case, we may assume that $\delta = 1$ in which case the assertion for $\Hom$ is immediate and the assertion for $\Ext^1$ follows from Lemma \ref{lemma:injExtV}. This proves the lemma.
\end{proof}

Let $\Modpro_{\GQp}(\O)[\B]$ be the full subcategory of $\Modpro_{\GQp}(\O)$ with object $\rho$ such that there exists $M \in \fC(\O)[\B]$ such that $\rho \cong \cV(M)$.

\begin{prop} \label{prop:equivB}
The functor $\cV$ induces an equivalence of categories between $\fD(\O)[\B]$ and $\Modpro_{\GQp}(\O)[\B]$.
\end{prop}

\begin{proof}
This is due to \cite{MR3150248, MR3544298} except the case that $\B = \{\1, \Sp \} \otimes \delta \circ \det$. Note that in the exceptional case, the proof of \cite[Proposition 10.36]{MR3150248} works verbatim with Lemma 10.35 in loc. cit. replaced by Lemma \ref{lemma:injExtV} above. This proves the proposition.
\end{proof}

\begin{prop} \label{prop:fgadm}
If $\pi \in \Modlfinz(k)$ is admissible, then $\cV(\pi^\vee)$ is  finitely generated as a $k \bracketGQp$-module.
\end{prop}

\begin{proof}
This follows from the proof of \cite[Proposition 2.8]{2018arXiv180307451T} with Lemma 2.6 in loc. cit. replaced by Lemma \ref{lemma:HomVinj} above.
\end{proof}
\section{Automorphic forms on $\GL_2(\A_F)$} \label{section:autoform}
We define the class of automorphic representations whose associated Galois representations we wish to study. Throughout this section, we let $F$ be a totally real field and fix an isomorphism $\iota: \bQp \cong \C$.

If $\lambda = (\lambda_{\kappa})_{\kappa: F \rightarrow \C} \in (\Z^2_+)^{\Hom(F, \C)}$, let $\Xi_{\lambda}$ denote the irreducible algebraic representation of $(\GL_2)^{\Hom(F, \C)}$ which is the tensor product over $\kappa \in \Hom(F, \C)$ of irreducible representations of $\GL_2$ with highest weight $\lambda_{\kappa}$. We say that $\lambda \in (\Z^2_+)^{\Hom(F, \C)}$ is an algebraic weight if it satisfies the parity condition, i.e. $\lambda_{\kappa, 1} + \lambda_{\kappa, 2}$ is independent of $\kappa$.

\begin{definition}
We say that a cuspidal automorphic representation $\pi$ of $\GL_2(\A_F)$ is regular algebraic if the infinitesimal character of $\pi_{\infty}$ has the same infinitesimal character as $\Xi_{\lambda}^{\vee}$ for an algebraic weight $\lambda$. 
\end{definition}

Let $\pi$ be a regular algebraic cuspidal automorphic representation of $\GL_2(\A_F)$ of weight $\lambda$. For any place $v | p$ of $F$ and any integer $a \geq 1$, let $\Iw_v(a, a)$ denote the subgroup of $\GL_2(\O_{F_v})$ of matrices that reduce to an upper triangular matrix modulo $\varpi_v^a$. We define the Hecke operator
\[
\U_{\varpi_v} = \bigg[\Iw_v(a, a) \begin{pmatrix} \varpi_v & 0 \\ 0 & 1 \end{pmatrix} \Iw_v(a, a) \bigg]
\]
and the modified Hecke operator
\[
\U_{\lambda, \varpi_v} = \bigg( \prod_{\kappa:F_v \hookrightarrow \bQp}  \kappa(\varpi_v)^{- \lambda_{\iota \kappa, 2}} \bigg) \U_{\varpi_v}.
\]

\begin{definition}
Let $v$ be a place of $F$ above $p$. We say that $\pi$ is $\iota$-ordinary at $v$, if there is an integer $a \geq 1$ and a nonzero vector in $(\iota^{-1} \pi_v)^{\Iw_v(a, a)}$ that is an eigenvector for $\U_{\lambda, \varpi_v}$ with an eigenvalue which is a $p$-adic unit. This definition does not depend on the choice of $\varpi_v$.
\end{definition}

The following theorem is due to the work of many people. We refer the reader to \cite{MR870690} and \cite{MR1016264} for the existence of Galois representations, to \cite{MR870690} for part (2) when $v \nmid p$, to \cite{MR2551990} for part (1) and part (2) when $v | p$, and to \cite{MR1463699, MR969243} for part (3).

\begin{thm}
Let $\pi$ be a regular algebraic cuspidal automorphic representation of $\GL_2(\A_F)$ of weight $\lambda$. Fix an isomorphism $\iota: \bQp \rightarrow \C$. Then there exists a continuous semi-simple representation
\[
\rho_{\pi, \iota}: G_F \rightarrow \GL_2(\bQp)
\]
satisfying the following conditions:
\begin{enumerate}
\item For each place $v | p$ of $F$, $\rho_{\pi, \iota} |_{G_{F_v}}$ is de Rham, and for each embedding $\kappa: F \rightarrow \bQp$, we have
\[
\HT_{\kappa}(\rho_{\pi, \iota} |_{G_{F_v}}) = \{\lambda_{\iota \kappa, 2}, \lambda_{\iota \kappa, 1}+1\}.
\]
\item For each finite place $v$ of $F$, we have $\WD(\rho_{\pi, \iota}|_{G_{F_v}})^{F-ss} \cong r_p(\iota^{-1} \pi_v)$.
\item If $\pi$ is $\iota$-ordinary at $v|p$, then there is an isomorphism
\[
\rho |_{G_{F_v}} \sim \begin{pmatrix} \psi_{v, 1} & * \\ 0 & \psi_{v, 2} \end{pmatrix},
\]
where for $i=1, 2$, $\psi_{v, i}: G_{F_v} \rightarrow \bQp^{\times}$ is a continuous character satisfying
\[
\psi_{v, i}(\Art_{F_v}(\sigma)) = \prod_{\kappa: F_v \hookrightarrow \bQp} \kappa(\sigma)^{-(\lambda_{\iota \kappa, 3-i} + i - 1)}
\]
for all $\sigma$ in some open subgroup of $\O_{F_v}^{\times}$.
\end{enumerate}

These conditions characterize $\rho_{\pi, \iota}$ uniquely up to isomorphism.
\end{thm}

\begin{definition}
We call a Galois representation $\rho: G_F \rightarrow \GL_2(\bQp)$ automorphic of weight $\iota^* \lambda = (\lambda_{\iota^{-1} \kappa, 1},  \lambda_{\iota^{-1} \kappa, 2}) \in (\Z^2_+)^{\Hom(F, \bQp)}$ if there exists a regular algebraic cuspidal automorphic representation of $\GL_2(\A_F)$ of weight $\lambda := (\lambda_{\kappa, 1},  \lambda_{\kappa, 2}) \in (\Z^2_+)^{\Hom(F, \C)}$ such that $\rho \cong \rho_{\pi, \iota}$. Moreover, if $\pi$ is $\iota$-ordinary at a place $v|p$ then we say $\rho$ is $\iota$-ordinary at $v$.
\end{definition}
\section{Galois deformation theory} \label{section:Galoisdeform}
\subsection{Global deformation problems}
Let $F$ be a number field and $p$ be a prime. We fix a continuous absolutely irreducible $\rhobar : G_F \rightarrow \GL_2(k)$ and a continuous character $\psi: G_F \rightarrow \O^{\times}$ such that $\chi \varepsilon$ lifts $\det \rhobar$. We fix a finite set $S$ of places of $F$ containing those above $p, \infty$ and the places at which $\rhobar$ and $\psi$ are ramified. For each $v \in S$, we fix a ring $\Lambda_v \in \CNL_\O$ and define $\Lambda_S = \ctimes_{v \in S, \O} \Lambda_v \in \CNL_{\O}$. 

For each $v \in S$, we denote $\rhobar|_{G_{F_v}}$ by $\rhobar_v$ and write $\Df_v : \CNL_{\Lambda_v} \rightarrow \Set$ (resp. $\Dfz_v : \CNL_{\Lambda_v} \rightarrow \Set$) for the functor associates $R \in \CNL_{\Lambda_v}$ the set of all continuous homomorphisms $r: G_{F_v} \rightarrow \GL_2(R)$ such that $r$ mod $\m_R = \rhobar_v$ (resp. and $\det r$ agrees with the composition $G_{F_v} \rightarrow \O^{\times} \rightarrow R^{\times}$ given by $\psi \varepsilon |_{G_{F_v}}$), which is represented by an object $\Rf_v \in \CNL_{\Lambda_v}$ (resp. $\Rfz_v \in \CNL_{\Lambda_v}$). We will write $\rhof_v: G_{F_v} \rightarrow \GL_2(\Rf_v)$ for the universal lifting of $\rhobar_v$.

\begin{definition}
Let $v \in S$, a local deformation problem for $\rhobar_v$ is a subfunctor $\D_v \subset \Df_v$ satisfying the following conditions:
\begin{itemize}
\item $\D_v$ is represented by a quotient $R_v$ of $\Rf_v$.
\item For all $R \in \CNL_{\Lambda_v}$, $a \in \Ker(\GL_2(R) \rightarrow \GL_2(k))$ and $r \in \D_v(R)$, we have $a r a^{-1} \in \D_v(R)$.
\end{itemize}
\end{definition}

\begin{definition}
A global deformation problem is a tuple
\[
\cS = (\rhobar, S, \{\Lambda_v \}_{v \in S}, \{\D_v \}_{v \in S})
\]
where
\begin{itemize}
\item the object $\rhobar$, $S$ and $\{\Lambda_v \}_{v \in S}$ are defined as above.
\item for each $v \in S$, $\D_v$ is a local deformation problem for $\rhobar_v$.
\end{itemize}
\end{definition}

\begin{definition}
Let $\cS = (\rhobar, S, \{\Lambda_v \}_{v \in S}, \{\D_v \}_{v \in S})$ be a global deformation problem. Let $R \in \CNL_{\Lambda_S}$, and let $\rho: G_F \rightarrow \GL_2(R)$ be a lifting of $\rhobar$. We say that $\rho$ is of type $\cS$ if it satisfies the following conditions:
\begin{enumerate}
\item $\rho$ is unramified outside $S$.
\item For each $v \in S$, $\rho_v := \rho |_{G_{F_v}}$ is in $\D_v(R)$, where $R$ has a natural $\Lambda_v$-algebra structure via the homomorphism $\Lambda_v \rightarrow \Lambda_S$.
\end{enumerate}
\end{definition}

We say that two liftings $\rho_1, \rho_2: G_F \rightarrow \GL_2(R)$ are strictly equivalent if there exists $a \in \Ker(\GL_2(R) \rightarrow \GL_2(k))$ such that $\rho_2 = a \rho_1 a^{-1}$. It's easy to see that strictly equivalence preserves the property of being type $\cS$. 

We write $\Df_\cS$ for the functor $\CNL_{\Lambda_S} \rightarrow \Set$ which associates to $R \in \CNL_{\Lambda_S}$ the set of liftings $\rho: G_F \rightarrow \GL_2(R)$ which are of type $\cS$, and write $\D_{\cS}$ for the functor $\CNL_{\Lambda_S} \rightarrow \Set$ which associates to $R \in \CNL_{\Lambda_S}$ the set of strictly equivalence classes of liftings of type $\cS$.

\begin{definition}
If $T \subset S$ and $R \in \CNL_{\Lambda_S}$, then a $T$-framed lifting of $\rhobar$ to $R$ is a tuple $(\rho, \{\alpha_v \}_{v \in T})$, where $\rho$ is a lifting of $\rhobar$, and for each $v \in T$, $\alpha_v$ is an element of $\Ker(\GL_2(R) \rightarrow \GL_2(k))$. Two $T$-framed liftings $(\rho, \{\alpha_v \}_{v \in T})$ and $(\rho', \{\alpha_v' \}_{v \in T})$ are strictly equivalent if there is an element $a \in \Ker(\GL_2(R) \rightarrow \GL_2(k))$ such that $\rho' = a \rho a^{-1}$ and $\alpha_v' = a \alpha_v$ for each $v \in T$.
\end{definition}

We write $\D^T_\cS$ for the functor $\CNL_{\Lambda_S} \rightarrow \Set$ which associates to $R \in \CNL_{\Lambda_S}$ the set of strictly equivalence classes of $T$-framed liftings $(\rho, \{\alpha_v\}_{v \in T})$ to $R$ such that $\rho$ is of type $\cS$. Similarly, we may consider liftings of type $\cS$ with determinant $\psi \varepsilon$, and we denote the corresponding functor by $\Dz_\cS$, $\Dfz_\cS$ and $\D^{T, \psi}_\cS$.

\begin{thm}
Let $\cS = (\rhobar, S, \{\Lambda_v \}_{v \in S}, \{\D_v \}_{v \in S})$ be a global deformation problem. Then the functor $\D_\cS$, $\Df_\cS$, $\D^T_\cS$, $\Dz_\cS$, $\Dfz_S$ and $\D^{T, \psi}_S$ are represented by objects $R_\cS$, $\Rf_\cS$, $R^T_\cS$, $\Rz_\cS$, $\Rfz_\cS$ and $R^{T, \psi}_\cS$, respectively, of $\CNL_{\Lambda_S}$.
\end{thm}

\begin{proof}
For $\D_\cS$, this is due to \cite[Theorem 9.1]{MR1860043}. The representability of the functors $\Df_\cS$, $\D^T_\cS$, $\Dz_\cS$, $\Dfz_S$ and $\D^{T, \psi}_S$ can be deduced easily from this.
\end{proof}

\begin{lemma} \label{lemma:Galoisframe}
Let $\cS$ be a global deformation problem. Choose $v_0 \in T$, and let $\cT = \O \llbracket X_{v, i, j} \rrbracket_{v \in T, 1 \leq i, j \leq 2} / (X_{v_0, 1, 1})$. There is a canonical isomorphism $R^T_{\cS} \cong R_\cS \ctimes_\O \cT$.
\end{lemma}

\begin{proof}
Let $\rho_{\cS}: G_F \rightarrow \GL_2(R_{\cS})$ be a universal solution of deformations of type $\cS$. Note that the centralizer in $id_2 + M_2(\m_{R_{\cS}})$ of $\rho_S$ is the scalar matrices, Thus the $T$-framed lifting over $R_\cS \ctimes_{\O} \cT$ given by the tuple $(\rho_S, \{id_2 + (X_{v, i, j})\}_{v \in T})$ is a universal framed deformation of $\rbar$ over $R_\cS \ctimes_{\O} \cT$. This shows that the induced map $R^T_\cS \rightarrow R_\cS \ctimes_{\O} \cT$ is an isomorphism.
\end{proof}

Let $\cS = (\rhobar, S, \{\Lambda_v \}_{v \in S}, \{\D_v \}_{v \in S})$ be a global deformation problem and denote $R_v \in \CNL_{\Lambda_v}$ the representing object of $\D_v$ for each $v \in S$. We write $A^T_\cS = \ctimes_{v \in T, \O} R_v$ for the completed tensor product of $R_v$ over $\O$ for each $v \in T$, which has a canonical $\Lambda_T := \ctimes_{v \in T, \O} \Lambda_v$ algebra structure. The natural transformation $(\rho, \{\alpha_v \}_{v \in T}) \mapsto (\alpha_v^{-1} \rho |_{G_{F_v}} \alpha_v)_{v\in T}$ induces a canonical homomorphism of $\Lambda_T$-algebras $A^{T}_\cS \rightarrow R^T_{\cS}$. Moreover, Lemma \ref{lemma:Galoisframe} allows us to consider $R_{\cS}$  as an $A^{T}_\cS$-algebra via the map $A^T_\cS \rightarrow  R^T_{\cS} \twoheadrightarrow R_{\cS}$.

\begin{prop} \label{prop:finitedefbc}
Let $\cS$ be a global deformation problem as before and $F'$ be a finite Galois extension of $F$. Suppose that 
\begin{itemize}
\item $\End_{G_{F'}}(\rhobar) = k$.
\item $\cS' = (\rhobar |_{G_{F'}}, S', \{\Lambda_w \}_{v \in S'}, \{\D_w \}_{w \in S'})$ is a deformation problem where 
\begin{itemize}
\item $S'$ is the set of places of $F'$ above $S$;
\item $T'$ is the set of places of $F'$ above $T$;
\item for each $w | v$, $\Lambda_w = \Lambda_v$ and $\D_w$ is a local deformation problem equipped with a natural map $R_w \rightarrow R_v$ induced by restricting deformations of $\rhobar_v$ to $G_{F'_w}$.
\end{itemize}
\end{itemize}
Then the natural map $R^{T', \psi}_{\cS'} \rightarrow R^{T, \psi}_\cS$ induced by restricting deformations of $\rhobar$ to $G_{F'}$, make $R^{T, \psi}_\cS$ into a finitely generated $R^{T', \psi}_{\cS'}$-module.
\end{prop}

\begin{proof}
Let $\m'$ be the maximal ideal of $R^{T', \psi}_{\cS'}$. It follows from \cite[Lemma 3.6]{MR2480604} and Nakayama's lemma that it is enough to show the image of $G_{F, S} \rightarrow \GL_2(R^{T, \psi}_\cS) \rightarrow \GL_2(R^{T, \psi}_\cS / \m' R^{T, \psi}_\cS)$ is finite. Since $G_{F', S'}$ is of finite index in $G_{F, S}$ and it gets mapped to the finite subgroup $\rhobar(G_{F', S'})$, we are done.
\end{proof}

\subsection{Local deformation problems}
In this section, we define some local deformation problems we will use later.

\subsubsection{Ordinary deformations}
\label{section:orddeform}
We define ordinary deformations following \cite[\S 1.4]{MR3252020}.

Suppose that $v | p$ and that $E$ contains the image of all embeddings $F_v \hookrightarrow \bQp$. We will assume throughout this subsection that there is some line $\overline{L}$ in $\rhobar_v$ that is stable by the action of $G_{F_v}$. Let $\overline{\eta}$ denote the character of $G_{F_v}$ giving the action on $\overline{L}$. Note that the choice of $\overline{\eta}$ is unique unless $\rhobar_v$ is the direct sum of two distinct characters. In this case we simply make a choice of one of these characters.

We write $\O_{F_v}^{\times}(p)$ for the maximal pro-$p$ quotient of $\O_{F_v}^{\times}$. Set $\Lambda_v = \O \llbracket \O^{\times}_{F_v}(p) \rrbracket$ and write $\psi^{\univ}: G_{F_v} \rightarrow \Lambda_v^{\times}$ for the universal character lifting $\overline{\psi}$. Note that $\Art_{F_v}$ restricts to an isomorphism $\O^{\times}_{F_v} \cong I^{\ab}_{F_v}$, where $I^{\ab}_{F_v}$ is the inertial subgroup of the maximal abelian extension of $F_v$.

Let $\P^1$ be the projective line over $\O$. We denote $\L_{\Delta}$ the subfunctor of $\P^1 \times_{\O} \Spec \Rfz_v$, whose $A$-points for any $\O$-algebra $A$ consist of an $\O$-algebra homomorphism $\Rfz_v \rightarrow A$ and a line $L \in \P^1(A)$ such that the filtration is preserved by the action of $G_{F_v}$ on $A^2$ induced from $\rho_v^{\Box}$ and such that the action of $G_{F_v}$ on $L$ is given by pushing forward $\psi^{\univ}$. This subfunctor is represented by a closed subscheme (c.f. \cite[Lemma 1.4.2]{MR3252020}), which we denote by $\L_{\Delta}$ also. We define $\Rord_v$ to be the maximal reduced, $\O$-torsion free quotient of the image of the map $\Rfz_v \rightarrow H^0(\L_{\Delta}, \O_{\L_{\Delta}})$.

\begin{prop} \label{prop:orddeform}
The ring $\Rord_v$ defines a local deformation problem. Moreover,
\begin{enumerate}
\item An $\O$-algebra homomorphism $x: \Rfz_v \rightarrow \bQp$ factors through $\Rord_v$ if and only if the corresponding Galois representation is $\GL_2(\bQp)$-conjugate to a representation
\[
\begin{pmatrix} \psi_1 & * \\ 0 & \psi_2 \end{pmatrix}
\]
where $\psi_1 |_{G_{F_v}} = x \circ \psi^{\univ}$.
\item Assume the image of $\bar{\rho} |_{G_{F_v}}$ is either trivial or has order $p$, and that if $p=2$, then either $F_v$ contains a primitive fourth roots of unity or $[F_v : \Q_2] \geq 3$. Then for each minimal prime $Q_v \subset \Lambda_v$, $\Rord_v / Q_v$ is an integral domain of relative dimension $3+ 2 [F_v : \Qp]$ over $\O$, and its generic point is of characteristic 0.
\end{enumerate}
\end{prop}

\begin{proof}
The first assertion follows from \cite[Proposition 1.4.4]{MR3252020} and the second assertion is due to \cite[Proposition 1.4.12]{MR3252020}.
\end{proof}

We define $\Dord_v$ to be the local deformation problem represented by $\Rord_v$.

\subsubsection{Potentially semi-stable deformations}
Suppose that $v | p$ and that $E$ contains the image of all embeddings $F_v \hookrightarrow \bQp$. Let $\Lambda_v = \O$.

\begin{prop}
For each $\lambda_v \in (\Z^2_+)^{\Hom(F_v, E)}$ and inertial type $\tau_v : I_v \rightarrow \GL_2(E)$, there is a unique (possibly trivial) quotient $\Rss_v$ (resp. $\Rcr_v$) of the universal lifting ring $\Rfz_v$ with the following properties:
\begin{enumerate}
\item $\Rss_v$ (resp. $\Rcr_v$) is reduced and $p$-torsion free, and all the irreducible components of $\Rss_v[1/p]$ (resp. $\Rcr_v[1/p]$) are formally smooth and of relative dimension $3+[F_v:\Qp]$ over $\O$.
\item If $E' / E$ is a finite extension, then an $\O$-algebra homomorphism $\Rfz_v \rightarrow E'$ factors through $\Rss_v$ (resp. $\Rcr_v$) if and only if the corresponding Galois representation $G_{F_v} \rightarrow \GL_2(E')$ is potentially semi-stable (resp. potentially crystalline) of weight $\lambda$ and inertial type $\tau$.
\item $\Rss_v / \varpi$ (resp. $\Rcr_v / \varpi$) is equidimensional.	
\end{enumerate}
\end{prop}

\begin{proof}
This is due to \cite{MR2373358} (see also \cite[Corollary 1.3.5]{Allen2014DeformationsOH}).
\end{proof}

In the case that $\Rss_v \neq 0$ (resp. $\Rcr_v \neq 0$), we define $\Dss_v$ (resp. $\Dcr_v$) to be the local deformation problem represented by $\Rss_v$ (resp. $\Rcr_v$).

\subsubsection{Fixed weight potentially semi-stable deformations}

For $\lambda_v \in (\Z^2_+)^{\Hom(F_v, E)}$, we define characters $\psi^{\lambda_v}_i: I_{F_v} \rightarrow \O^{\times}$ for $i = 1, 2$ by
\[
\psi^{\lambda_v}_i: \sigma \mapsto \varepsilon(\sigma)^{-(i-1)} \prod_{\kappa_v: F_v \hookrightarrow E} \kappa_v(\Art_{F_v}^{-1}(\sigma))^{-\lambda_{\kappa_v, 3 - i}}.
\]

\begin{definition} \label{def:ordinary}
Let $\lambda_v \in (\Z^2_+)^{\Hom(F_v, E)}$ and $\rho_v: G_{F_v} \rightarrow \GL_2(\O)$ be a continuous representation. We say $\rho$ is ordinary of weight $\lambda_v$ if there is an isomorphism
\[
\rho_v \sim \begin{pmatrix} \psi_{v, 1} & * \\ 0 & \psi_{v, 2} \end{pmatrix},
\]
where for $i=1, 2$, $\psi_{v, i}: G_{F_v} \rightarrow \O^{\times}$ is a continuous character agrees with $\psi^{\lambda_v}_i$ on an open subgroup of $I_{F_v}$.
\end{definition}

\begin{prop} \label{prop:Rpst}
For each $\lambda_v$, $\tau_v$ there is a unique (possibly trivial) reduced and $p$-torsion free quotient $\Rssord_v$ of $\Rord_v$ satisfying the following properties:
\begin{enumerate}
\item If $E' / E$ is a finite extension, then the $\O$-algebra homomorphism $\Rfz_v \rightarrow E'$ factors through $\Rssord_v$ if and only if the corresponding Galois representation $G_v \rightarrow \GL_2(E')$ is ordinary and potentially semi-stable of Hodge type $\lambda$ and inertial type $\eta$.
\item $\Spec \Rssord_v$ is a union of irreducible components of $\Spec R^{\lambda_v, \tau_v}_v$.	
\end{enumerate}
\end{prop}

\begin{proof}
This follows from \cite[Lemma 3.3.3]{MR2941425}.
\end{proof}

\begin{lemma} \label{lemma:ordchar}
If $\Rssord_v$ is non-zero, then $\tau = \alpha_1 \oplus \alpha_2$ is a sum of smooth characters of $I_v$. Moreover, the natural surjection $\Rord_v \twoheadrightarrow \Rssord_v$ factors through $\Rord_v \otimes_{\O \llbracket \O^{\times}_{F_v}(p) \rrbracket, \eta} \O$, where $\eta: \O \llbracket \O^{\times}_{F_v}(p) \rrbracket \rightarrow \O$ is given by $u \mapsto \alpha_1(\Art_{F_v}(u)) \prod_{\kappa_v: F_v \hookrightarrow E} \kappa_v(\Art_{F_v}^{-1}(\sigma))^{-\lambda_{\kappa_v, 2}}$ for $u \in \O^{\times}_{F_v}(p)$.
\end{lemma}

\begin{proof}
The first assertion is due to \cite[Lemma 3.3.2]{MR2941425}. For the second assertion, consider the following diagram
\[
  \begin{tikzcd}
  &\L_{\lambda_v, \tau_v}  \arrow[r] \arrow[d] &\L \arrow[r] \arrow[d] & \L_{\Delta} \arrow[d]\\
  &\Spec \Rss_v  \arrow[r, hookrightarrow] &\Spec \Rfz_v \arrow[r, hookrightarrow] &\Spec \tRfz_v,
  \end{tikzcd}
\]
where $\tRfz_v = \Rfz_v \ctimes_{\O} \O \llbracket \O^{\times}_{F_v}(p) \rrbracket$, $\Spec \Rfz_v \hookrightarrow \Spec \tRfz_v$ is induced by the surjection $\tRfz_v \twoheadrightarrow \Rfz_v$ given by $\eta$, $\L_{\lambda_v, \tau_v}$ is the closed subscheme of $\P^1 \times_{\O} \Spec \Rss_v$, whose $R$-valued points, $R$ an $\Rss_v$-algebra, consist of a $R$-line $L \subset R^2$ on which $I_{F_v}$ acts via the character $\eta$ composed with $\Art_{F_v}$, and $\L$ is the closed subscheme of $\P^1 \times_{\O} \Spec \Rfz_v$ defined in the same way using $\Rfz_v$ instead of $\Rss_v$.

It's easy to see that the left square (induced by the quotient $\Rfz_v \twoheadrightarrow \Rss_v$) is cartesian and the right square is commutative. This proves the proposition since $\Rord_v$ is the scheme theoretical image of $\L$ in $\tRfz_v$ and $\Rssord_v$ is the scheme theoretical image of $\L_{\lambda_v, \tau_v}$ in $\Spec \Rss_v$ (c.f. \cite[\S 3.3]{MR2941425}).
\end{proof}

\subsubsection{Irreducible components of potentially semi-stable deformations}
Suppose that $\cC_v$ is an irreducible component of $\Spec \Rss_v[1/p]$. Then we write $R^{\cC_v}_v$ for the maximal reduced, $p$-torsion free quotient of $\Rss_v$ such that $\Spec R^{\cC_v}_v [1/p]$ is the component $\cC_v$.

\begin{lemma}
Say that a lifting $\rho: G_{F_v} \rightarrow \GL_2(R)$ is of type $\D^{\cC_v}_v$ if the induced map $\Rfz_v \rightarrow R$ factors through $R^{\cC_v}_v$. Then $\D^{\cC_v}_v$ is a local deformation problem.
\end{lemma}

\begin{proof}
This follows from \cite[Lemma 1.2.2]{MR3152941} and \cite[Lemma 3.2]{MR2827723}.
\end{proof}

We say that an irreducible component $\cC_v$ of $\Spec \Rss_v$ is ordinary if it lies in the support of $\Spec \Rssord_v$, and non-ordinary otherwise.

\subsubsection{Odd deformations}
Assume that $F_v = \R$ and $\rhobar |_{G_{F_v}}$ is odd, i.e. $\det \rhobar(c) = -1$ for $c$ the complex conjugation. Let $\Lambda_v = \O$.

\begin{prop} \label{prop:Rodd}
There is a reduced and $p$-torsion free quotient $\Rodd_v$ of $\Rfz_v$ such that if $E' / E$ is a finite extension, a $\O$-homomorphism $\Rfz_v \rightarrow E'$ factors through $\Rodd_v$ if and only if the corresponding Galois representation is odd. Moreover,
\begin{itemize}
\item $\Rodd_v$ is a complete intersection domain of relative dimension $2$ over $\O$.
\item $\Rodd_v[1/p]$ is formally smooth over $E$.
\item $\Rodd_v \otimes_{\O} k$ is a domain.
\end{itemize}
\end{prop}

\begin{proof}
See \cite[Proposition 3.3]{MR2551764}.
\end{proof}

We write $\Dodd_v$ for the local deformation problem defined by $\Rodd_v$.

\subsubsection{Irreducible components of unrestricted deformations}
Let $v \nmid p$ and $\Lambda_v = \O$.

\begin{lemma} \label{lemma:ltype}
Let $x,y: \Rfz_v \rightarrow \bQp$ with $\rho_x, \rho_y: G_{F_v} \rightarrow \GL_2(\bQp)$ be the associated framed deformations.
\begin{enumerate}
\item If $x$ and $y$ lie on the same irreducible component of $\Spec \Rfz_v \otimes \bQp$, then
\[
(\rho_x) |^{ss}_{I_{F_v}} \cong (\rho_y) |^{ss}_{I_{F_v}}.
\]
\item Suppose that moreover neither $x$ nor $y$ lie on any other irreducible component of $\Spec \Rfz_v \otimes \bQp$. Then
\[
(\rho_x) |_{I_{F_v}} \cong (\rho_y) |_{I_{F_v}}.
\]
\end{enumerate}
\end{lemma}

\begin{proof}
See \cite[Lemma 1.3.4]{MR3152941}.
\end{proof}

Suppose that $\cC_v$ is an irreducible component of $\Spec \Rfz_v[1/p]$. Then we write $R^{\cC_{v}}_v$ for the maximal reduced, $p$-torsion free quotient of $\Rfz_v$ such that $\Spec R^{\cC_{v}}_v[1/p]$ is supported on the component $\cC_v$, which defines a local deformation problem $\D^{\cC_v}_v$ by \cite[Lemma 3.2]{MR2827723}. Moreover, it follows from Lemma \ref{lemma:ltype} that all points of $\Spec R^{\cC_{v}}_v[1/p]$ are of the same inertial type if $E$ is large enough.

\subsubsection{Unramified deformations}
Let $v \nmid p$ and $\Lambda_v = \O$.

\begin{prop}
Suppose $\rhobar|_{G_{F_v}}$ is unramified and $\psi$ is unramified at $v$. There there is a reduecd, $\O$-torsion free quotient $\Rur_v$ of $\Rfz_v$ corresponding to unramified deformations. Moreover, $\Rur_v$ is formally smooth over $\O$ of relative dimension 3.
\end{prop}

\begin{proof}
This is due to \cite[prop 2.5.3]{MR2551765}.
\end{proof}

We denote $\Dur_v$ the local deformation problem defined by $\Rur_v$.

\subsubsection{Special deformations}
Let $v \nmid p$ and $\Lambda_v = \O$.

\begin{prop} \label{prop:Rst}
There is a reduced, $\O$-torsion free quotient $\RSt_v$ of $\Rfz_v$ satisfying the following properties: 
\begin{enumerate}
\item If $E'/E$ is a finite extension then an $\O$-algebra homomorphism $\Rfz_v \rightarrow E'$ factors through $\RSt_v$ if and only if the corresponding Galois representation is an extension of $\gamma_v$ by $\gamma_v(1)$, where $\gamma_v: G_{F_v} \rightarrow \O^{\times}$ is an unramified character such that $\gamma_v^2 = \psi|_{G_{F_v}}$.  
\item $\RSt_v$ is a domain of relative dimension 3 over $\O$ and $\RSt_v[1/p]$ is regular.
\end{enumerate}
\end{prop}

\begin{proof}
This follows from \cite[Proposition 2.6.6]{MR2600871} and \cite[Theorem 3.1]{MR2551764}.
\end{proof}

We denote $\DSt_v$ the local deformation problem defined by $\RSt_v$.

\subsubsection{Taylor-Wiles deformations} \label{section:TWdef}
Suppose that $q_v \equiv 1$ mod $p$, that $\rhobar |_{G_{F_v}}$ is unramified, and that $\rhobar(\Frob_v)$ has distinct eigenvalues $\alpha_{v, 1}, \alpha_{v, 2} \in k$. Let $\Delta_v = k(v)^{\times}(p)$ be the maximal $p$-power order quotient of $k(v)^{\times}$ and $\Lambda_v = \O [\Delta_v^{\oplus 2}]$.

\begin{prop} \label{prop:TWdef}
$\Rf_v$ is a formally smooth $\Lambda_v$-algebra. Moreover, $\rho^{\Box}_v \cong \chi_{v, 1} \oplus \chi_{v, 2}$ with $\chi_{v, i}$ a character satisfying $\chi_{v, i}(\Frob_v) \equiv \alpha_{v, i}$ mod $\m_{\Rf_v}$ and $\chi_{v, i} |_{I_{F_v}}$ agrees, after the composition with the Artin map, with the character $k(v)^{\times} \rightarrow \Delta_v^{\oplus 2} \rightarrow \Lambda_v^{\times}$ defined by mapping $k(v)^{\times}$ to its image in the $i$-th component of $\Delta_v$.
\end{prop}

\begin{proof}
This follows from the proof of \cite[Lemma 2.44]{MR1474977} (see \cite[Proposition 5.3]{MR3554238} for an explicit computation of $\Rf_v$).
\end{proof}

In this case, we write $\DTW_v$ for $\Df_v$.

\subsection{Irreducible component of $p$-adic framed deformation rings of $G_{\Q_2}$} \label{section:irredcomp}
Assume $p=2$. Let $\rbar: G_{\Qp} \rightarrow \GL_2(k)$ and $\zeta: G_{\Qp} \rightarrow \O^{\times}$ be a lifting of $\det \rbar \varepsilon^{-1}$. We write $R_{\rbar}$ (resp. $\Rzp$) for the universal lifting ring of $\rbar$ (resp. universal lifting ring of $\rbar$ with determinant $\zeta \varepsilon$). Denote $R_{\overline{\zeta}}$ the universal deformation ring of $\overline{\zeta} = \det \rbar$ (note that $\overline{\varepsilon} = \1$).

\begin{thm} \label{thm:Rzpirred}
The morphism $\Spec R_{\rbar} \rightarrow \Spec R_{\overline{\zeta}}$ given by mapping a deformation of $\rbar$ to its determinant induces a bijection between the irreducible components of $\Spec R_{\rbar}$ and those of $\Spec R_{\overline{\zeta}}$.
\end{thm}

\begin{rmk}
When $p > 2$ and $\rbar: G_L \rightarrow \GL_2(k)$ with $L$ an arbitrary finite extension of $\Qp$, the theorem is proved in \cite[Theorem 1.9]{MR3406170}.
\end{rmk}

\begin{proof}
This is proved in \cite[Proposition 4.1]{Chenevier2009varcarapadique} when $\rbar$ absolutely irreducible or reducible indecomposable with non-scalar semi-simplification. Assume that $\rbar$ is split reducible with non-scalar semi-simplification (i.e. $\rbar \cong (\begin{smallmatrix}
\bar{\chi}_1 & 0 \\ 0 & \bar{\chi_2}
\end{smallmatrix})$ with $\bar{\chi}_1 \bar{\chi}_2^{-1} \neq \1$). It is proved in \cite[Proposition 5.2]{MR3671561} that $R^{\ver} \cong R^{\ps} \llbracket x, y \rrbracket / (xy -c)$, where $R^{\ver}$ is the versal deformation ring of $\rbar$, $R^{\ps}$ is the pseudo deformation ring of (the pseudo-character associated to) $\rbar$, and $c \in R^{\ps}$ is the element generating the reducibility ideal. Since $R^{\ps}$ is isomorphic to the universal deformation ring of $\rbar' = (\begin{smallmatrix}
\bar{\chi}_1 & * \\ 0 & \bar{\chi_2}
\end{smallmatrix})$ with $* \neq 0$ by \cite[Proposition 3.6]{MR3671561} and $xy-c$ is irreducible in $R^{\ps} \llbracket x, y \rrbracket$, it follows that the irreducible components of $\Spec R^{\ver}$ are in bijection with the irreducible components of $\Spec R_{\overline{\zeta}}$. This implies the theorem since $R_{\rbar}$ is formally smooth over $R^{\ver}$ \cite[Proposition 2.1]{MR2551764}. For $\rbar$ reducible with scalar semi-simplification, this is due to \cite[Theorem 9.4]{MR3384443} when $\rbar$ is split and \cite[Satz 5.4]{2015arXiv151209277B} when $\rbar$ is non-split.
\end{proof}

We will write $R_\1$ for the universal deformation ring of the trivial character $\1: G_{\Q_2} \rightarrow k^{\times}$ and $\1^{\univ}: G_{\Q_2} \rightarrow R_{\1}^{\times}$ for its universal deformation. Note that the map $\zeta \mapsto \zeta \chi$ with $\chi$ any lifting of $1$ induces an isomorphism $R_{\overline{\zeta}} \cong R_\1 \cong \O \llbracket x,y,z \rrbracket / ((1+z)^2 - 1)$, which has two irreducible components determined by $\zeta(\Art_{\Q_2}(-1)) \in  \{ \pm 1\}$. It follows that two points $x$ and $y$ of $\Spec R_{\rbar}$ lie in the same irreducible component if and only if the associated liftings $r_x$ and $r_y$ satisfying $\det r_x (\Art_{\Q_2}(-1) = \det r_y (\Art_{\Q_2}(-1))$. We denote $R^{\sign}_{\rbar}$ the the complete local noetherian $\O$-algebra pro-represents the functor sending $R \in \CNL_{\O}$ to the set of liftings $r$ of $\rbar$ to $R$ such that $\det r(\Art_{\Q_2}(-1)) = \zeta(\Art_{\Q_2}(-1))$. Thus $\Spec R^{\sign}_{\rbar}$ is an irreducible component of $\Spec R_{\rbar}$. 

\begin{corollary} \label{cor:Rsign}
$R^{\sign}_{\rbar}[\frac{1}{2}]$ is an integral domain.
\end{corollary}

\begin{proof}
If $\rbar$ absolutely irreducible or reducible indecomposable with non-scalar semi-simplification, it can be shown that $R^{\sign}_{\rbar} \cong \O \llbracket X_1, \cdots, X_5 \rrbracket$ using \cite[Proposition 4.1]{Chenevier2009varcarapadique}. The assertion for $\rbar$ split reducible non-scalar follows from the non-split case by the same arguments in the proof of Theorem \ref{thm:Rzpirred}. For $\rbar$ reducible with scalar semi-simplification, it is proved in \cite[Theorem 9.4]{MR3384443} when $\rbar$ is split and \cite[Satz 5.4]{2015arXiv151209277B} when $\rbar$ is non-split that $R^{\sign}_{\rbar}[1/2]$ is an integral domain. 
\end{proof}

\begin{prop} \label{prop:localetale}
The morphism $\Spec (\Rzp \ctimes_{\O} R_\1) \rightarrow \Spec R^{\sign}_{\rbar}$ induced by $(r, \chi) \mapsto r \otimes \chi$ is finite and becomes \'{e}tale after inverting $2$.
\end{prop}

\begin{proof}
Following the proof of \cite[Proposition 1.1.11]{Allen2014DeformationsOH}, we consider the following cartesian product
\[
  \begin{tikzcd}
   &\Spec R^{\sign}_{\rbar} \times_{\Spec R_\1} \Spec R_\1 \arrow[r] \arrow[d] &\Spec R_\1 \arrow[d, "s"] \\
   &\Spec R^{\sign}_{\rbar} \arrow[r, "\delta"] &\Spec R_\1,
  \end{tikzcd}
\]
where $s$ is given by the functor representing $\chi \mapsto \chi^2$ and $\delta$ is given by the functor representing $r \mapsto (\zeta \varepsilon)^{-1} \det r$. It follows that the points of $\Spec R^{\sign}_{\rbar} \times_{\Spec R_\1} \Spec R_\1$ are given by pairs $(r, \chi)$ with $r$ a framed deformation of $\rbar$ and $\chi$ a lifting of $1$ satisfying $\det r = \zeta \varepsilon \chi^2$. Thus the map $(r, \chi) \mapsto (r \otimes \chi^{-1}, \chi)$ induces an isomorphism $\Spec R^{\sign}_{\rbar} \times_{\Spec R_\1} \Spec R_\1 \cong \Spec (\Rzp \ctimes_{\O} R_\1)$. Note that the morphism $s$ is given by $x \mapsto (1+x)^2-1$, $y \mapsto (1+y)^2-1$, $z \mapsto 0$, which is finite and becomes \'{e}tale after inverting $2$. The assertion follows from base change.
\end{proof}

\begin{rmk}
Note that the map $(r, \chi) \mapsto r \otimes \chi$ defines a morphism $\Spec (\Rzp \ctimes_{\O} R_\1) \rightarrow \Spec R_{\rbar}$ for all $p$, which is an isomorphism when $p>2$ (by Hensel's lemma) and has image in $\Spec R^{\sign}_{\rbar}$ if $p=2$ (since $\det (r \otimes \chi) (\Art_{\Q_2}(-1)) = \det r (\Art_{\Q_2}(-1)) = \chi(\Art_{\Q_2}(-1))$).
\end{rmk}
\section{The patching argument} \label{section:patchingargument}
In this section, we first introduce completed cohomology for quaternionic forms and then patch completed cohomology following \cite{MR3529394, 2016arXiv160906965G}. In the rest of the paper we assume $p=2$.

\subsection{Quaternionic forms and completed cohomology} \label{section:QMF}
Let $F$ be a totally real field and $D$ be a quaternion algebra with center $F$, which is ramified at all infinite places and at a set of finite places $\Sigma$, which does not contain any primes dividing $p$. We will write $\Sigma_p = \Sigma \cup \{v|p\}$. We fix a maximal order $\O_D$ of $D$, and for each finite places $v \notin \Sigma$ an isomorphism $(\O_D)_v \cong M_2(\O_{F_v})$. For each finite place $v$ of $F$, we will denote by $\Norm(v)$ the order of the residue field at $v$, and by $\varpi_v \in F_v$ a uniformizer.

Denote by $\Af \subset \A_F$ the finite adeles and adeles respectively. Let $U = \prod_{v} U_{v}$ be a compact open subgroup contained in $\prod_v (\O_D)_v^{\times}$. We may write
\begin{equation} \label{equation:decomp}
(D \otimes_F \Af)^{\times} = \bigsqcup_{i \in I} D^{\times} t_i U \Aft
\end{equation}
for some $t_i \in (D \otimes_F \Af)^{\times}$ and a finite index $I$. We say $U$ is sufficiently small if it satisfies the following condition:
\begin{equation} \label{equation:suffsmall}
(U(\A^f_F)^{\times} \cap t^{-1} D^{\times} t) / F^{\times} =1 \quad \text{for all }  t \in (D \otimes_F \Af)^{\times}.
\end{equation}
For example, $U$ is sufficiently small if for some  place $v$ of $F$, at which $D$ splits and not dividing $2M$ with $M$ being the integer defined in \cite[Lemma 3.1]{MR3544298}, $U_{v}$ is the pro-$v$ Iwahori subgroup (i.e. the subgroup whose reduction modulo $\varpi_{v}$ are the upper triangular unipotent matrices). We will assume this is the case from now on and denote the place by $v_1$.

Write $U = U^p U_p$, where $U_p = \prod_{v | p} U_v$ and $U^p = \prod_{v \nmid p} U_v$. If $A$ is a topological $\O$-algebra, we let $S(U^p, A)$ be the space of continuous functions
\[
f : D^{\times} \backslash (D \otimes_F \Af)^{\times} / U^p \rightarrow A.
\]
The group $G_p = (D \otimes_{\Z} \Zp)^{\times} \cong \prod_{v|p} \GL_2(F_v)$ acts continuously on $S(U^p, A)$. It follows from (\ref{equation:suffsmall}) that there is an isomorphism of $A$-modules
\begin{align} \label{equation:Sdecomp}
S(U^p, A) &\xrightarrow{\sim} \bigoplus_{i \in I} C( F^{\times} \backslash K_p \Aft, A) \\
f &\mapsto (u \mapsto f(t_i u))_{i \in I},
\end{align}
where $C$ denotes the space of continuous functions, $K_p = \prod_{v|p} \GL_2(\O_{F_v})$, and $I$ is the finite index set in the decomposition (\ref{equation:decomp}). Let $\psi: \Aft / F^{\times} \rightarrow \O^{\times}$ be a continuous character such that $\psi$ is trivial on $\Aft \cap U^p$. We may view $\psi$ as an $A$-valued character via $\O^{\times} \rightarrow A^{\times}$. Denote $S_{\psi}(U^p, A)$ be the $A$-submodule of $S(U, A)$ consisting of functions such that $f(gz) = \psi(z)f(g)$ for all $z \in \Aft$. The isomorphism (\ref{equation:Sdecomp}) induces an isomorphism of $U_p$-representations:
\begin{align} \label{equation:Szdecomp}
S_{\psi}(U^p, A) \xrightarrow{\sim} \bigoplus_{i \in I} C_{\psi}(K_p, A),
\end{align}
where $C_{\psi}$ denotes the continuous functions on which the center acts by the character $\psi$. One may think of $S_{\psi}(U^p, A)$ as the space of algebraic automorphic forms on $D^{\times}$ with tame level $U^p$ and no restrictions on the weight or level at places dividing $p$.

Let $\sigma$ be a continuous representation of $U_p$ on a free $\O$-module of finite rank, such that $\Aft \cap U_p$ acts on $\sigma$ by the restriction of $\psi$ to this group. We let
\[
S_{\psi, \sigma}(U, A) := \Hom_{U_p}(\sigma, S_{\psi}(U^p, A)).
\]
We will omit $\sigma$ as an index if it is the trivial representation. If the topology on $A$ is discrete (e.g. $A = E / \O$ or $A = \O / \varpi^s$), then we have
\[
S_{\psi}(U^p, A) \cong \varinjlim_{U_p} S_{\psi}(U^p U_p, A),
\]
where $U_p$ runs through compact open subgroups of $K_p$. The module $S_{\psi}(U^p, A)$ is naturally equipped with an $A$-linear action of $G_p := (D \otimes_{\Z} \Zp)^{\times} \cong \prod_{v | p} \GL_2(F_v)$, which extends the $K_p$-action. To be precise, for $g \in G_p$, right multiplication by $g$ induces an map 
\[
\cdot g : S_{\psi}(U^p U_p, A) \rightarrow  S_{\psi}(U^p U_p^g, A)
\]
for each $U_p$, where $U^g_p = g^{-1} U_p^g g$. As $U_p$ runs through the cofinal subset of open subgroups of $K_p$ with $U_p^g \subset K_p$, the subgroups $U_p^g$ also runs through a cofinal subset of open subgroups of $K_p$, so we may identify $\varinjlim_{U_p} S_{\psi}(U^p U_p^g, A)$  with $S_{\psi}(U^p, A)$. 

Denote $F_p = F \otimes_{\Q} \Qp \cong \prod_v F_v$ and $\O_{F_p} = \O_F \otimes_{\Z} \Zp \cong \prod_v \O_{F_v}$. Let $\zeta: F_p^{\times} \rightarrow \O^{\times}$ be the character obtained restricting $\psi$ to $F_p^{\times}$.

\begin{lemma} \label{lemma:Szinj}
The representation $S_{\psi}(U^p, E/\O)$ lies in $\Modladmz(\O)$. Moreover, $S_{\psi}(U^p, E/\O)$ is admissible and injective in $\Mod^{\sm}_{K_p, \zeta}(\O)$.
\end{lemma}

\begin{proof}
This follows from (\ref{equation:Szdecomp}).
\end{proof}

Let $S_p$ be the set of places of $F$ above $p$, $S_\infty$ be the set of places of $F$ above $\infty$, and let $S$ be a union of the places containing $\Sigma_p$, $S_{\infty}$, and all the places $v$ of $F$ such that $U_v \neq (\O_D)_v^{\times}$. Write $W = S - (\Sigma_p \cup S_{\infty})$. We will assume that for $v \in W$, $U_v \subset \GL_2(\O_{F_v})$ is contained in the Iwahori subgroup and contains the pro-$v$ Iwahori subgroup.

We denote $\T^S = \O[T_v, S_v, \U_{\varpi_w}]_{v \notin S, w \in W}$ be the commutative $\O$-polynomial algebra in the indicated formal variables. If $A$ is a topological $\O$-algebra then $S_{\psi}(U^p, A)$ and $S_{\psi, \sigma}(U^p, A)$ become $\T^S$-modules with $S_v$ acting via the double coset operator $[U_v \big( \begin{smallmatrix} \varpi_v  & 0 \\ 0 & \varpi_v \end{smallmatrix} \big) U_v ]$, $T_v$ acting via $[U_v \big( \begin{smallmatrix} \varpi_v  & 0 \\ 0 & 1 \end{smallmatrix} \big) U_v] $, and $\U_{\varpi_w}$ acting via $[U_w \big( \begin{smallmatrix} \varpi_w  & 0 \\ 0 & 1 \end{smallmatrix} \big) U_w]$. Note that the operators $T_v$ and $S_v$ do not depend on the choice of $\varpi_v$ but $\U_{\varpi_w}$ does.

\subsection{Completed homology and big Hecke algebras} \label{section:completedcohomology}
Let $S = S_p \cup S_{\infty} \cup \Sigma \cup \{v_1\}$, where $S_p$ be the set of places of $F$ above $p$ and $S_\infty$ be the set of places of $F$ above $\infty$. We define an open compact subgroup $U^p = \prod_{v \nmid p} U_v$ of $G(\A^{\infty, p}_F)$ as follows:
\begin{itemize}
\item $U_v = G(\O_{F_v})$ if $v \notin S$ or $v \in \Sigma$.
\item $U_{v_1}$ is the pro-$v_1$ Iwahori subgroup.
\end{itemize}
Due to the choice of $v_1$, $U^p U_p$ is sufficiently small for any open compact subgroup $U_p$ of $G(F_p)$. It follows that the functor $V \mapsto S_{\psi}(U^p U_p, V)$ is exact by (\ref{equation:Szdecomp}).

\begin{definition}
We define the completed homology groups $M_{\psi}(U^p)$ by
\[
M_{\psi}(U^p) := \varprojlim_{U_p} S_{\psi}(U^p U_p, \O)^{d}
\]
equipped with an $\O$-linear action of $G_p$ extending the $K_p$-action coming from the $\O \llbracket K_p \rrbracket$-module structure.
\end{definition}

Following from the definition, there is a natural $G_p$-equivariant homeomorphism
\[
M_{\psi}(U^p) \cong S_{\psi}(U^p, E/\O)^{\vee}.
\]

\begin{corollary} \label{cor:Mzproj}
The representation $M_{\psi}(U^p)$ is a projective object in $\Modpro_{K_p, \zeta}(\O)$.
\end{corollary}

\begin{proof}
Note that we have natural $G_p$-equivariant homeomorphism
\[
M_{\psi}(U^p) \cong S_{\psi}(U^p, E/\O)^{\vee}
\]
by definition. Thus the corollary follows from Lemma \ref{lemma:Szinj}.
\end{proof}

For $U = U^p U_p$, we write $S_{\psi}(U, s)$ for $S_{\psi}(U, \O / \varpi^s)$. Define $\T_{\psi}^S(U, s)$ to be the image of the abstract Hecke algebra $\T^S$ in $\End_{\O / \varpi^s [K_p / U_p]}(S_{\psi}(U, s))$.

\begin{definition}
We define the big Hecke algebra $\T_{\psi}^S(U^p)$ by
\[
\T_{\psi}^S(U^p) = \varprojlim_{U_p, s} \T_{\psi}^S(U^p U_p, s)
\]
where the limit is over compact open normal subgroups $U_p$ of $K_p$ and $s \in \Z_{\geq 1}$, and the surjective transition maps come from
\[
\End_{\O / \varpi^{s'} [K_p / U_p']}(S_{\psi}(U_p' U^p, s')) \rightarrow \End_{\O / \varpi^{s} [K_p / U_p]}(\O/\varpi^s[K_p/U_p] \otimes_{\O / \varpi^{s'}[K_p / U_p']} S_{\psi}(U_p' U^p, s'))
\]
for $s' \geq s$ and $U_p' \subset U_p$ and the natural identification
\[
\O/\varpi^s[K_p/U_p] \otimes_{\O / \varpi^{s'} [K_p / U_p']} S_{\psi}(U_p' U^p, s') \cong S_{\psi}(U_p U^p, s).
\]

We equip $\T_{\psi}^S(U^p)$ with the inverse limit topology. It follows from the definition that the action of $\T_{\psi}^S(U^p)$ on $M_{\psi}(U^p)$ is faithful and commutes with the action of $G_p$.
\end{definition}

\begin{lemma} \label{lemma:bigHeckemax}
$\T_{\psi}^S(U^p)$ is a profinite $\O$-algebra with finitely many maximal ideals. Denote its finitely many maximal ideals by $\m_1, \cdots, \m_r$ and let $J = \cap_{i} \m_i$ denote the Jacobson radical. Then $\T_{\psi}^S(U^p)$ is $J$-adically complete and separated, and we have
\[
\T_{\psi}^S(U^p) = \T_{\psi}^S(U^p)_{\m_1} \times \cdots \times \T_{\psi}^S(U^p)_{\m_r}.
\]
For each $i$, $\T_{\psi}^S(U^p) / \m_i$ is a finite extension of $k$.
\end{lemma}

\begin{proof}
This is indeed \cite[Lemma 2.1.14]{2016arXiv160906965G}. It suffices to prove when $U_p' \subset U_p$ are open normal pro-$p$ subgroups such that $\psi|_{U_p' \cap \O_{F, p}^{\times}}$ is trivial modulo $\varpi^{s'}$, the map
\[
\T^S_{\psi}(U^pU_p', s') \rightarrow \T^S_{\psi}(U^pU_p, 1)
\]
induces a bijection of maximal ideals.

Let $\m$ be a maximal ideal of the artinian ring $\T^S_{\psi}(U^pU_p', s')$. Since $\T^S_{\psi}(U^pU_p', s')$ acts faithfully on $S_{\psi}(U^pU_p', s')$, we know that
\[
S_{\psi}(U^pU_p', s')[\m] \neq 0.
\]
The $p$-group $U_p / U_p'$ acts naturally on this $k$-vector space, hence has a non-zero fixed vector, which belongs to $S_{\psi}(U^pU_p, 1)$. Thus $S_{\psi}(U^pU_p, 1)[\m] \neq 0$ and $\m$ is also a maximal ideal of $T^S_{\psi}(U^pU_p, 1)$.
\end{proof}

Let $\m \subset \T_{\psi}^S(U^p)$ be a maximal ideal with residue field $k$. There exists a continuous semi-simple representation $\rhobar_\m : G_{F, S} \rightarrow \GL_2(k)$ such that for any finite place $v \notin S$ of $F$, $\rhobar_\m(\Frob_v)$ has characteristic polynomial
$X^2 - T_v X + q_v S_v \in k[X]$. If $\rhobar_\m$ is absolutely reducible, we say that the maximal ideal $\m$ is Eisenstein; otherwise, we say that $\m$ is non-Eisenstein.

We define a global deformation problem
\[
\cS = (\rhobar_\m, F, S, \{\O\}_{v \in S}, \{\Dfz_v\}_{v \in S_p} \cup \{\Dodd_v\}_{v \in S_{\infty}} \cup \{\DSt_v \}_{v \in \Sigma} \cup \{\Dfz_{v_1}\}).
\]

\begin{prop} \label{prop:GalHecke}
Suppose that $\m$ is non-Eisenstein. Then there exists a lifting of $\rhobar_\m$ to a continuous homomorphism 
\[
\rho_\m: G_{F, S} \rightarrow \GL_2(\T_{\psi}^S(U^p)_\m)
\] 
such that for any finite place $v \notin S$ of $F$, $\rhobar_\m(\Frob_v)$ has characteristic polynomial $X^2 - T_v X + q_v S_v \in \T_{\psi}^S(U^p)_\m [X]$. Moreover, $\rho_\m$ is of type $\cS$ and has determinant $\psi \varepsilon$.
\end{prop}

\begin{proof}
By the proof of Lemma \ref{lemma:bigHeckemax}, the surjective map $\T_{\psi}^S(U^p) \twoheadrightarrow \T_{ \psi}^S(U^p U_p,s)$ induces bijection of maximal ideals for $U_p$ small enough. By taking projective limit, it suffices to show that there exist continuous homomorphism $\rhobar_{\m, U_p, s} : G_{F, S} \rightarrow \GL_2(\T_{\psi}^S(U^p U_p,s) / \m)$ and $\rho_{\m, U_p, s}: G_{F, S} \rightarrow \GL_2(\T_{\psi}^S(U^p U_p,s)_\m)$ satisfies the same conditions as in the statement, which follows from the well-known assertion for $S_{\psi}(U^p U_p, \O)$ (c.f. \cite[\S 1]{MR2290604}).
\end{proof}

\subsection{Globalization} \label{section: globalization}
Keeping the setting of Sect. \ref{section:completedcohomology}. Fix a continuous representation
\[
\rhobar: G_{F, S} \rightarrow \GL_2(k)
\]
which comes from a non-Eisenstein maximal ideal of $\T^S_{ \psi}(U^p)$ (i.e. $\rhobar \cong \rhobar_\m$). Assume $\rhobar$ satisfies the following properties:
\begin{enumerate}[label=(\roman*)]
\item $\rhobar$ has non-solvable image.
\item $\rhobar$ is unramified at all finite places $v \nmid p$;
\item $\rhobar(\Frob_{v_1})$ has distinct eigenvalues.
\end{enumerate}

In application to the modularity lifting theorem, assumption (ii) is satisfied after a solvable base change. The following lemma will allow us to reduce to situations where (iii) holds.

\begin{lemma} \label{lemma:disteigen}
Suppose $\rhobar$ has non-solvable image. Then there exists a place $v_1$ of $F$ not dividing $2Mp$ such that the eigenvalues of $\rhobar(\Frob_{v_1})$ are distinct.
\end{lemma}

\begin{proof}
By Dickson's theorem, the projective image of $\rhobar$ is conjugate to $\PGL_2(\F_{2^r})$ for some $r > 1$, which contains elements with distinct eigenvalues, e.g. $(\begin{smallmatrix}1 & 1 \\ 1 & 0 \end{smallmatrix})$. Thus by Chebotarev density theorem, there are infinite many places $v$ of $F$ with distinct Frobenius eigenvalues. This proves the lemma.
\end{proof}

\begin{definition}
Let $L$ be a finite extension of $\Qp$. Given a continuous representation $\rbar: G_L \rightarrow \GL_2(k)$, we will say that $\rbar$ has a suitable globalization if there is a totally real field $F$ and a continuous representation $\rhobar: G_F \rightarrow \GL_2(k)$ satisfying the properties $(i) - (iii)$ above and moreover,
\begin{itemize}
\item $\rhobar|_{G_{F_v}} \cong \rbar$ for each $v | p$ (hence $F_v \cong L$);
\item $[F:\Q]$ is even;
\item there exists a regular algebraic cuspidal automorphic representation $\pi$ of $\GL_2(\A_F)$ of weight $(0, 0)^{\Hom(F, \C)}$ and level prime to $p$ satisfying $\overline{\rho}_{\pi, \iota} \cong \rhobar$.
\end{itemize}
\end{definition}

Given a suitable globalization of $\rbar$, we set $S = S_p \cup S_{\infty} \cup \{v_1\}$, $\Sigma = \emptyset$, $D$ the quaternion algebra with center $F$ which is ramified exactly at $S_{\infty}$, and $U^p$ as in Sect. \ref{section:completedcohomology}. Let $\psi: G_{F, S} \rightarrow \O^{\times}$ be the totally even finite order character such that $\det \rho_{\pi, \iota} = \psi \varepsilon$ and view $\psi$ as a character of $\Aft / F^{\times} \rightarrow \O^{\times} $ via global class field theory. Let $\m$ be the maximal ideal of $\T^S_{\psi}(U^p)$ corresponding to $\rhobar$ and $\gamma$ be the character given by $\pi$. Together with the last property, we are in the same situation as Sect. \ref{section:completedcohomology}.

\begin{lemma} \label{lemma:suitglob}
Given $\rbar: G_{\Qp} \rightarrow \GL_2(k)$, there exists a suitable globalization.
\end{lemma}

\begin{proof}
By \cite[Proposition 3.2]{MR2869026}, we may find $F$ and $\rhobar$ satisfying all but the last two conditions. If $[F':\Q]$ is odd, we make a further quadratic extension $F''$ linearly disjoint from $\overline{F}^{\ker \rhobar}$ over $F$, and in which all primes above $p$ splits completely. The result follows by replacing $F$ with $F''$.

It is proved in \cite[Proposition 8.2.1]{2009arXiv0905.4266S} that when $p$ is odd, there is a finite Galois extension $F' / F$ in which all places above $p$ split completely such that $\rhobar |_{G_{F'}}$ is modular. This assumption can be removed using the proof of \cite[Theorem 6.1]{MR2551764}, which shows the existence of points for some Hilbert-Blumenthal abelian varieties with values in local fields when $p=2$.
\end{proof}

The following lemma says we may change the weight of a globalization $\rhobar$ when $p$ splits completely in $F$.

\begin{lemma}
Assume that $p$ splits completely in $F$ and that $\rhobar: G_F \rightarrow \GL_2(k)$ is automorphic. Then $\rhobar$ is automorphic of weight $\lambda = (0, 0)_{v|p}$, i.e. there is a regular algebraic cuspidal automorphic representation $\pi$ of weight $\lambda = (0, 0)_{v|p}$ such that $\rhobar \cong \rhobar_{\pi, \iota}$. Moreover,
\begin{enumerate}
\item at each $v|p$, $\rho_{\pi, \iota} |_{G_{F_v}}$ is semi-stable;
\item $\pi$ is $\iota$-ordinary at those $v | p$ for which $\rhobar |_{G_{F_v}}$ is reducible.
\end{enumerate}
\end{lemma}

\begin{proof}
It is proved in \cite[Lemma 3.29]{MR3544298} that if $\rhobar$ is automorphic, then it is automorphic of weight $(0 , 0)^{\Hom(F, \C)}$ and semi-stable at each $v|p$. The assertion (2) follows from \cite[Lemma 3.5]{MR2551764}, which proves that for a continuous representation $r: G_{\Qp} \rightarrow \GL_2(E)$,
\begin{itemize}
\item if $r$ is crystalline of weight $(0, 0)$, then it is ordinary if and only if residually it is ordinary;
\item if $r$ is semi-stable non-crystalline of weight $(0, 0)$, then it is ordinary.
\end{itemize}
This finishes the proof.
\end{proof}

\subsection{Auxiliary primes} \label{subdection:Auxprimes}
Let $Q$ be a set of places disjoint from $S$, such that for each $v \in Q$, $q_v \equiv 1$ mod $p$ and $\rhobar(\Frob_v)$ has distinct eigenvalues. For each $v \in Q$, we fix a choice of eigenvalue $\alpha_v$. We refer to the tuple $(Q, \{\alpha_v\}_{v \in Q})$ as a Taylor-Wiles datum. Denote $\Delta_Q = \prod_{v \in Q} \Delta_v = \prod_{v \in Q} k(v)^{\times}(p)$, and define the augmented deformation problem
\begin{align*}
\cS_Q = (&\rhobar, S \cup Q, \{\O \}_{v \in S} \cup \{\O[\Delta_v]\}_{v \in Q}, \{\Dfz_v\}_{v \in S_p} \cup \{\Dodd_v\}_{v \in S_{\infty}} \cup \{\DSt_v \}_{v \in \Sigma} \cup \{\Dfz_{v_1}\} \\
&\cup \{\DTW_v\}_{v \in Q}).
\end{align*}

Thus $R_{\cS_Q}$ is naturally a $\O[\Delta_Q]$-algebra. If $\a_Q \subset \O[\Delta_Q]$ is the augmentation ideal, then there is a canonical isomorphism $R_{\cS_Q} / \a_Q R_{\cS_Q} \cong R_\cS$ (resp. $R^T_{\cS_Q} / \a_Q R^T_{\cS_Q} \cong R^T_\cS$). 

\begin{lemma} \label{lemma:auxiliary}
Let $T = S$. For every $N \gg 0$, there exists a Taylor-Wiles datum  $(Q_N, \{\alpha_v\}_{v \in Q_N})$ satisfying the following conditions:
\begin{enumerate}
\item $\# Q_N := q = \dim_k H^1(G_{F, S}, \ad \rhobar) -2$.
\item For each $v \in Q_N$, $q_v \equiv 1$ (mod $p^N$).
\item The ring $R^{S, \psi}_{\cS_{Q_N}}$ is topologically generated by $2q+1$ elements over $A^S_{\cS}$.
\item Let $G_{Q_N}$ be the Galois group of the maximal abelian $2$-extension of $F$ over $F$ which is unramified outside $Q_N$ and is split at primes in $S$. Then we have $G_{Q_N} / 2^N G_{Q_N} \cong (\Z / 2^N \Z)^t$ with $t := 2 - |S| + q$.
\end{enumerate}
\end{lemma}

\begin{proof}
See \cite[Lemma 5.10]{MR2551764}.
\end{proof}

\subsubsection{Action of $\Theta_Q$}
If $Q$ is a finite set of finite primes of $F$ disjoint from $S$, we denote by $\Theta_Q$ the Galois group of the maximal abelian $2$-extension of $F$ which is unramified outside $Q$ and in which every prime in $S$ splits completely. Let $\Theta^*_{Q}$ be the formal group scheme defined over $\O$ whose $A$-valued points is given by the group $\Hom(\Theta_{Q}, A)$ of continuous characters on $\Theta_{Q}$ that reduce to the trivial character modulo $\m_A$.

It follows that $\Spf R_{\cS_{Q}}$ (resp. $\Spf R^T_{\cS_{Q}}$) has a natural action by $\Theta^*_{Q}$ given by $\chi_A \times V_A \mapsto V_A \otimes \chi_A$ on $A$-valued points, which is free if $\rhobar$ has non-solvable image \cite[Lemma 5.1]{MR2551764}. Moreover, there is a $\Theta^*_{Q}$-equivariant map
\begin{align} \label{equation:deltaQ}
    \delta_{Q} : \Spf R^T_{\cS_{Q}} \rightarrow \Theta^*_{Q}; \qquad V_A \mapsto \det V_A \cdot (\psi \varepsilon)^{-1}
\end{align}
where $\Theta^*_{Q}$ acts on itself via the square of the identity map, and $ \Spf R^{T, \psi}_{\cS_{Q}} = \delta_{Q}^{-1}(1)$.

\subsection{Auxiliary levels}
A choice of Taylor-Wiles datum $(Q, \{\alpha_v\}_{v \in Q})$ having been fixed, we have defined an auxiliary deformation problem $\cS_{Q}$.

Let $U^p$ be the open compact subgroup of $G(\A^{\infty, p}_F)$ in Sect. \ref{section:completedcohomology}. We define compact open subgroups $U^p_0(Q) = \prod_{v \nmid p} U_0(Q)_v$ and $U^p_1(Q) = \prod_{v \nmid p} U_1(Q)_v$ of $U^p = \prod_{v \nmid p}U_v$ by:
\begin{itemize}
\item if $v \notin Q$, then $U_0(Q)_v = U_1(Q)_v = U_v$.
\item if $v \in Q$, then $U_0(Q)_v$ is the Iwahori subgroup of $\GL_2(\O_{F_v})$ and $U_1(Q)_v$ is the set of $g = (\begin{smallmatrix} a & b \\ c & d \end{smallmatrix}) \in U_0(Q)_v$ such that $ad^{-1}$ maps to $1$ in $\Delta_v$.
\end{itemize}
In particular, $U_1(Q)_v$ contains the pro-$v$ Iwahori subgroup of $U_0(Q)_v$, so we may identify $\prod_{v \in Q}U_0(Q)_v / U_1(Q)_v$ with $\Delta_{Q}$.

Let $\m_Q$ denote the ideal of $\T^{S \cup Q}$ generated by $\m \cap \T^{S \cup Q}$ and the elements $\U_{\varpi_v} - \tilde{\alpha}_v$ for $v \in Q$, where $\tilde{\alpha}_v$ is an arbitrary lift of $\alpha_v$. We denote by $\T^{S \cup Q}_{\psi}(U^p_i(Q)U_p, s)$ the image of $\T^{S \cup Q}$ in $\End_{\O / \varpi^s}(S_{\psi}(U^p_i(Q)U_p, s))$. Exactly as \cite[\S 2.1]{MR2505297}, we have the following:
\begin{enumerate}
\item The maximal ideal $\m_Q$ induces proper, maximal ideals in $\T^{S \cup Q}_{\psi}(U^p_i(Q)U_p, s)$. Moreover, the map
\[
S_{\psi}(U^pU_p, s)_\m \rightarrow S_{\psi}(U^p_0(Q)U_p, s)_{\m_Q}
\]
is an isomorphism.
\item $S_{\psi}(U^p_1(Q)U_p, s)_{\m_Q}$ is a finite projective $\O/ \varpi^s[\Delta_Q]$-module with
\[
S_{\psi}(U^p_1(Q)U_p, s)_{\m_Q}^{\Delta_Q} \xrightarrow{\sim} S_{\psi}((U^p_0(Q)U_p, s)_{\m_Q}.
\]
\item There is a deformation
\[
\rho_{\m, Q, s}: G_F \rightarrow \GL_2(\T^{S \cup Q}_{\psi}(U^p_1(Q)U_p, s))
\]
of $\rhobar$ which is of type $\cS_Q$ and has determinant $\psi \varepsilon$. In particular, $S_{\psi}(U^p_1(Q)U_p, s)_{\m_Q}$ is a finite $R^{\psi}_{\cS_Q}$-module.
\end{enumerate}

The following proposition is an immediate consequence of $(3)$.

\begin{prop}
Let $(Q, \{ \alpha_v\}_{v \in Q})$ be a Taylor-Wiles datum. Then there exists a lifting of $\rhobar_\m$ to a continuous morphism
\[
\rho_{\m, Q}: G_{F, S \cup Q} \rightarrow \GL_2(\T^{S \cup Q}_{\psi}(U^p_1(Q))_{\m_{Q, 1}})
\]
satisfying the following conditions:
\begin{itemize}
\item for each place $v \notin S \cup Q$ of $F$, $\rho_{\m, Q}(\Frob_v)$ has characteristic polynomial $X^2 - T_v X + q_v S_v \in \T^{S \cup Q}_{\psi}(U^p_1(Q))_{\m_{Q, 1}} [X]$;
\item for each place $v \in Q$, $\rho_{\m, Q} |_{G_{F_{v}}} \sim \big(\begin{smallmatrix}
\chi_{v} & * \\ 0 & *
\end{smallmatrix} \big)$ such that $\chi_{v} \circ \Art_{F_v}(\varpi^{-1}_v) = \U_{\varpi_{v}}$.
\end{itemize}
In particular, $\rho_{\m, Q}$ is of type $\cS_Q$ and has determinant $\psi \varepsilon$.
\end{prop}

It follows that we have an $\O[\Delta_Q]$-algebra surjection
\begin{align} \label{equation:GaloisHecke}
\Rz_{\cS_Q} \twoheadrightarrow \T^{S \cup Q}_{\psi}(U_1^p(Q))_{\m_Q}
\end{align}
such that for $v \notin S$ the trace of $\Frob_v$ on the universal deformation of type $\cS_Q$ maps to $T_v$ and $\chi_v(\varpi_v)$ maps to $\U_{\varpi_v}$ for $v \in Q$.

\subsubsection{Action of $\Theta_{Q}$}
Let $\chi \in \Theta^*_{Q}(\O)[2]$ be a character of $G_Q$ of order 2. As $\chi$ is split at infinite places, we can regard $\chi$ also as a character $\Aft$. Given $f \in S_{\psi}(U^p_1(Q) U_p, \O)$, we define
\[
f_{\chi}(g) := f(g) \chi(\det(g)),
\]
which also lies in $S_{\psi}(U^p_1(Q) U_p, \O)$. This induces an action of $\Theta^*_{Q}(\O)[2]$ on $S_{\psi}(U_p U^p_1(Q), s)$ for each $s \in \N$. By Proposition 7.6 of \cite{MR2551764}, we may also define an action $\chi$ on $\T^{S \cup Q}_{\psi}(U_1^p(Q))$ and $\O[\Delta_N]$ by sending $T_v$ to $\chi(\varpi_v) T_v$, $S_v$ to $\chi(\varpi_v) S_v$ and $\langle h \rangle$ to $\chi(h) \langle h \rangle$, which is compatible with the action of $\chi$ on $S_{\psi}(U_p U^p_1(Q), s)$. Moreover, the action of $\chi$ on $\T^{S \cup Q}_{\psi}(U_1^p(Q))$ preserves its maximal ideal $\m_Q$ and the homomorphism $\Rz_{\cS_{Q}} \rightarrow \T^{S \cup Q}_{\psi}(U_1^p(Q))_{\m_Q}$ is $\Theta^*_{Q}(\O)[2]$-equivariant.

\subsection{Patching} \label{section:patching}
We write $G_p$ for $\prod_{v|p} \GL_2(F_v)$, $K_p$ for $\prod_{v|p} \GL_2(\O_{F_v})$ and $Z_p \cong \prod_{v|p} F_v^{\times}$ for the center of $G_p$.

We let $(Q_N, \{\alpha_v\}_{v \in Q_N})$ be a choice of Taylor-Wiles datum for each $N \gg 0$ and $T = S$ be the subset as in Lemma \ref{lemma:auxiliary}. Choose $v_0 \in S$, and let $\cT = \O \llbracket X_{v, i, j} \rrbracket_{v \in S, 1 \leq i, j \leq 2} / (X_{v_0, 1, 1})$. By Lemma \ref{lemma:Galoisframe}, there is a canonical isomorphism $R^S_{\cS} \cong R_\cS \ctimes_\O \cT$ (resp. $R^{S, \psi}_{\cS} \cong \Rz_\cS \ctimes_\O \cT$). Let $\Deltainf = \Zp^{q}$, which is endowed with a natural surjection $\Deltainf \twoheadrightarrow \Delta_{Q_N}$ given by $(\Zp)^q \twoheadrightarrow (\Z / p^N \Z)^q \cong \prod_{v \in Q_N} k(v)^\times(p)$ for each $N$. This induces a surjection $\Oinf := \cT \llbracket \Deltainf \rrbracket \rightarrow \O_N := \cT \llbracket \Delta_N \rrbracket$ of $\cT$-algebras. Denote the kernel of the homomorphism $\Oinf \rightarrow \O$ which sends $\Deltainf$ to $1$ and all $4|S| - 1$ variables of $\cT$ to $0$ by $\a$.

We write $\Rloc$ for $A^S_{\cS}$ and denote $g = q + |S| - 1$. Fix a surjection $\Z_2^t \rightarrow \Theta_{Q_N}$ for each $N$. This induces an embedding of formal group scheme $\iota: \Theta^*_{Q_N} \hookrightarrow (\Gmc)^t$, where $\Gmc$ denotes the completion of the $\O$-group scheme $\Gm$ along the identity section. We define
\begin{itemize}
\item $\Rinf' = \Rloc \llbracket X_1, \ldots, X_{g+t} \rrbracket$. Then $\Spf \Rinf'$ is equipped with a free action of $(\Gmc)^t$, and a $(\Gmc)^t$-equivariant morphism $\delta: \Spf \Rinf' \rightarrow (\Gmc)^t$ induced by $\delta_{Q_N}$ (\ref{equation:deltaQ}), where $(\Gmc)^t$ acts on itself by the square of the identity map.
\item $\Rinf$ by $\Spf \Rinf = \delta^{-1}(1)$ and $\Rinvinf$ by $\Spf \Rinvinf := \Spf \Rinf' / (\Gmc)^t$ (cf. \cite[Proposition 2.5]{MR2551764}). By \cite[Lemma 9.4]{MR2551764}, $\Spf \Rinf'$ is a $(\Gmc)^t$-torsor over $\Spf \Rinvinf$.
\end{itemize}
We fix a $\Theta^*_{Q_N}$-equivariant surjective $\Rloc$-algebra homomoprhism $\Rinf' \twoheadrightarrow R^S_{\cS_{Q_N}}$ for each $N$, which induces a $\Theta^*_{Q_N}[2]$-equivariant surjective $\Rloc$-algebra homomorphism $\Rinf \twoheadrightarrow R^{S, \psi}_{\cS_{Q_N}}$.

\begin{definition}
Let $U_p$ be a compact open subgroup of $K_p$ and let $J$ be an open ideal in $\Oinf$. Let $I_J$ be the subset of $N \in \N$ such that $J$ contains the kernel of $\Oinf \rightarrow \O_N$. For $N \in I_J$, define
\[
M(U_p, J, N) := \Oinf / J \otimes_{\O_N} S_{\psi}(U_1^p(Q_N) U_p, \O)_{\m_{Q_N}}^d.
\]
\end{definition}

From the definition, it follows that $M(U_p, J, N)$ satisfies the following properties:
\begin{itemize}
\item We have a map
\begin{equation} \label{equation:RNtoTN}
R^{S, \psi}_{\cS_{Q_N}} \rightarrow \cT \ctimes_{\O} \T^{S}_{\psi}(U^p_1(Q_N))_{\m_{Q_N}},
\end{equation}
and a map 
\begin{equation} \label{equation:TNtoEndMN}
\cT \ctimes_{\O} \T^{S}_{\psi}(U^p_1(Q_N))_{\m_{Q_N}} \rightarrow \End_{\Oinf/J}(M(U_p, J, N)).
\end{equation}
In particular, for all $J$ and $N \in I_J$ we have a ring homomorphism
\[
\Rinf \rightarrow \End_{\Oinf/J} \big(M(U_p, J, N) \big)
\]
which factors through our chosen quotient map $\Rinf \rightarrow R^{S, \psi}_{\cS_{Q_N}}$ and the maps (\ref{equation:RNtoTN}), (\ref{equation:TNtoEndMN}). Moreover, it is $\Theta^*_{Q_N}[2]$-equivariant.
\item If $U'_p$ is an open normal subgroup of $U_p$, then $M(U'_p, J, N)$ is projective in the category of $\Oinf/J[U_p / U'_p]$-module with central character $\psi^{-1}|_{\O_{F_p}^{\times}}$.
\item Suppose that $\a \subset J$. Then $M(U_p, J, N) = S_{\psi}(U^p U_p, s(J))_{\m}^{\vee}$, where $\Oinf / J \cong \O / \varpi^{s(J)}$.
\end{itemize}

\begin{definition}
For $d \geq 1$, $J$ an open ideal in $\Oinf$ and $N \in I_J$, we define
\begin{align*}
R(d, J, N) &:= \Oinf/J \otimes_{\O_N} (R^{S, \psi}_{\cS_{Q_N}} / \m^d_{R^{T, \psi}_{\cS_{Q_N}}}).
\end{align*}
\end{definition}

We have the following properties:
\begin{itemize}
\item Each ring $R(d,J,N)$ is a finite commutative local $\Oinf / J$-algebra, equipped with a surjective $\O$-algebra homomorphism 
\[
\Rinf \twoheadrightarrow R(d, J, N).
\]
\item For $d$ sufficiently large, the map $\Rinf \rightarrow \End_{\Oinf/J} (M(U_p, J, N) )$ factors through $R(d, J, N).$
\item We have an isomorphism 
\begin{align*}
R(d, J, N) / \a R(d, J, N) \cong \Rz_{\cS} / (\m^d_{\Rz_{\cS}}, \varpi^{s(\a+J)}).
\end{align*}
\item For all open ideals $J' \subset J$ and open normal subgroups $U_p' \subset U_p$, we have a surjective map
\[
M(U_p', J', N) \rightarrow M(U_p, J, N)
\]
inducing an isomorphism
\[
\Oinf/J \otimes_{\Oinf/J' [U_p / U_p']} M(U_p', J', N) \rightarrow M(U_p, J, N).
\]
\item If $U_p$ is an open normal subgroup of $K_p$, then $\{M(U_p, J, N)\}_{N \in I_J}$ is a set of projective objects in the category of $\Oinf/J [K_p / U_p]$-modules with central character $\psi^{-1}|_{\O_{F_p}^{\times}}$.
\end{itemize}

We fix a non-principal ultrafilter $\filter$ on the set $\N$.

\begin{definition}
Let $(\Oinf/J)_{I_J} = \prod_{i \in I_J} \Oinf/J$ and $x \in \Spec \big((\Oinf/J)_{I_J} \big)$ given by $\filter$. We define
\begin{align*}
M(U_p, J, \infty) := (\Oinf/J)_{I_J, x} \otimes_{(\Oinf/J)_{I_J}} \bigg( \prod_{N \in I_J} M(U_p, J, N) \bigg), \\
R(d, J, \infty) := (\Oinf/J)_{I_J, x} \otimes_{(\Oinf/J)_{I_J}} \bigg( \prod_{N \in I_J} R(d, J, N) \bigg).
\end{align*}
\end{definition}

We have the following
\begin{itemize}
\item If $U_p$ is an open normal subgroup of $K_p$, then $M(U_p, J, \infty)$ is projective in the category of $\Oinf/J[K_p / U_p]$-module with central character $\psi^{-1}|_{\O_{F_p}^{\times}}$.

\item If $\a \subset J$, there is a natural isomorphism 
\begin{equation} \label{equation:Minfmodaisom}
M(U_p, J, \infty) / \a M(U_p, J, \infty) \cong S_{\psi}(U^p U_p, s(J))_{\m}^{\vee}.
\end{equation} 
\item For $d$ sufficiently large, the map 
\begin{equation} \label{equation:RinftoMinf}
\Rinf \rightarrow \End_{\Oinf/J}(M(U_p, J, \infty))
\end{equation} 
factors through $R(d, J, \infty)$ and the map
\begin{equation} \label{equation:RdinftoMinf}
R(d, J, \infty) \rightarrow \End_{\Oinf/J}(M(U_p, J, \infty))
\end{equation}
is an $\Oinf$-algebra homomorphism. Moreover, both (\ref{equation:RinftoMinf}) and (\ref{equation:RdinftoMinf}) are $\Theta^*_{Q_N}[2]$-equivariant.
\item We have an isomorphism
\begin{equation} \label{equation:Rinfmodaisom}
R(d, J, \infty) / \a \cong R_{\cS} / (\m^d_{R_{\cS}}, \varpi^{s(\a+J)}).
\end{equation}
\item For all open ideals $J' \subset J$ and open normal subgroups $U_p' \subset U_p$, the natural map
\[
M(U_p', J', \infty) \rightarrow M(U_p, J, \infty)
\]
is surjective, and induces an isomorphism of $\Oinf/J$-modules
\begin{equation} \label{equation:isom}
\Oinf/J \otimes_{\Oinf/J' [U_p / U_p']} M(U_p', J', \infty) \rightarrow M(U_p, J, \infty).
\end{equation}
\end{itemize}

\begin{definition}
We define an $\Oinf \llbracket K_p \rrbracket$-module
\[
\Minf := \varprojlim_{J, U_p} M(U_p, J, \infty).
\]
\end{definition}

We claim the following hold.

\begin{itemize}
\item $\Minf$ is endowed with an action of $\Rinf$ via the map $\alpha: \Rinf \rightarrow \varprojlim_{J, d} R(d, J, \infty)$. Since the image of $\alpha$ contains the image of $\Oinf$, $\alpha(\Rinf)$ is naturally an $\Oinf$-algebra. Since $\Oinf$ is formally smooth, we can choose a lift of the map $\Oinf \rightarrow \alpha(\Rinf)$ to a map $\Oinf \rightarrow \Rinf$. We make such a choice, and regard $\Rinf$ as an $\Oinf$-algebra and $\alpha$ as an $\Oinf$-algebra homomorphism.
\item The module $\Minf$ is naturally equipped with an $\Oinf$-linear action of $G_p$, which extends the $K_p$-action coming from the $\Oinf \llbracket K_p \rrbracket$-structure. To be precise, for $g \in G_p$, right multiplication by $g$ induces an map 
\[
\cdot g : M(U_p, J, N) \rightarrow  M(g^{-1} U_p g, J, N)
\]
for each $U_p, J, N$. Suppose that $g^{-1} U_p g \subset K_p$, our construction gives a map 
\[
\cdot g : M(U_p, J, \infty) \rightarrow  M(g^{-1} U_p g, J, \infty).
\]
As $U_p$ runs through the cofinal subset of open subgroups of $K_p$ with $g^{-1} U_p g \subset K_p$, the subgroups $g^{-1} U_p g$ also runs through a cofinal subset of open subgroups of $K_p$, so we may identify $\varprojlim_{J, U_p} M(g^{-1} U_p g, J, \infty)$  with $\Minf$. Taking the inverse limit over $J$ and $U_{\infty}$ gives the action of $g$ on $\Minf$.
\end{itemize}

\begin{prop} \label{prop:Minfprop}
\begin{enumerate}
\item For all open ideals $J$ and open compact subgroups $U_p$ of $K$, we have a sujective map
\[
\Minf \rightarrow M(U_p, J, \infty)
\]
inducing isomorphism
\[
\Oinf/J \otimes_{\Oinf/J[U_p]} \Minf \rightarrow M(U_p, J, \infty).
\]
\item There is a $\Theta^*_{Q_N}[2]$-equivariant homomorphism $\Rinf \rightarrow \End_{\Oinf \llbracket K \rrbracket}(\Minf)$ which factors as the composite of $\Oinf$-homomorphisms $\Rinf \rightarrow \varprojlim_{J, d} R(d, J, \infty)$ and $\varprojlim_{J, d} R(d, J, \infty) \rightarrow \End_{\Oinf \llbracket K_p \rrbracket}(\Minf)$ given by the homomorphisms above.
\item $\Minf$ is finitely generated over $\Oinf \llbracket K_p \rrbracket$ and projective in the category $\Modpro_{K_p, \zeta}(\Oinf)$, with $\zeta = \psi|_{\O_{F_p}^{\times}}$. In particular, it is finitely generated over $\Rinf \llbracket K_p \rrbracket$ and projective in $\Modpro_{K_p, \zeta}(\O)$.
\end{enumerate}
\end{prop}

\begin{proof}
The first assertion follows from the isomorphism (\ref{equation:isom}) and the second assertion can be deduced easily by the definition of $\Minf$. To show the third assertion, note that it is proved in \cite[Proposition 2.10]{MR3529394} (see \cite[Proposition 3.4.16 (1)]{2016arXiv160906965G} also) that $\Minf$ is finitely generated over $\Oinf \llbracket K_p \rrbracket$ and projective in the category $\Modpro_{K_p, \zeta}(\Oinf)$. We claim that the following conditions are equivalent for a compact module $M$ over a complete local $\varpi$-torsion free $\O$-algebra $R$:
\begin{align*}
&\text{$M$ is projective in $\Modpro_{K_p, \zeta}(R)$} \\
\Longleftrightarrow & \text{$M$ is $\varpi$-torsion free and $M / \varpi M$ is projective in $\Modpro_{K_p, \zeta}(R / \varpi)$} \\
\Longleftrightarrow & \text{$M$ is $\varpi$-torsion free and $M / \varpi M$ is projective in $\Modpro_{I_p, \zeta}(R / \varpi)$} \\
\Longleftrightarrow & \text{$M$ is $\varpi$-torsion free, and $M / \varpi M \cong \prod_{i \in J} R/\varpi \llbracket I_p / I_p \cap Z_p \rrbracket$}
\end{align*}
where $I_p$ is the pro-$p$ Iwahori subgroup of $G_p$ and $J$ is an index set. Given the claim, we see that $\Minf / \varpi \Minf \cong \prod_{J} \Oinf /\varpi \llbracket I_p / I_p \cap Z_p \rrbracket$. Since $\Oinf / \varpi \cong k \llbracket x_1, \dots, x_q \rrbracket \cong \prod_{J'} k$ for some index set $J'$ as $k$-vector spaces, we have $\Minf / \varpi \Minf \cong \prod_{J} \prod_{J'} k \llbracket I_p / I_p \cap Z_p \rrbracket$ as compact $I_p$-modules and thus $\Minf$ is projective in $\Modpro_{K_p, \zeta}(\O)$ by the claim.

To show the first equivalence, we first assume that $M$ is projective in $\Modpro_{K_p, \zeta}(R)$. Note that the map $K'_p \rightarrow (K'_p / K'_p \cap Z_p) \times \Gamma_p$, $g \mapsto (g (K'_p \cap Z_p), (\det g)^{-1})$, where $K'_p = \{g = sz \mid s =(s_v) \in \prod_{v | p} \SL_2(\O_{F_v}), \ s_v \equiv (\begin{smallmatrix} 1 & 0 \\ 0 & 1 \end{smallmatrix}) \text{ mod } \varpi_v^2, \ z \in \prod_{v | p} (1 + \varpi_v^2 \O_{F_v}) \}$ and $\Gamma_p = (K'_p \cap Z_p)^2$, is an isomorphism of groups. It follows that $R \llbracket K'_p \rrbracket \cong R \llbracket (K'_p / K'_p \cap Z_p) \rrbracket \ctimes_R R \llbracket \Gamma_p \rrbracket$. Viewing $M$ as compact $R \llbracket K'_p \rrbracket$-module, we see that it is a quotient of $\prod_{j} R \llbracket K'_p \rrbracket$ and thus a quotient of $\prod_{j} R \llbracket K'_p \rrbracket / (z - \zeta^{-1}(z))_{z \in K'_p \cap Z_p} \cong \prod_{j} R \llbracket (K'_p / K'_p \cap Z_p) \rrbracket$. Since $M$ is projective in $\Modpro_{K_p, \zeta}(R)$, it is projective in $\Modpro_{K'_p, \zeta}(R)$ and hence a direct summand of $\prod_{j} R \llbracket (K'_p / K'_p \cap Z_p) \rrbracket$. This shows that $M$ is $\varpi$-torsion free. Note that for every $N$ in $\Modpro_{K_p, \zeta}(R / \varpi)$ we have $\Hom(M, N) \cong \Hom(M/\varpi M, N)$ thus $M/\varpi M$ is projective in $\Modpro_{K_p, \zeta}(R / \varpi)$. On the other hand, suppose that $M$ is $\varpi$-torsion free and $M/\varpi M$ is projective in $\Modpro_{K_p, \zeta}(R / \varpi)$. Let $P$ be the projective envelope of $M / \varpi M$ in $\Modpro_{K_p, \zeta}(R)$. It follows that there is a morphism $P \rightarrow M$ lifting $P \twoheadrightarrow M/\varpi M$. This morphism is surjective by the Nakayama's lemma for compact modules ($P / \varpi P \cong M / \varpi M$). Denote $K$ to be the kernel of this morphism, we have $K / \varpi K = 0$ because $P / \varpi P \cong M / \varpi M$ and $0 \rightarrow K / \varpi K \rightarrow P / \varpi P \rightarrow M / \varpi M$ is exact (5-lemma). This implies $K = 0$ (by the Nakayama's lemma for compact modules) and thus $M \cong P$. The second equivalence is because $I_p$ is the pro-$p$ Sylow subgroup of $K_p$. Since $\zeta$ mod $\varpi$ is trivial on $I_p / I_p \cap Z_p$, $M / \varpi M$ is a compact module over $R / \varpi \llbracket I_p / I_p \cap Z_p \rrbracket$ and the third equivalence follows from the fact that a compact $R / \varpi \llbracket I_p / I_p \cap Z_p \rrbracket$-module is projective if and only if it is pro-free (because $R / \varpi \llbracket I_p / I_p \cap Z_p \rrbracket$ is local, projectivity coincides with freeness). This proves the proposition.
\end{proof}

\begin{prop} \label{prop:completecohomologymoda}
Let $\a = \ker(\Oinf \rightarrow \O)$ as before, we have a natural ($G$-equivariant) isomorphism
\[
\Minf / \a \Minf \cong M_{\psi}(U^p)_\m.
\]
There is a surjective map $\Rinf / \a \Rinf \rightarrow \Rz_{\cS} \rightarrow \T^S_{\psi}(U^p)_{\m}$ and the above isomorphism intertwines the action of $\Rinf$ on the left hand side with the action of $\T^S_{\psi}(U^p)_{\m}$ on the right hand side.
\end{prop}

\begin{proof}
Note that we have a isomorphism (\ref{equation:Minfmodaisom}). To prove the first part, it suffices to show that we have an isomorphism
\[
\Minf / \alpha \Minf \cong \varprojlim_{J, U_p} M(U_p, J, \infty) / \alpha M(U_p, J, \infty),
\]
which follows from \cite[Lemma A.33]{2016arXiv160906965G} (see also \cite[Corollary 2.11]{MR3529394}). The second part is an immediate consequence of isomorphism (\ref{equation:Rinfmodaisom}).
\end{proof}
\section{Patching and Breuil-M\'{e}zard conjecture} \label{section:patchBM}
We assume that $p$ ($=2$) splits completely in $F$. Equivalently, $F_v \cong \Qp$ for all $v|2$. Let $\rbar :G_{\Qp} \rightarrow \GL_2(k)$ be a continuous representation. We note that all the results in this section can be extended to arbitrary prime $p$ and general totally real field $F$ (by a similar method as in \cite{MR3134019}), we restrict ourself to this particular case since it is sufficient for our purpose.

\subsection{Local results}
\subsubsection{Locally algebraic type}
Fix a Hodge type $\lambda$, and inertia type $\tau$, and a continuous character $\zeta: G_{\Qp} \rightarrow \O^{\times}$ such that $\zeta |_{I_{\Qp}} = (\Art_{\Qp}^{-1})^{\lambda_1 + \lambda_2} \cdot \det \tau$. We define $\sigma(\lambda, \tau)= \sigma(\lambda) \otimes_E \sigma(\tau)$, where $\sigma(\lambda) = \sigma(\lambda) = M_{\lambda} \otimes_{\O} E$ and $\sigma(\tau)$ be the smooth type corresponding to $\tau$ (see Notations for the precise definition). Since $\sigma(\lambda, \tau)$ is a finite dimensional $E$-vector space and $K$ is compact and the action of $K$ on $\sigma(\lambda, \tau)$ is continuous, there is a $K$-stable $\O$-lattice $\sigma^{\circ}(\lambda, \tau)$ in $\sigma(\lambda, \tau)$. Then $\sigma^{\circ}(\lambda, \tau) / (\varpi)$ is a smooth finite length $k$-representation of $K$, we will denote by $\overline{\sigma(\lambda, \tau)}$ its semi-simplification. One may show that $\overline{\sigma(\lambda, \tau)}$ does not depends on the choice of a lattice. The same assertion holds for $\sigma^{cr}(\lambda, \tau) = \sigma(\lambda) \otimes \sigma^{cr}(\tau)$.

A locally algebraic type $\sigma$ is an absolutely irreducible representation of $\GL_2(\Qp)$ of the form $\sigma(\lambda, \tau)$ or $\sigma^{cr}(\lambda, \tau)$ for some inertial type $\tau$ and Hodge type $\lambda$. We say that a continuous representation $r: G_{\Qp} \rightarrow \GL_2(E)$ has type $\sigma = \sigma(\lambda, \tau)$ (resp. $\sigma^{cr}(\lambda, \tau)$) if it is potentially semi-stable (resp. potentially crystalline) of inertial type $\tau$ and Hodge type $\lambda$. Denote $\Rzp(\sigma)$ the local universal lifting ring of type $\sigma$ and determinant $\zeta \varepsilon$ for $\rbar$. 

If $x$ is a point of $\Spec \Rzp(\sigma)[1/p]$ with residue field $E_x$, we denote by $r_x : G_{\Qp} \rightarrow \GL_2(E_x)$ the lifting of $\rbar$ given by $x$. We define the locally algebraic $G$-representation $\pi_{\lalg}(r_x) = \pism(r_x) \otimes_{E_x} \pialg(r_x)$. Note that $\H(\sigma) := \End_G(\cInd^G_K(\sigma))$ acts via a character on the one-dimensional space $\Hom_{\GL_2(\Zp)}(\sigma, \pi_{\lalg}(r_x))$ (see the appendix to \cite{MR1944572}).

\begin{thm} \label{thm:HeckeGalois}
There is an $E$-algebra homomorphism
\[
\phi : \H(\sigma) \rightarrow \Rzp(\sigma)[1/p]
\]
which interpolates the local Langlands correspondence. More precisely, for any closed point $x$ of $\Spec \Rzp(\sigma)[1/p]$, the $\H(\sigma)$-action on $\Hom_{\GL_2(\Zp)}(\sigma, \pi_{\lalg, x})$ factors as $\phi$ composed with the evaluation map $\Rzp(\sigma)[1/p] \rightarrow E_x$.
\end{thm}

\begin{proof}
This follows from \cite[Theorem 4.1]{MR3529394} for $\sigma = \sigma^{cr}(\lambda, \tau)$ and \cite[Theorem 3.3]{2018arXiv180301610P} for $\sigma = \sigma(\lambda, \tau)$.
\end{proof}

\subsubsection{The Breuil-M\'{e}zard conjecture} \label{section:BM}
We now state the Breuil-M\'{e}zard conjecture \cite{MR1944572}.

\begin{conj}[Breuil-M\'ezard] \label{conj:Breuil-Mezardm}
There exist non-negative integers $\mu_{a}$ for each Serre weight $a$ of $\GL_2(k)$ such that for each locally algebraic type $\sigma$, we have
\begin{align*}
e(\Rzp(\sigma) / \varpi) = \sum_{a} m_{a}(\sigma) \mu_{a}(\rbar)
\end{align*}
where $a$ runs over all Serre weights (see Sect. \ref{RepGL2}), and $m_{a}(\sigma)$ is the multiplicity of $\bsigma_a$ as a Jordan-Holder factor of $\overline{\sigma}$.
\end{conj}

There is also a geometric version of the Breuil-M\'ezard conjecture due to \cite{MR3134019}.

\begin{conj} \label{conj:Breuil-Mezard}
For each Serre weight $a$ of $\GL_2(k)$, there exists a 4-dimensional cycle $\cC_{a}(\rbar)$ of $\Rzp$, independent of $\lambda$ and $\tau$, such that for each $\lambda, \tau$, we have equalities of cycles:
\begin{align*}
Z(\Rzp(\sigma) / \varpi) = \sum_{a} m_a(\sigma) \cC_{a}(\rbar)
\end{align*}
where $a$ runs over all Serre weights and $m_{a}(\sigma)$ is as in the previous conjecture.
\end{conj}

\begin{rmk}
Given two characters $\zeta, \zeta'$ lifting $\varepsilon^{-1} \det \rbar$, we have $\Rzp / \varpi \cong R^{\zeta'}_{\rbar} / \varpi$. Thus $\Rzp(\sigma) / \varpi \cong R^{\zeta'}_{\rbar}(\sigma) / \varpi$ if both characters are compatible with $\sigma$ (thus $\zeta = \zeta' \mu$ with $\mu$ an unramified charater). This implies that the two conjectures above are independent of the choice of $\zeta$.
\end{rmk}

\subsection{Local-global compatibility} \label{section:localglobal}
We now return to the global setting in Sect. \ref{section:patching}.

\subsubsection{Actions of Hecke algebras}
Let $\sigma$ be a representation of $K_p$ over $E$. Fix a $K_p$-stable $\O$-lattice $\sigma^{\circ}$ in $\sigma$. Let $\H(\sigma) = \End_{G_p}(\cInd^{G_p}_{K_p} \sigma)$ and $\H(\sigma^{\circ}) := \End_{G_p}(\cInd^{G_p}_{K_p} \sigma^{\circ})$, which is an $\O$-subalgebra of $\H(\sigma)$.

Since $\Minf$ is a pseudocompact $\Oinf \llbracket K_p \rrbracket$-module equipped with a compatible action of $G_p$, the $\Oinf$-module $\Minf(\sigma^{\circ}) := \sigma^{\circ} \otimes_{\O \llbracket K_p \rrbracket} \Minf$ has a natural action of $\H(\sigma^{\circ})$ commuting with the action of $\Rinf$ via isomorphisms
\[
(\sigma^{\circ} \otimes_{\O \llbracket K_p \rrbracket} \Minf)^d 
\cong \Homc_{\O \llbracket K_p \rrbracket}(\sigma^{\circ}, \Minf^d) 
\cong \Hom_{G_p}(\cInd^{G_p}_{K_p}(\sigma^{\circ}), (\Minf)^d),
\]
where the first isomorphism is induced by Schikhof duality and the second isomorphism is given by Frobenius reciprocity. In particular, $\Minf(\sigma^{\circ})$ is a $\O$-torsion free, profinite, linearly topological $\O$-module.

\subsubsection{Local-global compatibility}
We say a representation $\sigma$ of $K_p$ is a locally algebraic type if $\sigma = \otimes_{v|p} \sigma_v$, where $\sigma_v = \sigma(\lambda_v, \tau_v)$ or $\sigma^{cr}(\lambda_v, \tau_v)$ is a locally algebraic type of $\GL_2(F_v)$ for each $v|p$. We denote $\Rloc_p = \ctimes_{v|p} \Rfz_v$ and $\Rloc_p(\sigma) = \ctimes_{v|p} \Rfz_v(\sigma_v)$. Define $\Rloc(\sigma) = \Rloc \otimes_{\Rloc_p} \Rloc_p(\sigma)$, $\Rinf(\sigma) = \Rinf \otimes_{\Rloc_p} \Rloc_p(\sigma)$, $\Rinf'(\sigma)  = \Rinf' \otimes_{\Rloc_p} \Rloc_p(\sigma)$ and $\Rinvinf(\sigma)  = \Rinvinf \otimes_{\Rloc_p} \Rloc_p(\sigma)$.

\begin{lemma}  \label{lemma:Rinvinf}
\hfill
\begin{enumerate} 
\item There are $a_1, \cdots, a_t \in \m_{\infty}$ such that
\[ 
\Rinf(\sigma) = \frac{\Rinvinf(\sigma) \llbracket z_1 \rrbracket}{((1+z_1)^2 - (1+a_1))} \otimes_{\Rinvinf(\sigma)} \cdots \otimes_{\Rinvinf(\sigma)} \frac{\Rinvinf(\sigma) \llbracket z_t \rrbracket}{((1+z_t)^2 - (1+a_1))}. 
\]
In particular, $\Rinf(\sigma)$ is a free $\Rinvinf(\sigma)$-module of rank $2^t$.
\item Let $\p \in \Spec \Rinvinf(\sigma)$. The group $(\Gmc[2])^t(\O)$ acts transitively on the set of prime ideals of $\Rinf(\sigma)$ lying above $\p$.
\end{enumerate}
\end{lemma}

\begin{proof}
See \cite[Lemma 3.3]{MR3544298} for the first part and \cite[Lemma 3.4]{MR3544298} for the second part.
\end{proof}

\begin{prop} \label{prop:localglobalMinf}
\hfill
\begin{enumerate}
\item The action of $\Rinf$ on $\Minf(\sigma^{\circ})$ factors through $\Rinf(\sigma)$.
\item The action of $\H(\sigma)$ on $\Minf(\sigma^{\circ})[1/p]$ coincides with the composition
\[
\H(\sigma) \xrightarrow{\prod_{v|p} \phi_v} \Rloc_p(\sigma)[1/p] \rightarrow \Rinf(\sigma)[1/p],
\]
where $\phi_v$ is the map defined in Theorem \ref{thm:HeckeGalois}.
\item The module $\Minf(\sigma^{\circ})$ is finitely generated over $\Rinf(\sigma)$ and Cohen-Macaulay. Moreover, $\Minf(\sigma^{\circ})[1/p]$ is locally free of rank 1 over the regular locus of its support in $\Rinf(\sigma)[1/p]$.
\end{enumerate}
\end{prop}

\begin{proof}
This is an variance of \cite[Lemma 4.18, Theorem 4.19]{MR3529394}. The first assertion is an immediate consequence of local-global compatibility at $v | p$ at finite auxiliary levels. The second assertion follows from the first part and Theorem \ref{thm:HeckeGalois}. The first part of the third assertion is a consequence of numerical coincidence (cf. \cite[Lemma 3.5]{MR3544298}). The second part is due to \cite[Lemma 3.10]{MR3544298}. Note that the Hecke algebra in loc. cit. does not contain the Hecke algebra $U_{\varpi_{v_1}}$, thus their patched module is generically free of rank 2 instead of 1.
\end{proof}

\begin{definition}
It follows from Proposition \ref{prop:localglobalMinf} (3) that the support of $\Minf(\sigma^{\circ})[1/p]$ in $\Spec \Rinf(\sigma)[1/p]$ is a union of irreducible components, which we call the set of automorphic components of $\Spec \Rinf(\sigma)[1/p]$.
\end{definition}

\subsection{Breuil-M\'{e}zard via patching} \label{section:BMpatching}
Define $R^{S, \psi}_{\cS}(\sigma) = R^{S, \psi}_{\cS} \otimes_{\Rloc_p} \Rloc_p(\sigma)$ and $R^{\psi}_{\cS}(\sigma) = R^{\psi}_{\cS} \otimes_{\Rloc_p} \Rloc_p(\sigma)$.

\begin{prop} \label{prop:Stypepresentation}
For some $s \geq 0$, there is an isomorphism of $ \Rloc(\sigma)$-algebras
\[
R^{S, \psi}_{\cS}(\sigma) \cong \Rloc(\sigma) \llbracket x_1, \cdots, x_{s+ |S| -1} \rrbracket / (f_1, \cdots, f_s)
\]
for some elements $f_1, \cdots, f_s$. In particular, $\dim R^{S, \psi}_{\cS}(\sigma) \geq 4 |S|$ and $\dim \Rz_{\cS}(\sigma) \geq 1$.
\end{prop}

\begin{proof}
See \cite[Corollary 3.16]{MR3544298}.
\end{proof}

We define a Serre weight for $K_p$ to be an absolutely irreducible mod $p$ representations of $K_p = \prod_{v \in S_p} \GL_2(\O_{F_v}) \cong \prod_{v \in S_p} \GL_2(\Zp)$, which is of the form
\[
\bsigma_{a} = \otimes \bsigma_{a_v}
\]
with $\bsigma_{a_v}$ a Serre weight of $\GL_2(\O_{F_v})$ and $K_p$ acting on $\bsigma_{a}$ by reduction modulo $p$. 

For a Serre weight $\sigma_a$ for $K_p$, we write 
\begin{itemize}
    \item $\Minf^a := \Minf \otimes_{\O \llbracket K_p \rrbracket} \bsigma_a \cong \Homc_{\O \llbracket K_p \rrbracket}(\Minf, \bsigma_a^{\vee})^{\vee}$, which is an $\Rinf / \varpi$-module;
    \item $\mu_a'(\rhobar) := \frac{1}{2^t} e(\Minf^a, \Rinvinf / \varpi)$;
    \item $Z_a'(\rhobar) := \frac{1}{2^t} Z(\Minf^a)$ as a cycle on $\Rinvinf / \varpi$.
\end{itemize}
Suppose for each $v | p$, we have
\[
\overline{\sigma_v^{\circ}}
\xrightarrow{\sim} \oplus_{a_v} \bsigma_{a_v}^{m_{a_v}},
\]
then
\[
\overline{\sigma^{\circ}}
\xrightarrow{\sim} \oplus_{a} \bsigma_{a}^{m_a}
\]
with $m_a = \prod_v m_{a_v}$.

Due to \cite[Lemma 2.2.11]{MR2505297}, \cite[Lemma 4.3.9]{MR3292675}, \cite[Lemma 5.5.1]{MR3134019} and \cite[Proposition 3.17]{MR3544298}, we have the following equivalent conditions.

\begin{lemma} \label{lemma:BMequiv}
For any locally algebraic type $\sigma$, the following conditions are equivalent.
\begin{enumerate}
\item The support of $M(\sigma^\circ) \otimes_{\Zp} \Qp$ meets every irreducible component of $\Spec \Rloc(\sigma)[1/p]$.
\item $\Minf(\sigma^{\circ}) \otimes_{\Zp} \Qp$ is a faithful $\Rinf(\sigma)[1/p]$-module which is locally free of rank 1 over the regular locus of its support.
\item $R^{\psi}_{\cS}(\sigma)$ is a finite $\O$-algebra and $M(\sigma) \otimes_{\Zp} \Qp$ is a faithful $R^{\psi}_{\cS}(\sigma)[1/p]$-module.
\item $e(\Rinvinf(\sigma) / \varpi) = \sum_a m_a \mu'_a(\rhobar)$.
\item $Z(\Rinvinf(\sigma) / \varpi) = \sum_a m_a Z'_a(\rhobar)$.
\end{enumerate}
\end{lemma}

\begin{proof}
This is an analog of \cite[Proposition 3.17]{MR3544298} and \cite[Lemma 5.5.1]{MR3134019} in our setting.
\end{proof}

For each Serre weight $a_v$ ($\in \Z^2_+$) of $\GL_2(\O_{F_v})$, we have $M_{a_v} \otimes_{\O} k \cong \overline{\sigma}_{a_v}$ (see Notation for $M_{a_v}$). Define 
\[
\mu_{a_v}(\rhobar_v) = e(R_v^{a_v, \1, cr} / \varpi) \in \Z_{\geq 0}
\]
and
\[
\cC_{a_v}(\rhobar_v) = Z(R_v^{a_v, \1, cr} / \varpi)
\]
a 4-dimensional cycle of $\Spec \Rfz_v$. We obtain the following analogue of \cite[Theorem 5.5.2]{MR3134019}.

\begin{thm} \label{thm:BMequiv}
Suppose the equivalence conditions of Lemma \ref{lemma:BMequiv} hold for $\sigma = \otimes_{v|p} \sigma^{cr}(a_v, \1)$ with $a_v$ some Serre weights of $\GL_2(F_v)$. Then if $\sigma = \otimes_{v|p} \sigma_v$  is a locally algebraic type with $\sigma_v =  \sigma^*(\lambda_v, \tau_v)$ and $* \in \{\emptyset, cr\}$, and if we write
\[
\overline{\sigma^{\circ}} \xrightarrow{\sim} \oplus_{a} \bsigma_{a}^{m_a},
\]
then the following conditions are equivalent.
\begin{enumerate}
\item The equivalent conditions of Lemma \ref{lemma:BMequiv} hold for $\sigma$.
\item $e(R_v^{\lambda_v, \tau_v, *} / \varpi) = \sum_{a_v} m_{a_v} \mu_{a_v}(\rhobar_v)$ for each $v|p$.
\item $Z(R_v^{\lambda_v, \tau_v, *} / \varpi) = \sum_{a_v} m_{a_v} \cC_{a_v}(\rhobar_v)$ for each $v|p$.
\end{enumerate}
\end{thm}

\begin{proof}
Given Lemma \ref{lemma:BMequiv}, the proof of \cite[Theorem 5.5.2]{MR3134019} works verbatim in our setting.
\end{proof}

\subsection{The support at $v_1$}
Let $\sigma$ be a locally algebraic type for $G_p$. Suppose that $\Minf(\sigma^{\circ}) \neq 0$.

\begin{prop} \label{prop:suppv1}
The support of $\Minf(\sigma^{\circ}) \otimes_{\Zp} \Qp$ meets every irreducible component of $\Spec \Rfz_{v_1}[1/p]$.
\end{prop}

\begin{proof}
By assumption and Proposition \ref{prop:localglobalMinf} (3), $\Minf(\sigma^{\circ}) \otimes_{\Zp} \Qp$ is supported at an irreducible component $\cC$ of $\Spec \Rinf(\sigma)[1/p]$. We write $\cC_v$ for the corresponding irreducible component at $v \in S$. Let $\tilde{\cC}_{v_1}$ be an irreducible component of $\Spec \Rfz_{v_1}[1/p]$. It suffices to show that $\Minf(\sigma^{\circ}) \otimes_{\Zp} \Qp$ is supported at the irreducible component $\tilde{\cC}$ defined by $\{\cC_v \}_{v \in S - \{v_1\}}$ and $\tilde{\cC}_{v_1}$.

Choose a finite solvable totally real extension $F'$ of $F$ such that
\begin{itemize}
\item For each place $w$ of $F'$ above $v \in S_p$, $F'_w \cong F_v$;
\item For each place $w$ of $F'$ above $v_1$, the map $\Rfz_{w} \rightarrow \Rfz_{v_1}$ induced by restriction to $G_{F'_{w}}$ factors through $\Rur_{w}$.
\end{itemize}
Fix a place $w_1$ of $F'$ above $v_1$. Let $S' = S'_p \cup S'_{\infty} \cup \Sigma' \cup \{w_1\}$, where $S'_p$ is the set of places of $F'$ dividing $p$, $S'_{\infty}$ is the set of places of $F$ above $\infty$, and $\Sigma'$ is the set of places of $F'$ lying above $\Sigma$. Consider the following global deformation problems
\begin{align*}
\cR = &(\rhobar, S, \{\O\}_{v \in S} ,\{\D^{\cC_v}_v \}_{v \in S_p} \cup \{ \Dodd_v\}_{v \in \Sinf} \cup \{\DSt_v \}_{v \in \Sigma} \cup \{\D^{\tilde{\cC}_{v_1}}_{v_1} \}), \\
\cR' = &(\rhobar |_{G_{F'}}, S', \{\O\}_{w \in S'} ,\{\D^{\cC_w}_w \}_{w \in S_p'} \cup \{ \Dodd_w\}_{w \in \Sinf'} \cup \{\DSt_w \}_{w \in \Sigma'} \cup \{\Dur_{w_1} \}),
\end{align*}
where $\cC_w$ is the image of $\cC_v$. We claim that $\Rz_{\cR'}$ is a finite $\O$-algebra. Given this, since the morphism $\Rz_{\cR'} \rightarrow \Rz_{\cR}$ is finite by Proposition \ref{prop:finitedefbc}, $\Rz_{\cR}$ is a finite $\O$-module. On the other hand, $\Rz_{\cR}$ has a $\bQp$-point since it has Krull dimension at least 1 by Proposition \ref{prop:Stypepresentation}. This gives a lifting $\rho$ of $\rhobar$ of type $\cR$. Since $\rho|_{G_{F'}}$ lies in the automorphic component defined by $\cC$ restricted to $F'$, we obtain that  $\rho$ is automorphic by solvable base change. It follows that $\rho$ gives a point on $\tilde{\cC}$ and the theorem is proved.

To prove the claim, we denote the patched module constructed in the same way as $\Minf$ replacing $F$ with $F'$, $S$ with $S'$ and $v_1$ with $w_1$ by $\Minf'$, which is endowed with an $\Oinf'$-linear action $\Rinf'$. Note that by our assumption, the local deformation problem at $v_1$ (resp. $w_1$) of $\cS$ (resp. $\cS'$) is the Taylor-Wiles deformation defined in Sect. \ref{section:TWdef} and thus each irreducible component of $R_{v_1}$ (resp. $R_{w_1}$) can be realized by the level (pro-$v_1$ Iwahori) we choose in the patching process.

Write $\a'$ for the ideal of $\Oinf'$ defined by its formal variables, $\cS'$ for corresponding global deformation problem (as in Sect. \ref{section:completedcohomology}) and $\sigma'$ for the locally algebraic type defined by $\sigma$ restricting to $F'$. It follows that $\Minf'(\Sigma^{\prime, {\circ}}) \otimes_{A^{S'}_{\cS'}} A^{S'}_{\cR'}$ is a faithful $\Rinf'(\sigma') \otimes_{A^{S'}_{\cS'}} A^{S'}_{\cR'}$-module by Proposition \ref{prop:localglobalMinf} (3) and the irreducibility of $\Spec \Rinf'(\sigma') \otimes_{A^{S'}_{\cS'}} A^{S'}_{\cR'}$ (which is an automorphic component of $\Spec \Rinf'(\sigma')$). Thus $\Rz_{\cR'} \cong (\Rinf'(\sigma') \otimes_{A^{S'}_{\cS'}} A^{S'}_{\cR'}) / \a' (\Rinf'(\sigma') \otimes_{A^{S'}_{\cS'}} A^{S'}_{\cR'})$ is a finite $\O$-algebra by the same reason as in the proof of Lemma \ref{lemma:BMequiv}.
\end{proof}
\section{Patching and $p$-adic local Langlands correspondence} \label{section:patchpLLC}
Throughout this section, we will use freely the notations in Sect. \ref{section:patchingargument} and Sect. \ref{section:patchBM}. We fix a place $\p$ of $F$ lying above $p(=2)$. Let $G = \GL_2(F_\p) \cong \GL_2(\Qp)$, $K = \GL_2(\O_{F_\p}) \cong \GL_2(\Zp)$, $T$ be the subgroup of diagonal matrices in $G$, and $T_0$ be the subgroup of diagonal matrices in $K$.

\subsection{Patching and Banach space representations} \label{section:patchBanach}
For each place $v \neq \p$ above $p$, we fix a locally algebraic type $\sigma_v$ compatible with $\psi$ and an irreducible component $\cC_v$ of the corresponding deformation ring $R_v^{\lambda_v, \tau_v, *}$, where $* \in \{ss, cr\}$. Write $\sigma^{\p} = \otimes_{v \in S_p - \{\p\}} \sigma_v$, which is a representation of $K^{\p} = \prod_{v \in S_p - \{\p\}} \GL_2(\O_{F_v})$. 

We denote $R^{\loc, \p} = \ctimes_{\O, v \in S_p - \{\p\}} \Rfz_v \ctimes_{\O, v \in S - S_p} R_v$, $R^{\loc, \p}(\sigma^{\p}) = \ctimes_{\O, v \in S_p - \{\p\}} R_v^{\lambda_v, \tau_v, *}\ctimes_{\O, v \in S - S_p} R_v $ and $R^{\loc, \p}(\cC^{\p}) = \ctimes_{\O, v \in S_p - \{\p\}} R_v^{\cC_v} \ctimes_{\O, v \in S - S_p} R_v$, where $R_v$ is the local deformation ring at $v$ defined by the global deformation problem $\cS$ in Sect. \ref{section:completedcohomology}. Define
\[
\tMinf' := \Minf \otimes_{\O \llbracket K^{\p} \rrbracket} (\sigma^{\p})^{\circ}
\]
and
\[
\tMinf := \tMinf' \otimes_{R^{\loc, \p}} R^{\loc, \p}(\cC^{\p}).
\]
Thus $\tMinf'$ is an $\Oinf \llbracket K \rrbracket$-module endowed with an $\Oinf$-linear action of 
\[
\tRinf' := \Rinf \otimes_{R^{\loc, \p}} R^{\loc, \p}(\sigma^{\p}),
\]
which is free over $\tRinvinfp := \Rinvinf \otimes_{R^{\loc, \p}} R^{\loc, \p}(\sigma^{\p})$ of rank $2^t$ (Lemma \ref{lemma:Rinvinf} (1)). Similarly, $\tMinf$ is an $\Oinf \llbracket K \rrbracket$-module endowed with an $\Oinf$-linear action of 
\[
\tRinf := \Rinf \otimes_{R^{\loc, \p}} R^{\loc, \p}(\cC^{\p}),
\]
which is free over $\tRinvinf = \Rinvinf \otimes_{R^{\loc, \p}} R^{\loc, \p}(\sigma^{\p})$ of rank $2^t$. Assume that $\tMinf[1/p]$ is non-zero.

\begin{rmk}
The assumption is satisfied when $\rhobar$ admits an automorphic lift $\rho$ whose associated local Galois representation $\rho|_{G_{F_v}}$ lies on $\cC_v$ for each $v \in S_p - \{\p\}$, $\rho|_{G_{F_v}}$ is of Steinberg type for each $v \in \Sigma$ and is unramified away from $S$ since the corresponding automorphic form is a specialization of $\tMinf$.
\end{rmk}

The following proposition is a direct consequence of Proposition \ref{prop:Minfprop} (3).

\begin{prop}
$\tMinf'$ is finitely generated over $\Oinf \llbracket K \rrbracket$ and projective in the category $\Modpro_{K, \zeta}(\Oinf)$, with $\zeta = \psi|_{\O_{F_\p}^{\times}}$. In particular, it is finitely generated over $\tRinf' \llbracket K \rrbracket$ and projective in $\Modpro_{K, \zeta}(\O)$.
\end{prop}

\begin{rmk}
$\tMinf'$ is the same as the patched module considered in \cite{MR3529394}.
\end{rmk}

Let us denote by $\Piinf := \Homc_\O(\tMinf', E)$. If $y \in \mSpec \tRinf'[1/p]$, then we have
\begin{gather*}
\Pi_y := \Homc_\O(\tMinf' \otimes_{\tRinf', y} E_y, E) = \Piinf[\m_y]
\end{gather*}
is an admissible unitary $E$-Banach space representation of $\GL_2(L)$ (by \cite[Proposition 2.13]{MR3529394}). The composition $\Rfz_{\p} \rightarrow \Rinf \xrightarrow{y} E_y$ defines an $E_y$-valued point $x \in \Spec \Rfz_{\p}[1/p]$ and thus a continuous representation $r_x : G_{\Q_2} \rightarrow \GL_2(E_y)$.

\begin{prop} \label{prop:genericlalg}
Let $y \in \mSpec \tRinf'[1/p]$ be a closed $E$-valued point whose the associated local Galois representation $r_x$ is potentially semi-stable of type $\sigma_{\p}$. Assume that $y$ lies on an automorphic component of $\Rinf(\sigma)$ with $\sigma = \sigma_{\p} \otimes \sigma^{\p}$ and $\pism(r_x)$ is generic. Then
\[
\Pi_y^{\lalg} \cong \pi_{\lalg}(r_x).
\]
\end{prop}

\begin{proof}
The proof of \cite[Theorem 4.35]{MR3529394} ($r_x$ potentially crystalline) and \cite[Theorem 7.7]{2018arXiv180301610P} ($r_x$ potentially semi-stable) works verbatim in our setting.
\end{proof}

\subsection{Patched eigenvarieties}
We write $R_{\1}$ for the universal deformation ring of the trivial character $\1: G_{\Q_2} \rightarrow k^{\times}$ and $\1^{\univ}$ for the universal character. Via the natural map $\O[Z] \rightarrow R_{\1}[Z]$, the maximal ideal of $R_{\1}[Z]$ generated by $\varpi$ and $z - \1^{\univ} \circ \Art_{L}(z)$ gives a maximal ideal of $\O[Z]$. If we denote by $\Lambda_Z$ the completion of the group algebra $\O[Z]$ at this maximal ideal, then the character $\1^{\univ} \circ \Art_{L}$ induces an isomorphism $\Lambda_Z \xrightarrow{\sim} R_{\1}$.

We define the patched eigenvarieties following \cite[\S 3]{MR3623233} and \cite[\S 6]{2018arXiv180906598E}. Denote $\Rfsign_{\p}$ the quotient corresponding to the irreducible component of $\Spec \Rf_{\p}$ given by $\psi(\Art_{\Q_2}(-1))$ (see Sect. \ref{section:irredcomp}). 

We define $\tAinf'$ (resp. $\tAinvinfp$, $\tAinf$ and $\tAinvinf$) in the same way $\tRinf'$ (resp. $\tRinvinfp$, $\tRinf$ and $\tRinvinf$) is defined in Sect. \ref{section:patching} and Sect. \ref{section:patchBanach}, but by replacing $\Rfz_{\p}$ with $\Rfsign_{\p}$ at $\p$ (and keeping all other places unchanged). Let $\Xinf := \Spf(\tAinvinfp)^{\rig}$, $\cX_{\p} = \Spf (\Rf_{\p})^{\rig}$, $\cX^{\p} = \Spf (R^{\loc, \p}(\sigma^{\p}))^{\rig}$ so that
\[
\Xinf = \cX_{\p} \times \cX^{\p} \times \cU^g,
\]
where $\cU := \Spf (\O_E \llbracket x \rrbracket)^{\rig}$ is the open unit disk over $E$.

We define $\tNinf = \tMinf' \ctimes_{\O} \1^{\univ}$ and $\tPiinf = \Hom(\tNinf, E)$, both of which are equipped with an $\tAinvinfp$-action (resp. $\tAinf'$-action) via $\tAinvinfp \rightarrow \tRinvinfp \ctimes_{\O} R_\1$ (resp. $\tAinf' \rightarrow \tRinf' \ctimes_{\O} R_\1$) induced by $\Rfsign_\p \rightarrow \Rfz_{\p} \ctimes_{\O} R_{\1}$ in Sect. \ref{section:irredcomp}. Note that $\GL_2(\Q_2)$ acts on $\1^{\univ}$ via $\GL_2(\Q_2) \xrightarrow{\det} \Q_2^\times \rightarrow \Lambda_Z^{\times} \xrightarrow{\sim} R_{\1}^{\times}$ and thus on $\tNinf$ diagonally, which commutes with the action of $\tAinvinfp$ (resp. $\tAinf'$).

\begin{prop} \label{prop:tMinfprop}
Let $K'$ be the open normal subgroup of $K$ defined by $\{g = sz \mid s \in \SL_2(\Z_2), \ s \equiv (\begin{smallmatrix} 1 & 0 \\ 0 & 1 \end{smallmatrix}) \text{ mod } 4, \ z \in 1 + 4 \Z_2 \}$. Then $\tNinf$ is projective in the category $\Modpro_{K'}(\O)$.
\end{prop}

\begin{proof}
Using the decomposition $K' \cong (K' / K' \cap Z) \times \Gamma$ as in the proof of Proposition \ref{prop:Minfprop}, the proof of \cite[Proposition 6.10]{MR3732208} works verbatim in our setting.
\end{proof}

Let $\hT$ be the rigid analytic space over $E$ parametrizing continuous characters of $T$ and $\hT^0$ be the rigid analytic space over $E$ parametrizing continuous characters of $T_0$. Define the patched eigenvariety $\Xtrinf$ as the support of the coherent $\O_{\Xinf \times \hT}$-module 
\[
J_B(\tPiinf^{\tAinf'-\an})'
\]
on $\Xinf \times \hT$, where $J_{B}$ is Emerton's Jacquet functor with respect to $B$ defined in \cite{MR2292633}, $\tPiinf^{\tAinf'-\an}$ is the subspace of $\tAinf'$-analytic vectors defined in \cite[Definition 3.2]{MR3623233}, and $'$ is the strong dual. This is a reduced closed analytic subset of $\Xinf \times \hT$ \cite[Corollary 3.20]{MR3623233} whose points are
\[
\{ x = (y, \delta) \in \Xinf \times \hT \ | \Hom_{T}(\delta, J_{B}\big(\tPiinf^{\tAinf'-\an}[\p_y] \otimes_{E_y} E_x) \big) \neq 0 \}
\]
with $\p_y \subset \tAinf'$ the prime ideal corresponding to the point $y \in \Xinf$ and $E_y$ the residue field of $\p_y$. 

Let $\W_{\infty} = \Spf(\O_{\infty})^{\rig} \times \hT^0$ be the weight space of the patched eigenvariety. We define the weight map $\omega_X : \Xtrinf \rightarrow \W_{\infty}$ by the composite of the inclusion $\Xtrinf \rightarrow \Xinf \times \hT$ with the map from $\Xinf \times \hT$ to $\Spf(\O_{\infty})^{\rig} \times \hT^0$ induced by the $\Oinf$-structure of $\tRinf$ and by the restriction $\hT \rightarrow \hT^0$.

\begin{prop}
The rigid analytic space $\Xtrinf$ is equidimensional of dimension $q + 4|S| + 1$ and has no embedded component.
\end{prop}

\begin{proof}
The proof of \cite[Proposition 3.11]{MR3623233}, which shows that the weight map $\omega_X$ is locally finite, works verbatim in our setting. Thus the dimension of $\Xtrinf$ is equal to the dimension of $\W_{\infty}$, which is given by
\begin{align*}
\dim \W_{\infty} &= \dim \Spf(\O_{\infty})^{\rig}  + \dim \hT^0 \\
&= q + 4|S| - 1 + 2.
\end{align*}
\end{proof}

Let $\iota$ be an automorphism of $\hT$ given by 
\[
\iota(\delta_{v, 1}, \delta_{v, 2}) = (\unr(q) \delta_{v, 1}, \unr(q^{-1}) \delta_{v, 2} (\cdot)^{-1}),
\]
which induces an isomorphism of rigid spaces
\begin{align*}
\Xinf \times \hT &\xrightarrow{\sim} \Xinf \times \hT \\
(x, \delta) &\mapsto (x, \iota^{-1}(\delta) ),
\end{align*}
and thus a morphism of reduced rigid spaces over $E$:
\[
\Xtrinf \rightarrow \Xtri_\p \times \cX^{\p} \times \cU^g,
\]
where $\Xtri_\p$ is the space of trianguline deformation of $\rhobar |_{G_{F_\p}}$ \cite[Definition 2.4]{MR3623233}.

\begin{thm} \label{thm:lgtri}
This morphism induces an isomorphism from $\Xtrinf$ to a union of irreducible components of $\Xtri_\p \times \cX^{\p} \times \cU^g$.
\end{thm}

\begin{proof}
This can be proved in the same way as in \cite[Theorem 3.21]{MR3623233}.
\end{proof}

\begin{prop} \label{prop:suppMinfp}
The support of $\tNinf$ in $\Spec \tAinf'$ is equal to a union of irreducible components in $\Spec \tAinf'$.
\end{prop}

\begin{proof}
Replacing \cite[Theorem 3.21]{MR3623233} with Theorem \ref{thm:lgtri}, the proof of \cite[Theorem 6.3]{2018arXiv180906598E} works verbatim in our setting.
\end{proof}

\begin{corollary} \label{cor:fixtypefaithful}
Let $\Sigma_{ps}$ be the set of principal series types. Then the Zariski closure in $\Spec \tAinf'$ of the set of points having types $\sigma \in \Sigma_{ps}$ and lying in the support of $\tNinf(\sigma) := \tNinf \otimes_{\O \llbracket K \rrbracket} \sigma$ is equal to a union of irreducible components of $\Spec \tAinf'$.
\end{corollary}

\begin{proof}
Since $\tNinf$ is projective in $\Modpro_{K'}(\O)$ by Proposition \ref{prop:tMinfprop}, it is captured by the family of principal series types by \cite[Proposition 3.11]{2018arXiv180906598E}.  Applying proposition \cite[Proposition 2.11]{2018arXiv180906598E} to $M = \tNinf$ and $R = \tAinf' / \Ann_{\tAinf'}(\tNinf)$, we see that the set of points having principal series types are dense in $\tAinf' / \Ann_{\tAinf'}(\tNinf)$, which is equal to a union of irreducible components in $\Spec \tAinf'$ by Proposition \ref{prop:suppMinfp}. This proves the corollary.
\end{proof}

\subsection{Relations with Colmez's functor}
\begin{lemma} \label{lemma:Minfadm}
$\tMinf$ lies in $\fC(\O)$.
\end{lemma}

\begin{proof}
This follows immediately from Proposition \ref{prop:Minfprop} (3).
\end{proof}

As a result, we may apply Colmez's functor $\cV$ to $\tMinf$ and obtain an $\tRinf \llbracket \GQp \rrbracket$-module $\cV(\tMinf)$.

\begin{prop} \label{prop:Minffg}
$\cV(\tMinf)$ is finitely generated over $\tRinf \bracketGQp$.
\end{prop}

\begin{proof}
Using Proposition \ref{prop:fgadm}, the proof of \cite[Proposition 3.4]{2018arXiv180307451T} works without any change.
\end{proof}

Let $\sigma$ be a locally algebraic type for $G$. We define $\tRinf(\sigma) = \tRinf \otimes_{\Rfz_{\p}} \Rfz_{\p}(\sigma)$ (resp. $\tRinf'(\sigma) = \tRinf' \otimes_{\Rfz_{\p}} \Rfz_{\p}(\sigma)$) and $\tMinf(\sigma^{\circ}) = \tMinf \otimes_{\O \llbracket K \rrbracket} \sigma^{\circ}$ (resp. $\tMinf'(\sigma^{\circ}) = \tMinf' \otimes_{\O \llbracket K \rrbracket} \sigma^{\circ}$), which satisfies a similar local-global compatibility as in Sect. \ref{section:localglobal}.

\begin{thm} \label{thm:MinfCH}
The action of $\tRinf \llbracket \GQp \rrbracket$ on $\cV(\tMinf)$ factors through $\tRinf \bracketGQp / J$, where $J$ is a closed two-sided ideal generated by $g^2-\tr\big(\rinf(g)\big)g+\det\big(\rinf(g)\big)$ for all $g \in \GQp$, where $\rinf : \GQp \rightarrow \GL_2(\tRinf)$ is the Galois representation lifting $\rbar$ induced by the natural map $\Rfz_{\p} \rightarrow \tRinf$.
\end{thm}

\begin{proof}
The proof of \cite[Theorem 3.7]{2018arXiv180307451T} works verbatim in our setting.
\end{proof}

\begin{corollary} \label{cor:VMfinite}
$\cV(\tMinf)$ is finitely generated over $\tRinf$.
\end{corollary}

\begin{proof}
See \cite[Corollary 3.8]{2018arXiv180307451T}.
\end{proof}

\begin{prop} \label{prop:RinfVMfaithful}
$\tRinf [1/p]$ acts on $\cV(\tMinf) [1/p]$ nearly faithfully, i.e. $\Ann_{\tRinf [1/p]}(\cV(\tMinf) [1/p])$ is nilpotent.
\end{prop}

\begin{proof}
Consider $V := \cV(\tMinf) \ctimes_{\O} 1^{\univ}$, which is an $\tAinf$-module (resp. $\tAinvinf$-module) via $\tAinf \rightarrow \tRinf \ctimes_{\O} R_1$ (resp. $\tAinvinf \rightarrow \tRinf \ctimes_{\O} R_1$) induced by the homomorphism $\Rfsign_\p \rightarrow \Rfz_{\p} \ctimes_{\O} R_{1}$ in Sect. \ref{section:irredcomp}. Note that irreducible components of $\Spec \tAinvinf$ are in bijection with irreducible components of $\Spec \Rfz_{v_1}$ if $\O$ is sufficiently large (in the sense that all irreducible components of local deformation rings are geometrically irreducible, see \cite[Appendix A]{hu2019crystabelline}). By Corollary \ref{cor:fixtypefaithful}, the set of points in $z \in \mSpec \tAinf[1/p]$ with a principal series types $\sigma$ lying in the support of $\tNinf(\sigma)$ are dense in a union of irreducible components of $\Spec \tAinf[1/p]$, which is equal to $\Spec \tAinf[1/p]$ by Lemma \ref{lemma:Rinvinf} (2) and Proposition \ref{prop:suppv1}.

On the other hand, for any point $z \in \mSpec \tAinf[1/p]$ as above, there is a $x \in \mSpec \tRinf \ctimes_{\O} R_1[1/p]$ lying in the preimage of $z$ satisfying $(\tMinf)_y \neq 0$, where $y \in \mSpec \tRinf[1/p]$ is the point given by $x$. Note that the point $y$ is also of principal series type. It follows that $\cV(\tMinf)_y \neq 0$ by Proposition \ref{prop:genericlalg} ($\Pi^{\lalg}_y \cong \pi_{\lalg}$), \cite[Theorem 4.3.1]{MR2642406} and \cite[Proposition 2.2.1]{MR2667890} ($\cV(\widehat{\pi_{\lalg}}) \neq 0$), which implies that $V_z \neq 0$. Hence $\tAinf[1/p]$ acts on $V[1/p]$ nearly faithfully.

Note that $V$ admits two actions of $R_1$, one via $R_1 \rightarrow \Rfz_{\p} \ctimes_{\O} R_{1}$ given by $(r, \chi) \mapsto \chi^2$ and the other via $R_1 \rightarrow \Rfsign_{\p}$ given by $r \mapsto (\zeta \varepsilon)^{-1}\det r$, which are compatible by the following commutative diagram
\[
\begin{tikzcd}
&R_1  \arrow[r, "s"] \arrow[d] &R_1 \arrow[d] \\
&\Rfsign_{\p}  \arrow[r] &\Rfz_{\p} \ctimes_{\O} R_{1},
\end{tikzcd}
\]
where $s$ is the map induced by $\chi \mapsto \chi^2$. Denote $\iota: R_1 \rightarrow \O$ the homomorphism given by the trivial lifting of $1$. It induces the following commutative diagram
\[
\begin{tikzcd}
&\Rfsign_{\p}  \arrow[r] \arrow[d, "\otimes_{R_1 , \iota} \O"] &\Rfz_{\p} \ctimes_{\O} R_{1} \arrow[d, "\otimes_{R_1 , \iota}  \O"] \\
&\Rfz_{\p} \arrow[r, equal] &\Rfz_{\p}
\end{tikzcd}
\]
and thus an $\tRinf$-module isomorphism $V \otimes_{R_1 , \iota} \O \cong \cV(\tMinf)$ (for both $R_1$-actions because $\iota \cong \iota \circ s$). Denote $I$ the kernel of the homomorphism $\tAinf \rightarrow \tRinf$ induced by $\iota$. Since $V$ is finite over $\tAinf$ ($V$ is finite over $\tRinf \ctimes_\O R_1$ by Corollary \ref{cor:VMfinite} and $\tRinf \ctimes_\O R_1$ is finite over $\tAinf$ by Proposition \ref{prop:localetale}), we see that $\cV(\tMinf) [1/p] \cong V / I V [1/p]$ is a nearly faithful $\tRinf [1/p] \cong \tAinf / I \tAinf [1/p]$-module by \cite[Lemma 2.2]{MR2470688}. This finishes the proof.
\end{proof}

\label{Cor:nonvanishing}
\begin{corollary}
For all $y \in \Spec \tRinf[1/p]$, we have $\cV(\Pi_y) \neq 0$. In particular, $\Pi_y \neq 0$.
\end{corollary}

\begin{proof}
See \cite[Corollary 3.10]{2018arXiv180307451T}.
\end{proof}

\begin{thm} \label{thm:absirredMinf}
For $y \in \mSpec \tRinf[1/p]$ whose associated Galois representation $r_x$ is absolutely irreducible, we have $\cV(\Pi_y) \cong r_x^{\oplus n_y}$ for some integer $n_y \geq 1$. In particular, $\Minf(\sigma^{\circ})[1/p]$ is supported on every non-ordinary (at $\p$) component of $\Rinf(\sigma)[1/p]$ for each locally algebraic type $\sigma$ for $G$.
\end{thm}

\begin{proof}
The proof of \cite[Theorem 4.1]{2018arXiv180307451T} works verbatim in our setting with Corollary 3.10 in loc. cit. replaced by Corollary \ref{Cor:nonvanishing}.
\end{proof}

\begin{corollary}
If moreover $r_x$ is potentially semi-stable except possibly in the following cases:
\begin{itemize}
    \item $\lambda = (a, b)$ with $a+b$ odd, $\tau = \eta \oplus \eta$, and $\pism(r_x)$ is non-generic;
    \item $\lambda = (a, b)$ with $a+b$ even, $r_x \otimes \chi$ is potentially crystalline of inertial type $\eta \oplus \eta$ with $\pism(r_x \otimes \chi)$ is non-generic, where $\chi= \sqrt{\pr(\varepsilon)}$ and $\pr: \O^{\times} \rightarrow 1 + \varpi \O$ given by projection,
\end{itemize}
then we have $n_y = 1$. In particular, $n_y = 1$ in an open dense subset of $\mSpec \tRinf[1/p]$.
\end{corollary}

\begin{proof}
Replacing Proposition 2.7 in \cite{2018arXiv180307451T} with Proposition \ref{prop:equivB}, the proof of Corollary 4.2 in loc. cit. works verbatim in our setting.
\end{proof}

\section{Patching argument: ordinary case} \label{section:ordpatch}
The goal of this section is to construct automorphic points on some partially ordinary irreducible components of $\Rinf(\sigma)$. We will follow the strategy in \cite{MR3252020, MR3327536, MR3904451, Sasaki2018II} and use freely the notations in Sect. \ref{section:QMF}. 

Let $p=2$ and $F$ be a totally real field ($p$ may not split completely). If $v$ is a finite place of $F$ above $2$ and $c \geq b \geq 0$ are integers, then we define an open compact subgroup $\Iw_v(b,c)$ of $\GL_2(\O_{F_v})$ by the formula
\[
\Iw_v(b,c) = \bigg\{ \begin{pmatrix} t_1 & * \\ 0 & t_2 \end{pmatrix} \ \text{mod } \varpi_v^c \mid t_1 \equiv t_2 \equiv 1 \ \text{mod } \varpi_v^b \bigg\}.
\]
Thus $\Iw_v(0,1)$ is the Iwahori subgroup of $\GL_2(\O_{F_v})$ and $\Iw_v(1, 1)$ is the pro-$v$ Iwahoric subgroup.  

Let $U_v = \Iw_v(b, c)$ for some integers $c \geq b \geq 1$. We define the operator $\U_{\varpi_v}$ by the double coset operator $\U_{\varpi_v} = [U_v ( \begin{smallmatrix} \varpi_w  & 0 \\ 0 & 1 \end{smallmatrix} ) U_v]$, and the diamond operator $\langle \alpha \rangle = [U_v ( \begin{smallmatrix} \alpha  & 0 \\ 0 & 1 \end{smallmatrix} ) U_v]$ for $\alpha \in \O_{F_v}^{\times}$.

\begin{lemma} \label{lemma: Iwahoricommute}
Let $v$ be a fixed place of $F$ above $p$. If $U' \subset U$ are open compact subgroups of $G(\A_F^{\infty})$ such that $U'_w = U_w$ if $w \neq v$, and $U_v' = \Iw_v(b', c') \subset U_v = \Iw_v(b, c)$ for some $b' \geq b \geq 1$, $c' \geq c$. Then for any topological $\O$-algebra $A$, the operators $\U_{\varpi_v}$ and $\langle \alpha \rangle$ for $\alpha \in \O_{F_v}^{\times}$ commute with each other and with the natural map
\[
S_{\psi}(U, A) \rightarrow S_{\psi}(U', A).
\]
\end{lemma}

\begin{proof}
See \cite[\S 1]{MR1097614}.
\end{proof}

\subsection{Partial Hida families} \label{section:Hidafamily}
Let $S = S_p \cup S_\infty \cup \Sigma \cup \{v_1\}$ be a set defined as in Sect. \ref{section:QMF}. Let $P \subset S_p$ be a subset. For each $v \in S_p - P$, we fix a locally algebraic type $\sigma_v$ compatible with $\psi$.  Define the open compact subgroup $U^P = \prod_{v} U_v$ of $(D \otimes_{F} \A^{\infty, P}_F)^{\times}$ by
\begin{itemize}
\item $U_v = (\O_D)^{\times}_v$ if $v \notin S$ or $v \in \Sigma \cup (S_p - P)$.
\item $U_{v_1}$ is the pro-$v_1$ Iwahori subgroup.
\end{itemize}
If $c \geq b \geq 0$ are two integers, then we set $U(b,c) = U^P \times \prod_{v \in P} \Iw_v(b, c)$. Let $\sigma^P(b, c) = \otimes_{v \in S_p - P} \sigma_v \bigotimes \otimes_{v \in P} 1$ be a continuous representation of $\prod_{v \in S_p -P} U_v \times \prod_{v \in P} \Iw_v(b, c)$. We will write $S_{\sigma^P, \psi}(U(b, c), \O)$ for $S_{\sigma^P(b, c), \psi}(U(b,c), \O)$.

We define $\O_P^{\times}(b,c) = \{ t \in (\O_{F_v} / \varpi_v^c)^{\times} | t \equiv 1 \ \text{mod } \varpi_v^b \}$. The group $U(0,c)$ acts on $S_{\sigma^P, \psi}(U(b,c), \O)$, which is uniquely determined by the diamond operator action of $\O_P^{\times}(0,c)$ via the embedding
\[
\O_P^{\times}(0,c) / \O_P^{\times}(b,c) \rightarrow U(0,c) /  U(b,c) \quad (y_v)_{v \in P} \ \text{mod } \O_P^{\times}(b,c) \mapsto \Big( (\begin{smallmatrix} y_v & 0 \\ 0 & 1 \end{smallmatrix}) \Big)_{v \in P} \ \text{mod } U(b,c).
\]
We define $\Lambda_P(b,c) = \O[\O_P^{\times}(0,c) / \O_P^{\times}(b,c)]$ and $\Lambda^b_P = \varprojlim_c \Lambda_P(b,c)$. If $b=1$, we write $\Lambda_P$ for $\Lambda^1_P$.

We write $\T^{\ord}_{S, P}$ for the polynomial algebra over $\Lambda_P[\Delta_{v_1}]$ in the indeterminates $T_v, S_v$ for $v \notin S$ and the indeterminates $\U_{\varpi_{v}}$ for $v \in P \cup \{v_1\}$. Define a $\T^{\ord}_{S, P}$-module structure on $S_{\sigma^P, \psi}(U(b, c), \O)$ by letting $\Lambda_P[\Delta_{v_1}]$ act via diamond operators and $T_v, S_v, \U_{\varpi_{v}}$ act as usual. Since for $v \in P$ the operators $U_{\varpi_v}$ and $\langle \alpha \rangle$ commutes with all inclusions $S_{\sigma^P, \psi}(U(b, c), \O) \rightarrow S_{\sigma^P, \psi}(U(b', c'), \O)$ for every $b' \geq b \geq 1$, $c' \geq c$, these maps become maps of $\T^{\ord}_{S, P}$-modules.

Denote $\U = \U_P := \prod_{v \in P} \U_{\varpi_v}$, it follows that $e =  \lim_{n \to \infty} (\U_P)^{n !}$ defines an idempotent in $\End_{\O}(S_{\sigma^P, \psi}(U(b,c), \O))$ (resp. $\End_{\O / \varpi^s}(S_{\sigma^P, \psi}(U(b,c), s))$) (c.f. \cite[Lemma 2.10]{MR3702498}). Define the ordinary subspace of $S_{\sigma^P, \psi}(U(b,c), \O)$ (resp.
$S_{\sigma^P, \psi}(U(b,c), s)$) by
\[
\Sord(U(b, c), \O) = e S_{\sigma^P, \psi}(U(b, c), \O) \quad (\text{resp. } \Sord(U(b, c), s) = e S_{\sigma^P, \psi}(U(b, c), s)).
\]

\begin{lemma}
For all $c \geq b \geq 1$, the natural map
\[
\Sord(U(b,b), \O) \rightarrow \Sord(U(b,c), \O)
\]
is an isomorphism.
\end{lemma}

\begin{proof}
See \cite[Lemma 2.3.2]{MR3252020} and \cite[Lemma 2.5.2]{MR2941425}.
\end{proof}

We now define the partial Hida family. By Lemma \ref{lemma: Iwahoricommute}, for $c' \geq c$ the natural maps
\[
S_{\psi}(U(c, c), \O) \rightarrow S_{\psi}(U(c',c'), \O)
\]
commute with the action of the Hecke operator $\U_P$ and $\langle \alpha \rangle$, $\alpha \in \O_P^{\times}(p)$.

\begin{definition}
We define
\[
\Mord(U^P) = \varprojlim_{c} \Sord(U(c,c), \O)^d,
\]
which is naturally a $\Lambda_P$-module.
\end{definition}

\begin{prop} {\ }
\begin{enumerate}
\item For every $s, c \geq 1$, there is an isomorphism
\[
\Mord(U^P) \otimes_{\Lambda_P} \Lambda_P(1,c) / (\varpi^s) \xrightarrow{\sim} \Sord(U(c,c), s)^{\vee}.
\]
\item For every $c \geq 1$, the $\Lambda_P^c$-module $\Mord(U^P)$ is finite free of rank equal to the $\O$-rank of $\Sord(U(c,c), \O)$.
\end{enumerate}
\end{prop}

\begin{proof}
See \cite[Proposition 2.3.3]{MR3252020}.
\end{proof}

The algebra $\T^{\ord}_{S, P}$ acts naturally on $\Sord(U(c,c), s)$. We write $\T^{S, \ord}_{\psi}(U(c,c), \O)$ for its image in $\End_{\Lambda_P}(\Sord(U(c,c), \O))$.

\begin{definition}
We define
\[
\T^{S, \ord}_{\psi}(U^P) := \varprojlim_{c} \T^{S, \ord}_{\psi}(U(c,c), \O)
\]
endowed with inverse limit topology. It follows immediately from the definition that $\T^{S, \ord}_{\psi}(U^P)$ acts on $\Mord(U^P)$ faithfully.
\end{definition}

\begin{lemma}
$\T^{S, \ord}_{\psi}(U^P)$ is a finite $\Lambda_P$-algebra with finitely many maximal ideals. Denote its finitely many maximal ideals by $\m_1, \cdots, \m_r$ and let $J = \cap_{i} \m_i$ denote the Jacobson radical. Then $\T^{S, \ord}_{\psi}(U^P)$ is $J$-adically complete and separated, and we have
\[
\T^{S, \ord}_{\psi}(U^P) = \T^{S, \ord}_{\psi}(U^P)_{\m_1} \times \cdots \times \T^{S, \ord}_{\psi}(U^P)_{\m_r}.
\]
For each $i$, $\T^{S, \ord}_{\psi}(U^P) / \m_i$ is a finite extension of $k$.
\end{lemma}

\begin{proof}
The proof is identical to Lemma \ref{lemma:bigHeckemax}.
\end{proof}

Let $\m \subset \T^{S, \ord}_{\psi}(U^P)$ be a maximal ideal with residue field $k$. There exists a continuous semi-simple representation $\rhobar^{\ord}_\m : G_{F, S} \rightarrow \GL_2(k)$
such that $\rhobar^{\ord}_\m$ is totally odd, and for any finite place $v \notin S$ of $F$, $\rhobar_\m(\Frob_v)$ has characteristic polynomial
$X^2 - T_v X + q_v S_v \in (\T^{S, \ord}_{\psi}(U^P) / \m) [X]$. If $\rhobar^{\ord}_\m$ is absolutely reducible, we say that the maximal ideal $\m$ is Eisenstein; otherwise, we say that $\m$ is non-Eisenstein.

Suppose that $\m$ is non-Eisenstein. For each $v \in S_p - P$, let $\lambda_v$ and $\tau_v$ be the Hodge type and inerital type given by $\sigma_v$. We define a global deformation problem
\begin{align*}
\cS^{P} = (&\rhobar^{\ord}_{\m}, F, S, \{ \O\llbracket \O_v^{\times}(p) \rrbracket \}_{v \in P} \cup \{\O\}_{v \in S - P}, \{\Dord_v\}_{v \in P} \cup \{\Dss_v\}_{v \in S_p - P} \cup \{\Dodd_v\}_{v \in S_{\infty}} \\
&\cup \{\DSt_v \}_{v \in \Sigma} \cup \{ \Dfz_{v_1} \}),
\end{align*}
where $\Dord_v$ is the ordinary deformation problem defined with respect to the character $\overline{\eta}_v$ given by $\overline{\eta}_v(\varpi_v) = U_{\varpi_v}$ mod $\m$ and $\overline{\eta}_v(\alpha) = \langle \alpha \rangle$ mod $\m$ for all $\alpha \in \O_{F_v}^{\times}$.

\begin{prop}
Suppose that $\m$ is non-Eisenstein. Then there exists a lifting of $\rhobar^{\ord}_\m$ to a continuous homomorphism 
\[
\rho^{\ord}_\m: G_{F, S} \rightarrow \GL_2(\T^{S, \ord}_{\psi}(U^P)_\m)
\] 
such that 
\begin{itemize}
\item for each place $v \notin S$ of $F$, $\rhobar^{\ord}_\m(\Frob_v)$ has characteristic polynomial $X^2 - T_v X + q_v S_v \in \T^{S, \ord}_{\psi}(U^P)_\m [X]$;
\item for each place $v \in P$, $\rhobar^{\ord}_\m |_{G_{F_{v}}} \sim \big(\begin{smallmatrix}
\chi_{v} & * \\ 0 & *
\end{smallmatrix} \big)$ such that $\chi_{v} \circ \Art_{F_v}(\varpi^{-1}_v) = \U_{\varpi_{v}}$ and $\chi_{v} \circ \Art_{F_v}(t)= \langle t \rangle$ for $t \in \O_{F_v}^{\times}$.
\end{itemize}
Moreover, $\rho^{\ord}_\m$ is of type $\cS^P$ and has determinant $\psi \varepsilon$.
\end{prop}

\begin{proof}
The proof of \cite[Proposition 2.4.4]{MR3252020} works verbatim in our setting.
\end{proof}

\subsection{Ordinary patching}
Let $\m$ be a non-Eisenstein maximal ideal of $\T^{S, \ord}_{ \psi}(U^P)$. Let $T = S - \{v_1\}$ and $(Q_N, \{\alpha_v\}_{v \in Q_N})$ be a Taylor-Wiles datum as in Lemma \ref{lemma:auxiliary}. There are isomorphisms $R^T_{\cS^P} \cong R_{\cS^P} \ctimes_\O \cT$ (resp. $R^{T, \psi}_{\cS^P} \cong R^{\psi}_{\cS^P} \ctimes_\O \cT$). Define $S_N = \O_N \ctimes_{\O} \Lambda_P$, $\Sinf = \Oinf \ctimes_{\O} \Lambda_P$. Denote $\Roinfp := A^T_{\cS^P} \llbracket x_1, \cdots, x_{g+t} \rrbracket$. Then $\Spf \Roinfp$ is equipped with a free action of $(\Gmc)^t$, and a $(\Gmc)^t$-equivariant morphism $\delta^{\Delta}: \Spf \Roinfp \rightarrow  (\Gmc)^t$, where $(\Gmc)^t$ acts on itself by the square of the identity map. Define $\Roinf$ by $\Spf \Roinf = (\delta^{\Delta})^{-1}(1)$ and $\Roinvinf$ by $\Spf \Roinvinf := \Spf \Roinfp / (\Gmc)^t$. We fix a $\Theta^*_{Q_N}$-equivariant surjective $A^T_{\cS^P}$-algebra homomoprhism $\Roinfp \twoheadrightarrow R^T_{\cS^P_{Q_N}}$ for each $N$, which induces a $\Theta^*_{Q_N}[2]$-equivariant surjective $A^T_{\cS^P}$-algebra map $\Rinf \twoheadrightarrow R^{T, \psi}_{\cS^P_{Q_N}}$.

Let $c \in \N$ and let $J$ be an open ideal in $\Sinf$. Let $I_J$ be the subset of $N$ such that $J$ contains the kernel of $\Sinf \rightarrow S_N$. For $N \in I_J$, define
\[
\Mord(c, J, N) := \Sinf / J \otimes_{S_N} \Sord(U_1(Q_N)(c, c), \O)_{\m_{Q_N, 1}}^d.
\]
Applying Taylor-Wiles method to $\Mord(c, J, N)$ by the same way as in Sect. \ref{section:patching} (with some choice of ultrafilter $\filter$), we obtain an $\Sinf$-module $\Moinf$, which is finite free over $\Sinf$ and endowed with a $\Sinf$-linear action of $\Roinf$. Moreover, we have $\Moinf / \a \Minf \cong \Mord(U^P)$ with $\a = \ker(\Oinf \rightarrow \O)$.

The following proposition is an analog of \cite[Theorem 4.3.1]{MR2941425} and \cite[Theorem 3]{MR3904451}.

\begin{prop} \label{prop:ordinarymod}
Assume that for each $v \in P$, the image of $\rhobar^{\ord}_{\m} |_{G_{F_v}}$ is either trivial or has order $p$, and that either $F_v$ contains a primitive fourth roots of unity or $[F_v : \Q_2] \geq 3$. We have $\Supp_{\Roinf} \Moinf = \Roinf$.
\end{prop}

\begin{proof}
Let $Q$ be a minimal prime ideal of $\Lambda_P$. Then $\Moinf / Q$ is a finite free $\Sinf / Q$-module. It follows that the depth of $\Moinf / Q$ as an $\Roinf$-module is at least $\dim \Sinf / Q$. Thus every minimal prime of $(\Roinf /Q) / \Ann(\Moinf/Q)$ has dimension at least $\dim \Sinf / Q$. On the other hand, by Proposition \ref{prop:orddeform}(2), $\Roinf / Q$ is irreducible of dimension 
\begin{align*}
&g + 1 + \sum_{v \in P} (3+2[F_v:\Qp]) + \sum_{v \in S_p - P} (3 + [F_v:\Qp]) + \sum_{v \in S_{\infty}} 2 + \sum_{v \in \Sigma} 3 \\
= &q + 4|T| + \sum_{v \in P} [F_v:\Qp]
\end{align*}
which is equal to $\dim \Sinf / Q$. Thus $\Moinf /Q$ is supported on all of $\Spec \Roinf / Q$ and the proposition follows.
\end{proof}

\begin{corollary} \label{cor:ordinarymod}
Under the assumption of Proposition \ref{prop:ordinarymod}, the homomorphism $\Rz_{\cS^P} \twoheadrightarrow \T^{S, \ord}_{\psi}(U^P)_{\m}$ induces isomorphisms
\[
(\Rz_{\cS^P})^{\red} \cong \T^{S, \ord}_{\psi}(U^P)_{\m}.
\]
\end{corollary}

\begin{proof}
Reducing modulo $\a$ we see that $\Sord(U^P)^d \cong \Moinf / \a$ is a nearly faithful $\Roinf / \a$-module. However, the action of $\Roinf / \a$ on $\Sord(U^P)$ factors through the homomorphism $\Roinf / \a \Roinf \twoheadrightarrow \Rz_{\cS^P} \twoheadrightarrow \T^{S, \ord}_{\psi}(U^P)_{\m}$. It follows that the induced map $(\Rz_{\cS^P})^{\red} \twoheadrightarrow \T^{S, \ord}_{\psi}(U^P)_{\m}$ is an isomorphism as required.
\end{proof}

\begin{corollary} \label{cor:ordidefin}
Under the assumption of Proposition \ref{prop:ordinarymod}, $\Rz_{\cS^P}$ is a finite $\Lambda_P$-module.
\end{corollary}

\begin{proof}
The proof of \cite[Corollary 8.7]{MR2979825} works verbatim in our setting. We include the proof for the sake of completeness. Corollary \ref{cor:ordinarymod} shows that $\Rz_{\cS^P}/J$ is a quotient of the finite $\Lambda_P$-module $\T^{S, \ord}_{\psi}(U^p)_{\m^{\ord}}$, for some nilpotent ideal $J$ of $\Rz_{\cS^P}$. This implies that $\Rz_{\cS^P} / \m'$ is a finite $k$-algebra, where $\m'$ is the maximal ideal of $\Lambda_P$. Thus the corollary follows from Nakayama's lemma.
\end{proof}

\subsection{Constructing Galois representations}
\begin{thm} \label{thm:ordglobalization}
Let $F$ be a totally real field and let
\[
\rhobar: G_F \rightarrow \GL_2(k)
\]
be a continuous representation unramified outside $p$. Suppose that $\bar{\rho}$ has non-solvable image.

Let $\Sigma$ be a finite subset of places of $F$ not containing those above $p$ and let $\Sigma_p = \Sigma \cup \{ v|p \}$. Given a subset $P$ of $\{v|p\}$ such that $\rhobar|_{G_{F_v}}$ is reducible, and an ordinary lift $\rho_v$ of $\rhobar|_{G_{F_v}}$ for each $v \in P$. 

Assume that there is a regular algebraic cuspidal automorphic representation $\pi$ of $\GL_2(\A_F)$ such that
\begin{itemize}
\item $\overline{\rho}_{\pi, \iota} \cong \rhobar$;
\item $\det \rho_{\pi, \iota}|_{G_{F_v}} = \det \rho_v$ for each $v \in P$;
\item $\pi_v$ is unramified outside $\Sigma_p$ and is special at $\Sigma$;
\item $\pi$ is $\iota$-ordinary at $v \in P$.
\end{itemize}

Then there is an automorphic lift $\rho: G_F \rightarrow \GL_2(\O)$ of $\rhobar$ such that
\begin{itemize}
\item $\rho$ is unramified outside $\Sigma_p$ and $\rho(I_v)$ is unipotent non-trivial at $v \in \Sigma$;
\item if $v \in S_p - P$, then $\rho|_{G_{F_v}}$ and $\rho_{\pi, \iota}|_{G_{F_v}}$ lies on the same irreducible component of the potentially semi-stable deformation ring given by $\rho_{\pi, \iota}|_{G_{F_v}}$;
\item if $v \in P$, then $\rho|_{G_{F_v}}$ and $\rho_v$ lies on the same irreducible component of the potentially semi-stable deformation ring (corresponding to $\rho_v$).
\end{itemize}
\end{thm}

\begin{proof}
This theorem is a variant of \cite[Theorem 10.2]{MR2979825}. Let $\psi = \varepsilon^{-1} \det \rho_{\pi, \iota}$. Choose a finite solvable totally real extension $F'$ of $F$ such that
\begin{itemize}
\item $[F' : \Q]$ is even;
\item $F'$ is linearly disjoint form $\overline{F}^{\Ker{\rhobar}}(\zeta_p)$;
\item $\rho_{\pi, \iota} |_{G_{F'}}$ is ramified at an even number of places outside $p$;
\item for every place $w$ of $F'$ lying above $P$, the image of $\rhobar |_{G_{F'_w}}$ is either trivial or has order $p$, and that either $F'_w$ contains a primitive fourth roots of unity or $[F'_w : \Q_p] \geq 3$. 
\end{itemize}

Let $D$ be the quaternion algebra with center $F'$ ramified exactly at all infinite places and all $w$ lying above $\Sigma$. Choose $w_1$ to be a place not in $\Sigma$ such that $v_1 \nmid 2Mp$ and $\Frob_{v_1}$ has distinct eigenvalues. Fix a place $v_1$ of $F$ dividing $w_1$. Let $S = S_p \cup S_{\infty} \cup \Sigma \cup \{v_1\}$ and $S' = S'_p \cup S'_{\infty} \cup \Sigma' \cup \{w_1\}$, where $S_p$ (resp. $S'_p$) is the set of places of $F$ (resp. $F'$) dividing $p$, $S_{\infty}$ (resp. $S'_{\infty}$) is the set of places of $F$ (resp. $F'$) above $\infty$, and $\Sigma'$ is the set of places of $F'$ lying above $\Sigma$. Denote $P'$ the set of places of $F'$ lying above $P$ and $U^{P'} = \prod_{w \notin P'} U_w$ the open compact subgroup of $G(\A^{\infty}_{F'})$ defined by $U_w = \O_D^{\times}$ if $w \notin P' \cup \{w_1\}$ and $U_{w_1}$ is the pro-$w_1$ Iwahori subgroup. Let $\sigma_v$ be the locally algebraic type given by $\rho_{\pi, \iota}$ if $v \in S'_p - P'$ and let $\m$ be the maximal ideal in $\T^{S',\ord}_{\psi}$ defined by $\pi|_{F'}$ and $\varpi$. Thus we are in the setting of previous sections.

Let $\lambda_v$ and $\tau_v$ be the type given by $\rho_v$ if $v \in P$ (resp. $\rho_{\pi, \iota}$ if $v \in S_p - P$) and let $\cC_v$ be an irreducible component of the potentially semi-stable deformation ring containing $\rho_v$ if $v \in P$ (resp. $\rho_{\pi, \iota}$ if $v \in S_p - P$). Define $\lambda_w, \tau_w$, $\cC_w$ similarly for $w \in S_p'$. Let $T = S - \{v_1\}$ and $T' = S' - \{w_1\}$. Let $\gamma$ be the character given by $\rho_{\pi, \iota} |_{G_{F_{v_1}}}$. Consider the following global deformation problems
\begin{align*}
\cR = &(\rhobar, S, \{\O\}_{v \in S} ,\{\D^{\cC_v}_v \}_{v \in S_p} \cup \{ \Dodd_v\}_{v \in \Sinf} \cup \{\DSt_v \}_{v \in \Sigma} \cup \{\Dur_{v_1} \}), \\
\cR' = &(\rhobar |_{G_{F'}}, S', \{\O\}_{w \in S'} ,\{\D^{\cC_w}_w \}_{w \in S_p'} \cup \{ \Dodd_w\}_{w \in \Sinf'} \cup \{\DSt_w \}_{w \in \Sigma'} \cup \{\Dur_{w_1} \}), \\
\cR^{P, \prime} = &(\rhobar|_{G_{F'}}, S', \{ \O \llbracket \O_{F'_w}^{\times}(p) \rrbracket \}_{w \in P'} \cup \{\O\}_{w \in S' - P'} ,\{\Dord_w \}_{w \in P} \cup \{\D^{\cC_w}_w \}_{w \in S_p' - P'} \cup \{ \Dodd_w\}_{w \in \Sinf'} \\ & \cup \{\DSt_w \}_{w \in \Sigma'} \cup \{\Dur_{w_1} \}).
\end{align*}
Then by Corollary \ref{cor:ordidefin}, $\Rz_{\cR^{P, \prime}}$ is a finite $\Lambda_{P'}$-module. Note that $\Rz_{\cR}$ is a quotient of $\Rz_{\cR^{P, \prime}} \otimes_{\Lambda_P} \O$ by Lemma \ref{lemma:ordchar}, thus a finite $\O$-module. Since the morphism $\Rz_{\cR^{\prime}} \rightarrow \Rz_{\cR}$ is finite by Proposition \ref{prop:finitedefbc} and $\Rz_{\cR^{\prime}}$ is a finite $\O$-module by Corollary \ref{cor:ordidefin}, we deduce that $\Rz_{\cR}$ is a finite $\O$-module. 

On the other hand, $\Rz_{\cR}$ has a $\bQp$-point since it has Krull dimension at least 1 by Proposition \ref{prop:Stypepresentation}. This gives the desired lifting $\rho$ of $\rhobar$. It remains to show that $\rho$ is automorphic, which follows from the automorphy of $\rho|_{G_{F'}}$ and solvable base change.
\end{proof}

\section{Main results} \label{section:main}
\begin{thm} \label{thm:main}
Suppose that $p$ splits completely in $F$ (i.e. $F_v \cong \Q_2$ for $v | p$). For each locally algebraic type $\sigma$, the support of $\Minf(\sigma^{\circ}) \otimes_{\Zp} \Qp$ meets every irreducible component of $\Rinf(\sigma)[1/p]$. 
\end{thm}

\begin{proof}
Given an arbitrary irreducible component $\cC$ of $\Rinf(\sigma)[1/p]$, we want to show that there is a point $y$ lying on $\cC$ such that $\Minf(\sigma^{\circ}) \otimes_{\Rinf(\sigma), y} E_y \neq 0$. 

For each $v|2$, let $\cC_v$ be the irreducible component of $\Rss_v$ given by $\cC$ and let $\cC'_v$ be the irreducible component of $R_v^{\lambda_v', \tau_v'}$ given by an automorphic lift of $\rhobar$ (which exists by assumption and $\cC'_v$ can be chosen to be ordinary of weight $(0, 0)^{\Hom(F, \bQp)}$ if $\rhobar_v$ is reducible).

Fix a place $\p$ of $F$ above $2$. We claim that the support of $\Minf(\sigma^{\circ}) \otimes_{\Zp} \Qp$ meets the irreducible component of $\Rinf(\sigma)[1/p]$ defined by $\cC_{\p}$ and $\cC_v'$ for $v \in S_p - \{\p \}$. In the case $\cC_\p$ is ordinary, this follows from Theorem \ref{thm:ordglobalization}, otherwise this is due to Theorem \ref{thm:absirredMinf}. Repeating the argument for each place $v | p$, we obtain a point lying on $\cC$. This proves the theorem.
\end{proof}

Due to the equivalent conditions in Theorem \ref{thm:BMequiv} and Lemma \ref{lemma:suitglob}, we obtain the following:

\begin{corollary} 
Conjecture \ref{conj:Breuil-Mezardm} and Conjecture \ref{conj:Breuil-Mezard} hold for each continuous representation $\rbar: G_{\Qp} \rightarrow \GL_2(k)$.
\end{corollary}

This gives a new proof of Breuil-M\'ezard conjecture when $p=2$, which is new in the case $\rbar \sim (\begin{smallmatrix} \chi & * \\ 0 & \chi \end{smallmatrix})$ with $\chi: \GQp \rightarrow k^{\times}$ a continuous character.

Another application of Theorem \ref{thm:main} is an improvement of a theorem in \cite{MR3544298} below, which is new in the case $\rhobar|_{G_{F_v}} \sim (\begin{smallmatrix} \chi & * \\ 0 & \chi \end{smallmatrix})$ for some $v | p$.

\begin{thm} \label{thm:FMC}
Let $F$ be a totally real field in which $p$ splits completely. Let $\rho : G_F \rightarrow \GL_2(\O)$ be a continuous representation. Suppose that
\begin{enumerate}
\item $\rho$ is ramified at only finitely many places;
\item $\bar{\rho}$ is modular;
\item $\bar{\rho}$ is totally odd;
\item $\bar{\rho}$ has non-solvable image;
\item for every $v | p$, $\rho \vert_{F_v}$ is potentially semi-stable with distinct Hodge-Tate weights.
\end{enumerate}
Then (up to twist) $\rho$ comes from a Hilbert modular form.
\end{thm}

\begin{proof}
Let $\psi = \varepsilon^{-1} \det \rho$. By solvable base change, it is enough to prove the assertion for the restriction of $\rho$ to $G_{F'}$, where $F'$ is a totally real solvable extension of $F$. Moreover, we can choose $F'$ satisfying
\begin{itemize}
\item $[F' : \Q]$ is even.
\item $F'$ is linearly disjoint form $\overline{F}^{\Ker{\rhobar}}(\zeta_p)$ and splits completely at $p$.
\item $\rhobar |_{G_{F'}}$ is unramified outside $p$.
\item If $\rho$ is ramified at $v \neq p$, then the image of inertia is unipotent.
\item $\rho$ is ramified at an even number of places outside $p$.
\end{itemize}

Let $\Sigma$ be the set of places outside $p$ such that $\rhobar |_{G_{F'}}$ is ramified. If $v \in \Sigma$, then
\[
\rhobar |_{G_{F'}} \cong \begin{pmatrix} \gamma_v(1) & * \\ 0 & \gamma_v \end{pmatrix},
\]
where $\gamma_v$ is an unramified character such that $\gamma_v^2 = \psi |_{G_{F'_v}}$.

Let $D$ be the quaternion algebra with center $F'$ ramified exactly at all infinite places and all $v \in \Sigma$. Choose a place $v_1$ of $F'$ as in the proof of Theorem \ref{thm:ordglobalization}. Let $S$ be the union of infinite places, places above $p$, $\Sigma$ and $v_1$. Let $U^p = \prod_{v \nmid p} = U_v$ be an open subgroup of $G(\A^{\infty, p}_{F'})$ such that $U_v = G(\O_{F'_v})$ if $v \neq v_1$ and $U_{v_1}$ is the pro-$v_1$ Iwahori subgroup. Let $\m$ be the maximal ideal in the Hecke algebra $\T^S_{\psi}(U^p)$ defined by $\rhobar |_{G_{F'}}$.  Thus we are in the setting of Sect. \ref{section: globalization}. 

By Theorem \ref{thm:main} and Lemma \ref{lemma:BMequiv} (3) with $\sigma$ the locally algebraic type associated to $\rho |_{G_{F'}}$, we see that $\rho |_{G_{F'}}$ is automorphic and this proves the theorem. 
\end{proof}

\medskip

\bibliographystyle{alpha}
\bibliography{Modularity}

\end{document}